\documentclass[a4paper,11pt]{article}

\usepackage{graphicx}
\usepackage{bm}

\usepackage{authblk}
\usepackage{amscd}
\usepackage{amsfonts}
\usepackage[numbers]{natbib}
\usepackage{amsthm}
\usepackage{amssymb}
\usepackage{amsmath}
\usepackage{algpseudocode}
\usepackage{algorithm}
\usepackage{cancel}
\usepackage{chngcntr}
\usepackage{caption}
\usepackage{enumerate}
\usepackage{epsfig}
\usepackage{epic}
\usepackage{eepic}
\usepackage{enumitem}
\usepackage{etoolbox}
\usepackage{fancyhdr}
\usepackage{fullpage}
\usepackage[a4paper, total={6in, 8in}]{geometry}
\usepackage{graphicx}
\usepackage{here}
\usepackage[utf8]{inputenc}
\usepackage{lscape}
\usepackage{listings}
\usepackage{longtable}
\usepackage{natbib}
\usepackage{nccmath}
\usepackage{parskip}
\usepackage{rotating}
\usepackage{subfigure}
\usepackage{tabularx}
\usepackage{textcomp}
\usepackage{tikz}
\usepackage{titlesec}
\usepackage{threeparttable}
\usepackage{varwidth}
\usepackage{framed,multirow}

\usepackage{xcolor}

\title{A finite-difference ghost-point multigrid method for multi-scale modelling of sorption kinetics of a surfactant past an oscillating bubble}

\author[1]{Clarissa Astuto}
\author[2]{Armando Coco}
\author[2]{Giovanni Russo}

\affil[1]{King Abdullah University of Science and Technology (KAUST), 4700, Thuwal, Saudi Arabia}
\affil[2]{Department of Mathematics and Computer Science, University of Catania, Viale Andrea Doria 6, 95125, Catania, Italy}

\begin{document}
\maketitle

\begin{abstract}
We propose a method for the numerical solution of a multiscale model describing sorption kinetics of a surfactant around an oscillating bubble. 
The evolution of the particles is governed by a convection-diffusion equation for the surfactant concentration $c$, with suitable boundary condition on the bubble surface, which models the action of the short range attractive-repulsive potential acting on them when they get sufficiently close to the surface \cite{multiscale_mod}. In the domain occupied by the fluid, the particles are  transported by the fluid motion generated by the bubble oscillations. 


The method adopted to solve the equation for $c$ is based on a finite-difference scheme on a uniform Cartesian grid and implemented in 2D and 3D axisymmetric domains. We use a level-set function to define the region occupied by the bubble, while the boundary conditions are discretized by a ghost-point technique to guarantee second order accuracy at the curved boundary. The sparse linear system is finally solved with a geometric multigrid technique designed \textit{ad-hoc\/} for this specific problem.
Several accuracy tests are provided to prove second order accuracy in space and time. 

The fluid dynamics generated by the oscillating bubble is governed by the Stokes equation solved with a second order accurate method based on a monolithic approach, where the momentum and continuity equations are solved simultaneously. 
Since the amplitude of the bubble oscillations are very small, a simplified model is presented where the computational bubble is actually steady and its oscillations are represented purely with time-dependent boundary conditions. A numerical comparison with the moving domain model confirms that this simplification is perfectly reasonable for the class of problems investigated in this paper.
\end{abstract}

\section{Introduction}
\begin{figure}[tb]
	\begin{minipage}
		{.48\textwidth}
		\centering
		\includegraphics[width=\textwidth]{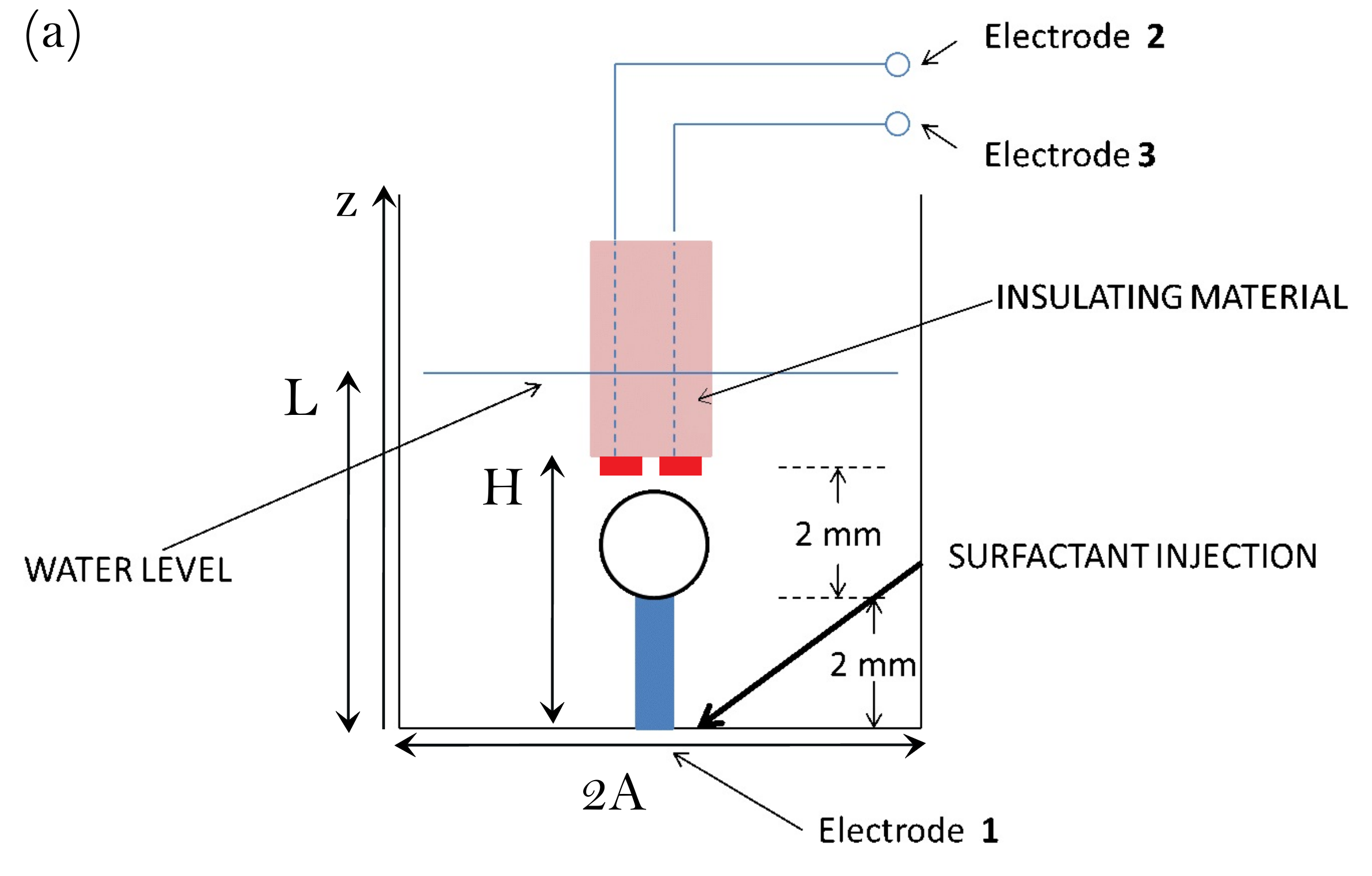}
	\end{minipage}
	\begin{minipage}
		{.48\textwidth}
		\centering
		\includegraphics[width=\textwidth]{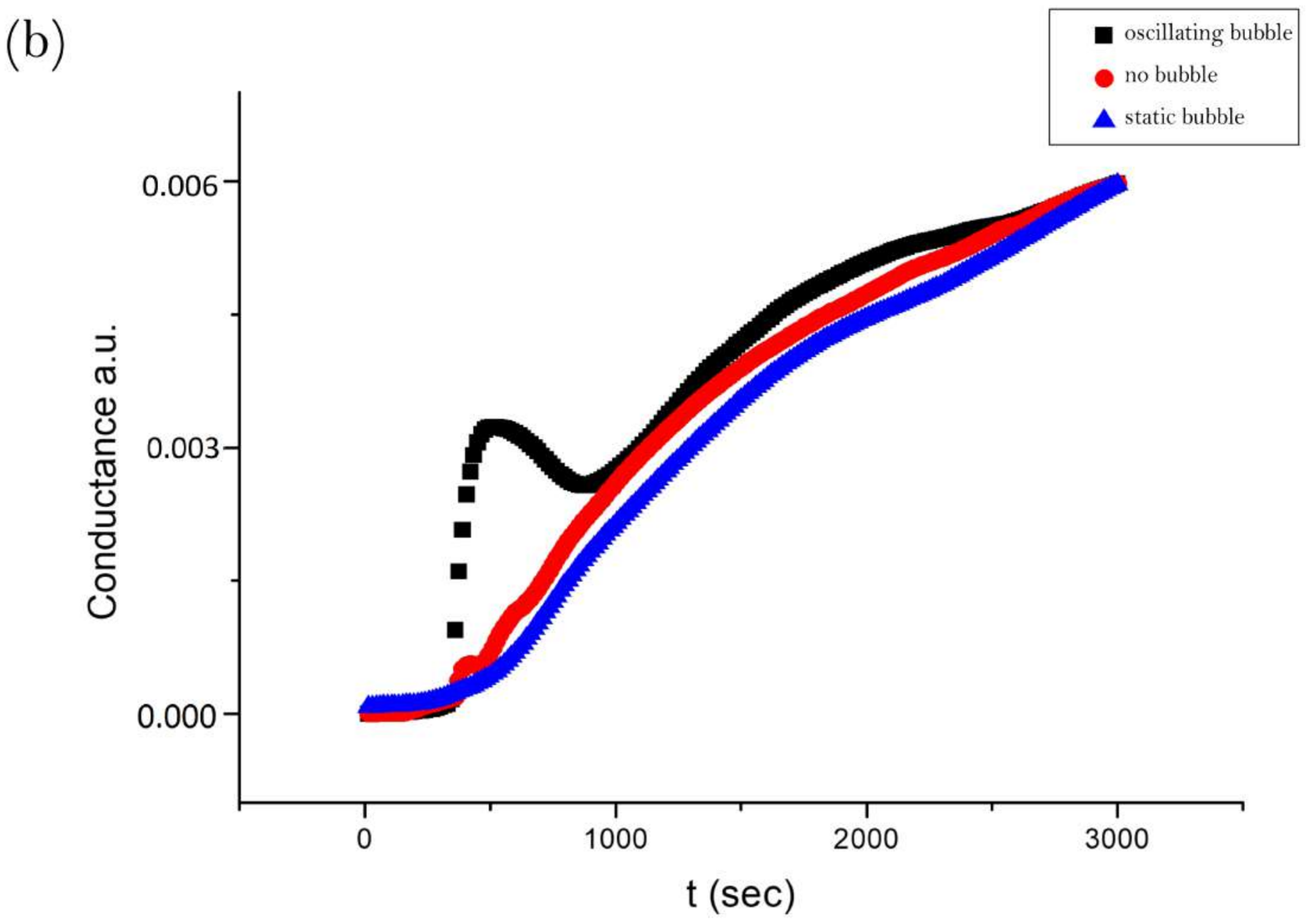}
	\end{minipage}
	\caption{\textit{Experimental domain and results from~\cite{Raudino20168574, CiCP-31-707}. (a) Schematic setup of the real apparatus. The central sphere represents the oscillating bubble. See~\cite{Raudino20168574, CiCP-31-707} for a detailed description and for experimental values of $H$, $A$, $L$. (b) Time evolution of the aqueous solution conductance measured above the bubble (electrodes 2 and 3 of (a)). Red line: no bubble; blue line: saturated bubble submitted to a flux of surfactants; black line: oscillating bubble submitted to a flux of surfactants (adapted from  \cite{Raudino20168574}). The line thicknesses represent an estimate of the experimental uncertainty of the conductivity measurements. }}
	\label{plot_setup}
\end{figure}

Diffusion equations and Stokes problems in a time dependent domain are present in countless research areas with many applications such as problems of temperature distribution \cite{doi10113716M1083207}, studies of biological pattern formation and cell motility on evolving surfaces \cite{doi:10.1098/rsif.2012.0276,PhysRevE.60.4588,garcke2013diffuse}. The most relevant topic that moves our research is the modelling of surfactants in two-phase flows using a diffuse interface, that is part of a long time project \cite{CiCP-31-707,Raudino20168574,variraudinoA,variraudinoB,variraudinoC,variraudinoD,variraudinoE,variraudinoF,variraudinoG}. In these works the authors report the experimental trapping kinetics of a diffusing flux of surfactants sticking at the surface of an oscillating gas bubble set in the middle of the vessel (see Fig.~\ref{plot_setup} (a)). The surfactant concentration past the oscillating bubbles is detected by conductivity measurements (see Fig.~\ref{plot_setup} (b)). A different and unexpected behavior is observed in presence of an empty bubble oscillating at resonance frequency (black curve in Fig.~\ref{plot_setup} (b)).  The phenomenon is particularly relevant when the bubbles are exposed to intense forced oscillations near resonance.


Surfactants are important for several industrial applications, such as processes of emulsification and mixing~\cite{wu2021new} or the production and stabilization of 2D nanomaterials~\cite{tyurnina2021environment}. They can be soluble in at least one of the fluid phases and the exchange of surfactants between the bulk phases and the fluid interfaces is governed by the process of adsorption and desorption. In \cite{garcke2013diffuse} different phase field models are derived for two-phase flow with a surfactant soluble in possibly both fluids. In \cite{doi:10.1063/1.1724167,Diamant_1996} they present models that account for both the diffusive transport inside the solution and the kinetics taking place at the interface using a free-energy formulation.

A theoretical work of Ward and Tordai \cite{ward1946time} 
formulated a time-dependent relation between the surface density of surfactants adsorbed at an interface and their concentration at the sub-surface layer of the solution, assuming a diffusive transport from the bulk solution.
Consequent theoretical works have focused on providing a second closure relation between
these two variables, as in \cite{multiscale_mod}, while in this paper we focus on the numerical aspects of the problem. 

Specifically, we study a diffusion equation in a bulk domain, with a dynamic time-dependent boundary condition derived by conservation arguments and stating that at the boundary the flux is proportional to the time derivative of the solution~\cite{PhysRevE.60.4588,Plaza}.
We formulate a finite-difference scheme for advection-diffusion equation with moving curved surfaces/boundaries. Time-discretization is performed with the Crank-Nicolson method. The bubble region is implicitly described by a level-set function, while the implementation of boundary conditions on complex-shaped boundaries/surfaces is based on a ghost-point method.

The ghost-point approach for domains described by level-set functions has been successfully proposed in several contexts~\cite{Fedkiw:GFM, Gibou:Ghost, Gibou:fourth_order, Gibou:fluid_solid, fernandez2020very, clain2021very}, as it has the advantage to allow an implicit representation of the boundary and then it can be employed on meshes that do not necessarily conform with a complex-shaped boundary. This advantage alleviates the computational burden that is associated with mesh generation steps, as observed in other approaches based on fitted-boundary methods.
A flexible ghost-point technique, suitable for different boundary condition types, is presented in~\cite{COCO2013464, COCO2018299} and applied to several contexts~\cite{COCO2020109623, chertock2018second}. In this paper, we extend the approach to accommodate time dependant boundary conditions and the presence of second order tangential derivatives.

A geometric multigrid method is employed to efficiently solve the sparse linear system arising from the discretization of the problem.
The multigrid approach is extended from~\cite{COCO2013464} in order to account for time-dependent boundary conditions and the presence of tangential derivatives. A suitable technique is presented in order to maintain the optimal efficiency of the multigrid method and contain the boundary effect degradation of the performance.
The method is second-order accurate in space and time, as confirmed by numerical tests.

The convection of particles is driven by the fluid motion around the oscillating bubble, usually modelled by incompressible Navier-Stokes equations. 
We assume the solute does not significantly change the density and rheology of the fluids, therefore the fluid motion is independent of the solute concentration (one way coupling).
We also assume that the motion of the bubble is assigned {\it a priori} and does not depend on the fluid dynamics. In a more realistic scenario, the bubble surface deformation would be influenced by the fluid motion, resulting in a two-way coupling as in fluid/membrane interaction problems~\cite{mavroyiakoumou2020large, mavroyiakoumou2021dynamics}.

For the physical parameters range adopted in the experiments, the Reynolds numbers are very small so that the convective terms of the Navier-Stokes equations can be neglected, and we can safely model the fluid-dynamics by the Stokes equations instead.
To numerically solve Stokes equations on moving domains we employ the method proposed in~\cite{COCO2020109623} based on a monolithic approach.

The paper is structured as follows. In Sect.~\ref{sect:mathmodel} we present the mathematical model for the diffusion of particles. Sect.~\ref{sect:FDdisc} describes the finite-difference ghost-point technique to implement the time-dependent boundary conditions.
In Sect.~\ref{sect:MG} we present the multigrid method to solve the sparse linear system arising from the ghost-point discretization. In Sect.~\ref{numtest} we perform several numerical tests to prove the second order accuracy in 2D and 3D axisymmetric geometries.
Sect.~\ref{sect:moving_bubble} deals with the fluid dynamics generated by the bubble oscillations. Small amplitude of the oscillations suggests that the bubble motion can be modelled purely from time-dependent boundary condition for the fluid velocity, while the computational bubble domain remains steady, saving then meaningful computational efforts.
This simplification is justified by numerical tests.
Finally, we couple the Stokes problem with the convection-diffusion equations to model the particle concentration evolution around an oscillating attracting bubble, proposing two types of oscillations (harmonic and ellipsoidal). 
Conclusions are drawn in Sect.~\ref{sec:conclusions}.


\section{Multiscale Model}\label{sect:mathmodel}
Modelling the diffusion in presence of a trap is challenging if multiple scales are involved. In recent papers, such as \cite{CiCP-31-707,Raudino20168574}, the range of the attractive-repulsive core of the trap is of the order of nanometers, a length that is 
several orders of magnitude smaller than the size of the domain. In order to overcome such difficulty, a \textit{multiscale model} was proposed in \cite{multiscale_mod}, which we briefly recall here, in the single carrier approximation.

\begin{figure}[!ht]
	\centering
	\begin{minipage}
		{.49\textwidth}
		\centering
		\includegraphics[width=\textwidth]{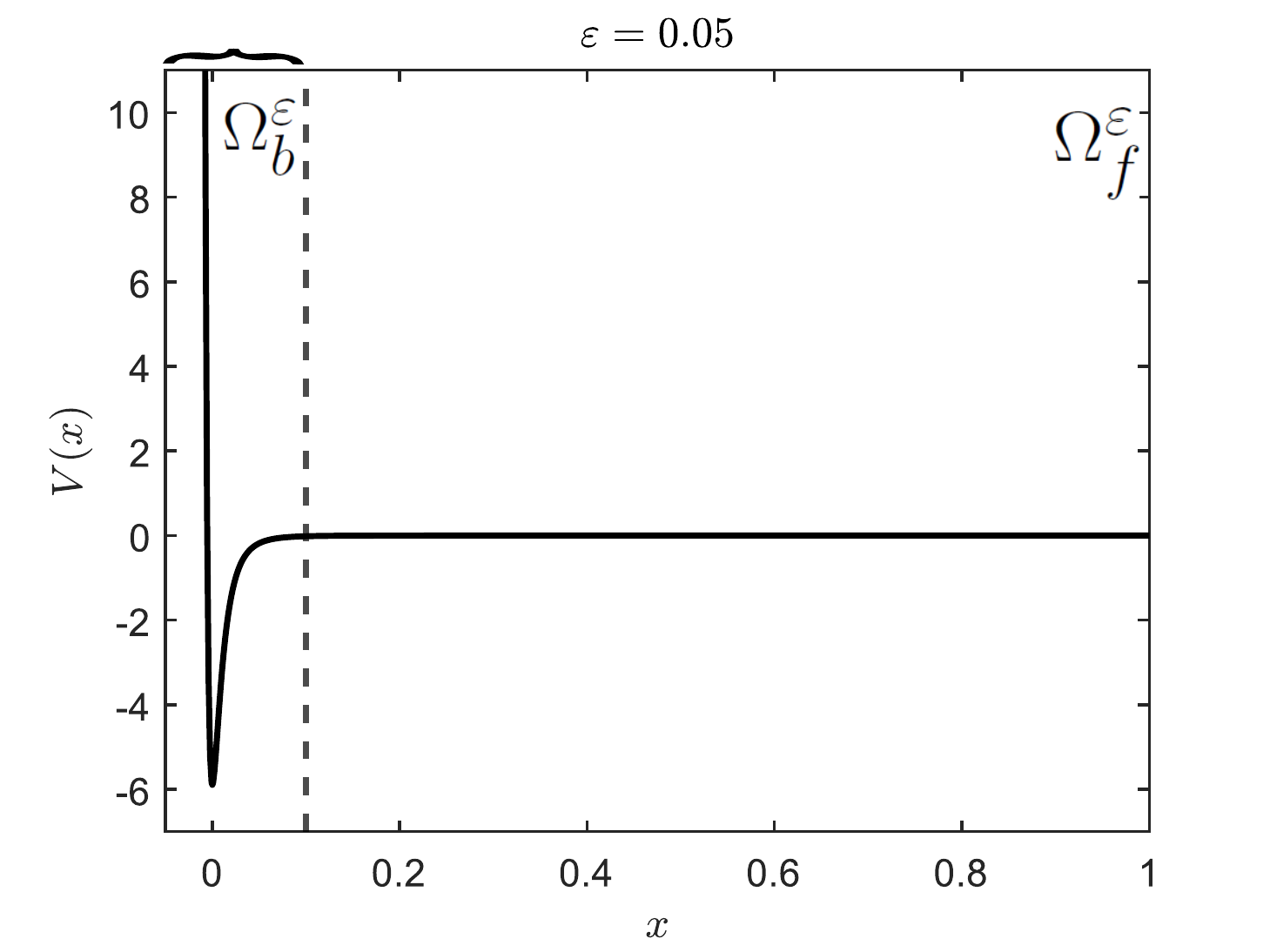}
	\end{minipage}
	\begin{minipage}
		{.49\textwidth}
		\includegraphics[width=\textwidth]{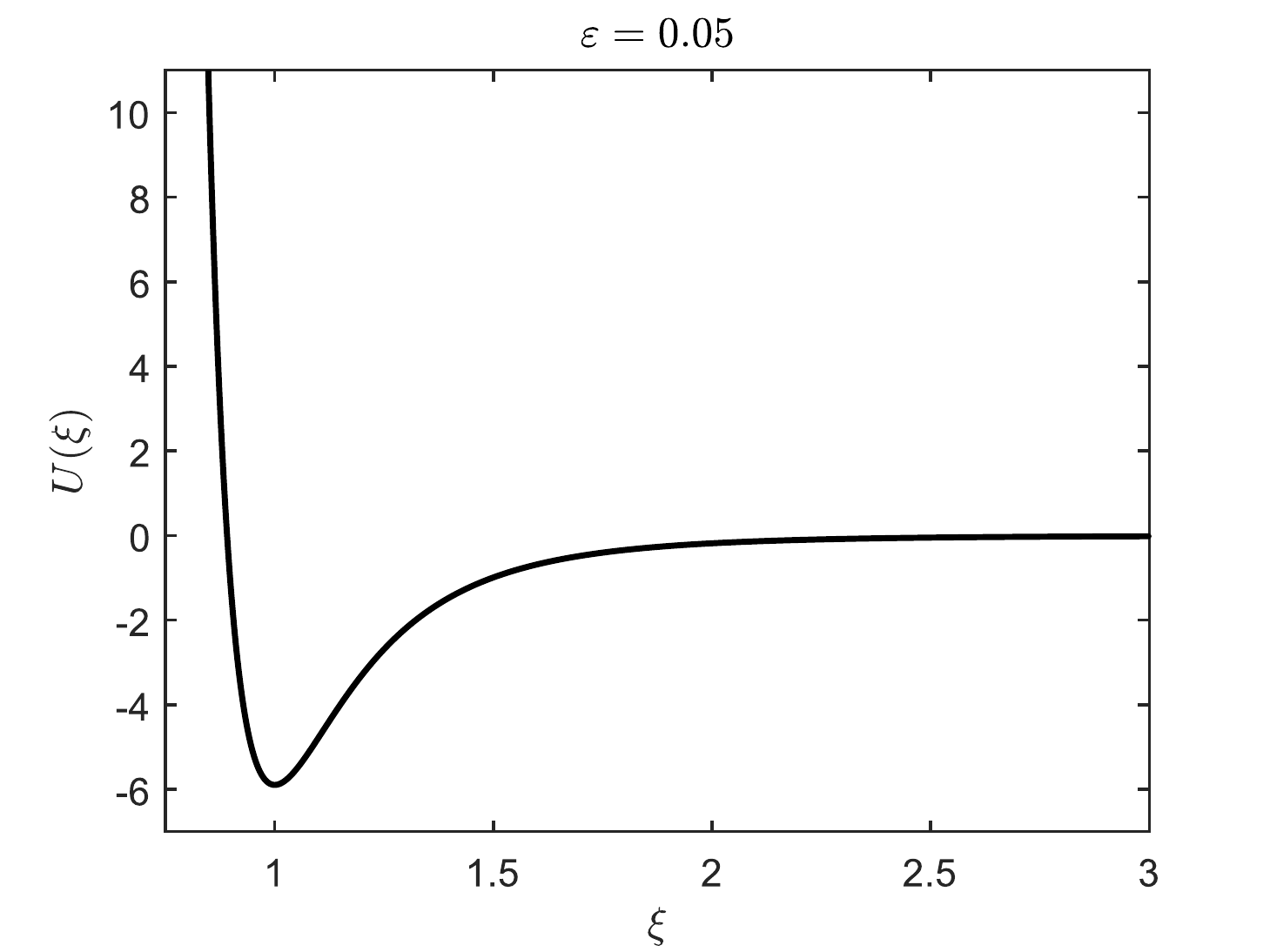}
	\end{minipage}	\caption{\textit{Representation of $V(x)$, on the left, and $U(\xi)$, on the right, for $\varepsilon = 0.05$ and $L=2$. On the left the dashed line $x=\varepsilon L$ denotes the right boundary of $\Omega_b^\varepsilon$. }}
	\label{figure_potential_V_1D}
\end{figure}

The time evolution of a local concentration of ions $c= c(\vec{x},t)$ diffusing in a steady fluid is governed by the conservation law \begin{equation}\label{eq:conservation}\frac{
\partial c}{\partial t}=-\nabla\cdot J.
\end{equation}
The flux term $J$ contains a diffusion and a drift term: 
\begin{equation} \label{eq:flux}
J=\ -D\left(\nabla c +\ \frac{1}{k_B T} c \nabla V\right)
\end{equation}
where $k_B$ is the Boltzmann's constant, $T$ is the absolute temperature and $V = V(\vec{x})$ is a suitable attractive-repulsive potential that models the particle trap. For simplicity, we describe the 1D model derivation only, referring the reader to~\cite{multiscale_mod} for a detailed derivation of higher order dimension models. In 1D, equations \eqref{eq:conservation} and \eqref{eq:flux} read:
\begin{align}
\label{eq_1D}
\displaystyle \frac{\partial  c}{\partial t} + \frac{\partial J}{\partial x} &=0 &\\
\label{eq_flux_1D}
\displaystyle J &= -D\left(\frac{\partial  c}{\partial x} + \frac{1}{k_B T} \,  c \, V' \right) &
\end{align}
Assume that the trap is located at $x=0$. Initially, particles are located within the region $x>0$.  If they get close to the trap, they are attracted towards $x=0$. This phenomenon is simulated by a potential $V(x)$ such that $V'(x)>0$ for $x>0$. On the other hand, if particles pass to the region $x<0$, they are repulsed towards $x=0$. Then, $V'(x)<0$ for $x<0$.
The attractive/repulsive mechanism is therefore modeled in the neighborhood of $x=0$ by a short range  potential $V(x)$ that is different from zero only in a thin region around the trap, say $\Omega^\varepsilon_b = [-\varepsilon,\varepsilon L]$ with $\varepsilon \ll 1$ and $L>0$. Therefore, $V(x)=0$ for $x\geq\varepsilon L$ and (assuming that the first derivative of $V(x)$ is continuous) $V'(\varepsilon L) = 0$.
A typical shape of the potential is reported in the left panel of Fig.~\ref{figure_potential_V_1D}.

Assuming that there is a wall at $x=1$, the fluid domain is represented by $\Omega_f^\varepsilon = [\varepsilon L , 1]$.
The problem consists of solving \eqref{eq:conservation} and \eqref{eq:flux} in $\Omega^\varepsilon =\Omega_b^\varepsilon \cup \Omega_f^\varepsilon = [-\varepsilon,1]$ with boundary conditions $J(-\varepsilon) = J(1)= 0$. This problem  presents a multiscale challenge because of the different spatial scales of $\Omega_b^\varepsilon = [-\varepsilon, \varepsilon L]$ and $\Omega_f^\varepsilon = [\varepsilon L, 1]$.
To overcome this difficulty, we aim at approximating the behaviour of the trap in $\Omega^\varepsilon_b$ with a suitable boundary condition at $x=0$, obtaining then a  simplified problem in $\Omega = [0,1]$ as follows.
Using a scaling variable $\xi = 1+x/\varepsilon$, the potential can be written in terms of $U(\xi)$ for $\xi \in [0,1+L]$ as $V(x) = U(\xi)$. In summary, we first choose a scaling potential $U(\xi)$ for $\xi \in [0,1+L]$ such that $U'(\xi)<0$ in $[0,1]$, $U'(\xi)>0$ in $[1,1+L]$, $U(1+L)=0$, $U'(1+L)=0$,
and then we study the behaviour of the trap when the potential is $V(x) = U(1+x/\varepsilon)$.
We assume that the solution $c_\varepsilon(\xi,t)$ of the scaled problem 
\begin{align}
\label{eq_eps}
\displaystyle \frac{\partial  c _\varepsilon}{\partial t} + \frac{1}{\varepsilon}\frac{\partial J_\varepsilon}{\partial \xi} &=0 &\\
\label{eq_flux_eps}
\displaystyle J_{\varepsilon} &= -D\frac{1}{\varepsilon}\left(\frac{\partial  c _\varepsilon}{\partial \xi} + \frac{1}{k_B T} \,  c _\varepsilon  \, U' \right) &
\end{align}
has the following expansion in $\Omega^\varepsilon_b$:
\begin{equation}\label{exprho}
 c_\varepsilon(\xi,t) =  c ^{(0)}(\xi,t)+\varepsilon  c ^{(1)}(\xi,t)+O(\varepsilon^2).
\end{equation} 
Since the flux $J_\varepsilon$ must be bounded for $\varepsilon \to 0$, from $\eqref{eq_flux_eps}$ we have that the coefficient of the term $\mathcal{O}(\varepsilon^{-1})$ in $J_\varepsilon$ has to vanish:
\begin{equation}
\frac{\partial  c^{(0)}}{\partial \xi}+ \frac{1}{k_B T} U'(\xi) c ^{(0)} = 0. 
\label{Jorder0}
\end{equation}
This equation can be solved for $c^{(0)}(\xi,t)$, yielding
\begin{equation}
 c ^{(0)}(\xi,t) =  c ^{(0)}(1+L,t)\exp \left(-\frac{U(\xi)}{k_B T}\right)
\end{equation}
since $U(1+L) = 0$.
Integrating \eqref{eq_1D} in $\Omega^\varepsilon_b$ we have:
\[
\frac{d}{dt}\int_{- \varepsilon}^{\varepsilon L} c (x,t) \, dx +J(\varepsilon L) - J(-\varepsilon) =  0
\]
and using the approximation $c(x,t) \approx c^{(0)}(\xi,t)$, the boundary condition $J(-\varepsilon)=0$ and that $V'(\varepsilon L)=0$ we obtain 
\begin{align}
	 \nonumber
 \varepsilon \;
 \frac{\partial c (\varepsilon L,t)}{\partial t} \;  \int_{0}^{1+L}\exp\left(-\frac{U(\xi)}{k_B T}\right) d \xi - D \frac{\partial  c (\varepsilon L,t)}{\partial x} &= 0 &
\end{align}
that represents a boundary condition of $c(x,t)$ at $x=\varepsilon L$. Using this boundary condition at $x=0$ instead of $x=\varepsilon L$, we finally obtain the following simplified problem related to the multiscale model:
\begin{align}
\frac{\partial  c }{\partial t} &= \frac{\partial}{\partial x} \left( D \frac{\partial c}{\partial x}  \right) \quad {\rm {for} }\, x \in [0,1] & \label{reduced1d} \\
\frac{\partial  c }{\partial x} &= 0  \quad {\rm {at} }\, x = 1 & \label{BCt} \\
M\frac{\partial  c }{\partial t} - D\frac{\partial  c }{\partial x} &= 0 \quad {\rm {at} }\, x = 0 & \label{BCeps}
\end{align}
where
\begin{equation}
\label{expr_M}
M=\varepsilon\int_{0}^{1+L}\exp\left(-\frac{U(\xi)}{k_B T}\right)d\xi.
\end{equation}
We observe that if the potential does not depend on $\varepsilon$, $M \rightarrow 0$ as $\varepsilon \rightarrow 0$ and then the condition \eqref{BCt} reduces to a zero Neumann boundary condition, therefore the interesting multiscale limit is obtained by letting $\varepsilon\to 0$, still maintaining $M$ finite.\footnote{The effective dependence of the small still finite size of $\varepsilon$ is studied in \cite{multiscale_mod}.}

\begin{figure}[htp]
\centering
\hfill
\begin{minipage}[b]
	{.3\textwidth}
	\centering
	\includegraphics[width=1.1\textwidth]{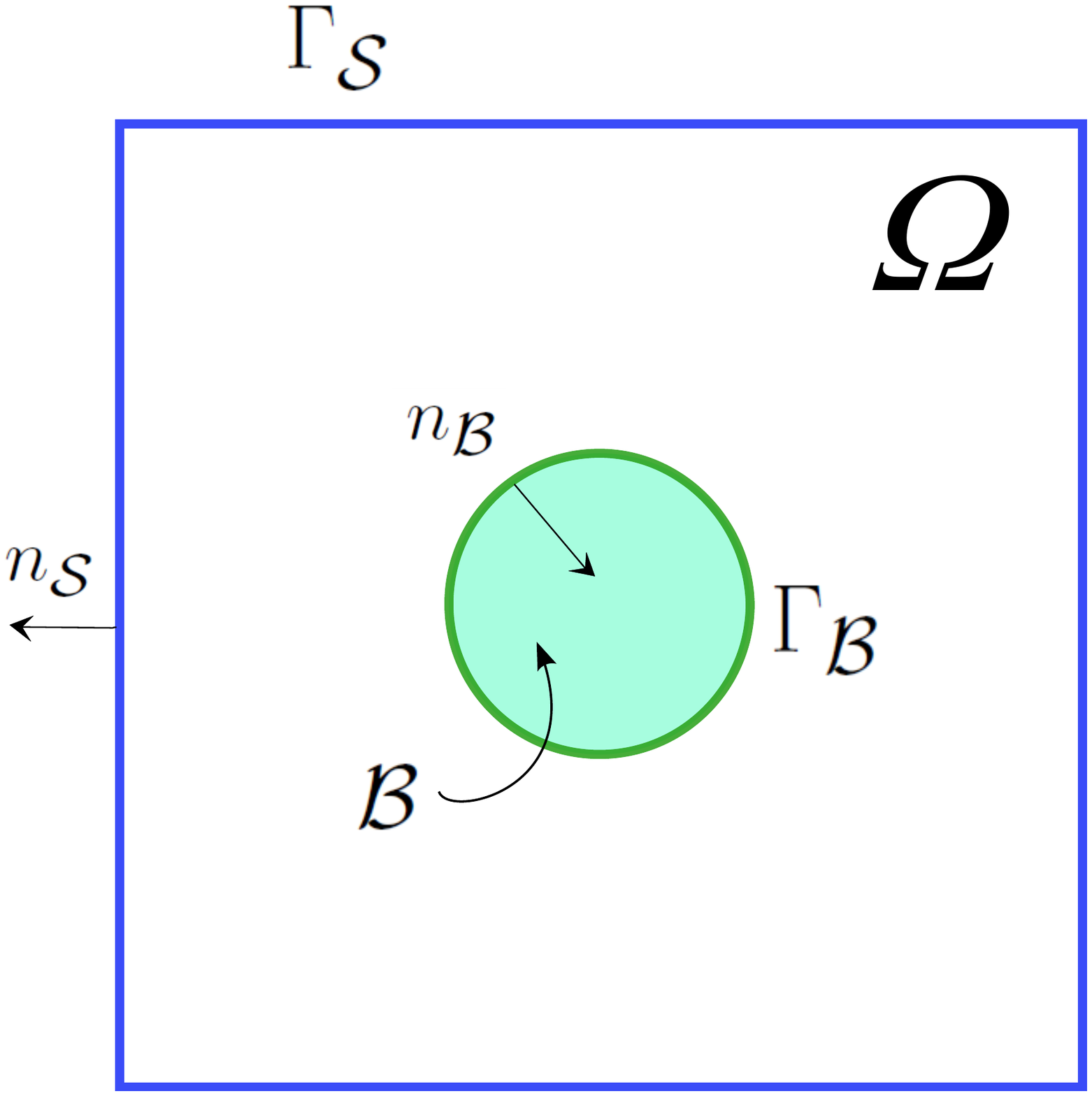}
\end{minipage}\hfill
\begin{minipage}[b]
	{.45\textwidth}
	\centering
	\includegraphics[width=0.6\textwidth]{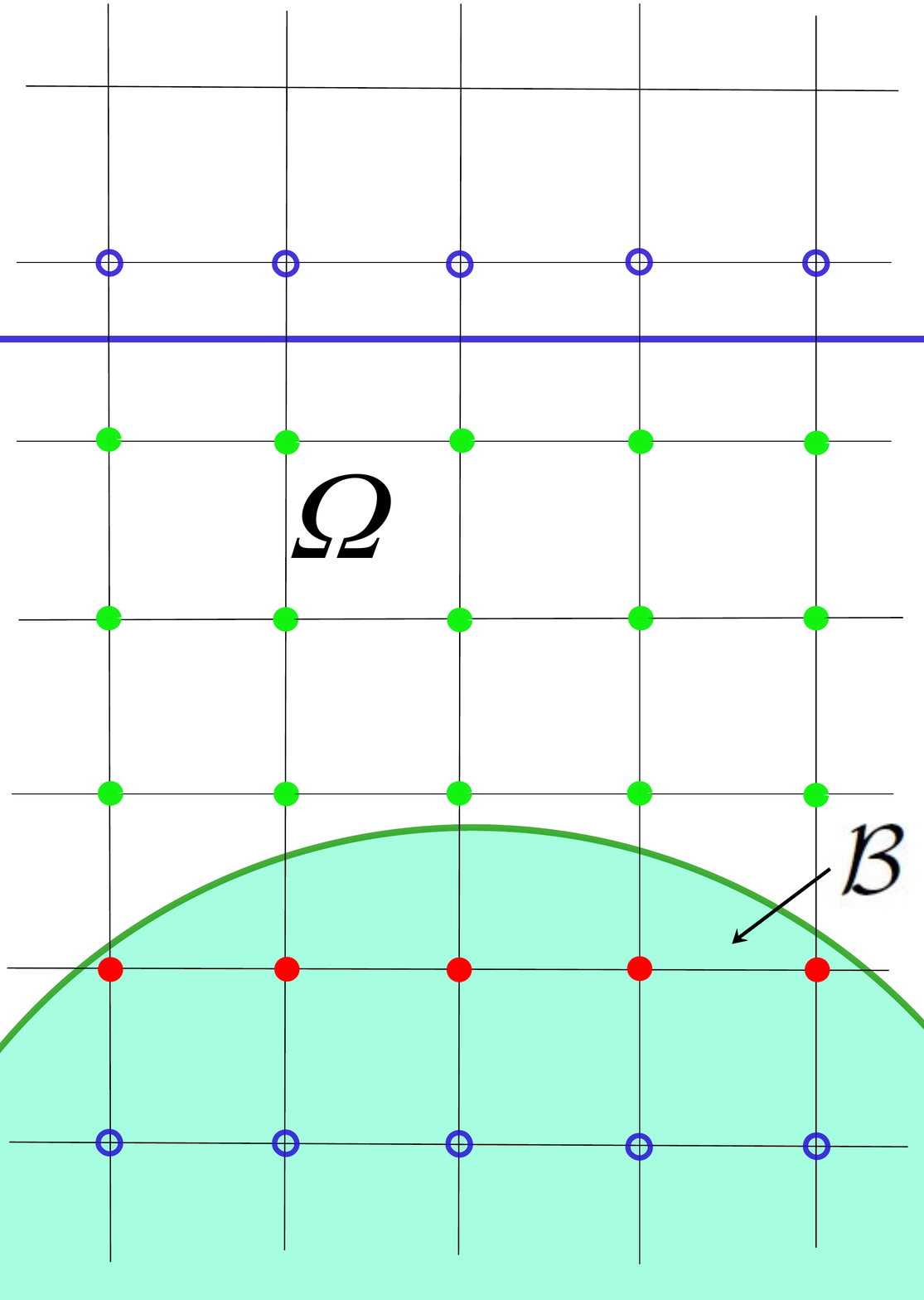}
\end{minipage}
\hspace*{\fill}
\caption{\textit{Representation of the domain on the left and classification of inside grid points (green), ghost points (red) and inactive points (blue circles) on the right.}}
\label{classification_points}
\end{figure}

Let us now describe the problem in higher dimensions. 
In two space dimensions, the fluid is contained in a domain $\Omega$ which is the a region external to a bubble $\mathcal{B}$ and internal to the square box $\mathcal{S} = (-a,a)^2 \subset \mathbb{R}^2$, with $a>0$ (see Figure 
\eqref{classification_points}, left panel).
Eq.~\eqref{reduced1d} reads:
\begin{equation}\label{pde2d}
\frac{\partial   c }{\displaystyle \partial t} = \nabla \cdot \left( D \nabla  c  \right) \text{ in } \Omega
\end{equation}
where $D$ is the diffusion coefficient.
Imposing zero flux at the wall results in homogeneous Neumann boundary conditions on $\Gamma_\mathcal{S} = \partial \mathcal{S}$:
\begin{equation}\label{bcwall}
\frac{\partial  c }{\partial {n}_\mathcal{S}} = 0 \quad \text{ on } \Gamma_\mathcal{S},
\end{equation}
where $n_\mathcal{S}$ denotes the unit normal vector on $\Gamma_\mathcal{S}$, pointing out of the domain $\Omega$.

In the presence of a steady bubble, the fluid domain is represented by $\Omega = \mathcal{S} \backslash \mathcal{B}$, where $\mathcal{B}$ is the region occupied by the bubble and represented by a sphere centred in the origin and with radius $R_\mathcal{B}$ such that $0 < R_\mathcal{B} < a$.
Similarly to \eqref{BCeps}, a suitable boundary condition is enforced on the boundary $\Gamma_\mathcal{B} = \partial \mathcal{B}$ to simulate the \textit{attractive-repulsive} mechanism of the bubble surface with the particles. In 2D, the analogue of boundary condition \eqref{BCeps} becomes (see~\cite{multiscale_mod} for more details): 
\begin{equation}\label{bcbubble}
M\frac{\partial  c }{\partial t} = MD\frac{\partial ^2 c }{\partial {\tau} ^2}-D\frac{\partial  c }{\partial {n}_\mathcal{B}}\quad \text{ on } \Gamma_\mathcal{B},
\end{equation}
where $M$ is given by the analogue of \eqref{expr_M} (the integration of the potential is performed along the direction normal to the surface of the bubble), $\tau$ denotes the unit vector tangential to $\Gamma_\mathcal{B}$, 
$n_\mathcal{B}$ is the unit normal vector on $\Gamma_\mathcal{B}$ pointing out of the domain $\Omega$ and $\partial^k/\partial\tau^k$ denotes the $k$-th derivative along such tangential direction.

In 3D the region $\mathcal{S}$ is the cube $(-a,a)^3$, the static bubble $\mathcal{B}$ is a sphere centered at the origin with radius $R_{\mathcal{B}}<a$, and the boundary condition on $\mathcal{B}$ becomes:
\begin{equation}\label{bcbubble3d}
M\frac{\partial  c }{\partial t} = MD\Delta_\perp  c -D\frac{\partial  c }{\partial n_\mathcal{B}}\quad \text{ in } \Gamma_{\mathcal{B}}
\end{equation}
where $\Delta_\perp$ is the Laplacian-Beltrami operator on the surface of the bubble $\Gamma_\mathcal{B}$.

\section{Finite-Difference discretization}\label{sect:FDdisc}
\subsection{Discretization in time}
Eqs.~\eqref{pde2d} and \eqref{bcbubble} can be written in compact form
\begin{equation}\label{compactQ}
\frac{\partial  c  }{\partial t} = Q \, c 
\end{equation}
where $Q$ is the following (linear) differential operator
\begin{equation}
Q \, c =
\left\{
\begin{matrix}
	D \Delta  c  & \text{ in } \Omega\\
	\\
	\displaystyle D\frac{\partial ^2 c  }{\partial \tau ^2}-DM^{-1} \frac{\partial  c  }{\partial n} \quad & \text{ on } \Gamma_\mathcal{B}
\end{matrix}
\right.
\label{Qexpr}
\end{equation}
with homogeneous Neumann boundary condition \eqref{bcwall} on $\Gamma_\mathcal{S}$.

Eq.~\eqref{compactQ} is discretized in time by using the \textit{Crank-Nicolson} method, which is second order accurate:
\begin{align} \nonumber
\frac{ c  ^{n+1}- c  ^n}{k} &= \frac{1}{2} \left( Q \, c ^n + Q \, c  ^{n+1}  \right) &\\  \label{CNdisc}
\left(I - \frac{k}{2}Q\right) \, c  ^{n+1} &= \left(I + \frac{k}{2}Q\right) \, c  ^n& 
\end{align}
where $k$ is the time step and $I$ is the identity operator.

\subsection{Discretization in space}\label{sect:discspace}

The computational domain $\mathcal{S}$ is discretized through a uniform Cartesian mesh with spatial step $h = 2a/N = \Delta x = \Delta y$, where $N^2$ is the number of cells.
We choose to use a cell-centered discretization to facilitate the implementation of homogeneous Neumann boundary conditions on $\Gamma_\mathcal{S}$. However, the accuracy of the method does not rely on this choice and a vertex-centered discretization would produce similar results.
Therefore, the set of grid points is $\mathcal{S}_h = \{(x_i,y_j)=(-a - h/2 + ih,-a - h/2 + jh), (i,j) \in \{1,\cdots,N\}^2 \}$. Within the set of grid points we define the set of internal points $\Omega_h = \mathcal{S}_h \cap \Omega$, the set of bubble points $\mathcal{B}_h =\mathcal{S}_h \cap \mathcal{B}$ and the set of ghost points $\mathcal{G}_h$
as grid points inside the bubble with at least one neighbor point inside $\Omega$:
\begin{equation}
(x_i,y_j) \in \mathcal{G}_h \iff (x_i,y_j) \in \mathcal{B}_h \text{ and } \{(x_i \pm h,y_j),(x_i,y_j\pm h) \} \cap \Omega_h \neq \emptyset.
\end{equation}
The remaining grid points that are neither inside nor ghost are called inactive points. See Fig.~\ref{classification_points} (right panel) for a classification of inside, ghost, and inactive points.
Let $N_I = |\Omega_h|$ and $N_G = |\mathcal{G}_h|$ be the number of internal and ghost points, respectively. We aim at approximating the solution $c$ at grid points of $\Omega_h \cup \mathcal{G}_h$, then our numerical solution can be represented as a column vector $ c _h = (\ldots,  c _{i,j}, \ldots)^T \in \mathbb{R}^{N_I+N_G}$, after choosing a bijective map between $\left\{ 1,\ldots,N_I+N_G \right\}$ and the grid points of $\Omega_h \cup \mathcal{G}_h$ (the overall numerical method does not rely on the particular choice of this map).
The problem \eqref{CNdisc} is then discretized in space, leading to a linear system
\begin{equation}\label{CNdiscspace}
\left(I_h - \frac{k}{2}Q_h\right)  c _h^{n+1} = \left(I_h + \frac{k}{2}Q_h\right)  c _h^n,
\end{equation}
to be solved at each time step, where $I_h$ and $Q_h$ are the $(N_I+N_G) \times (N_I+N_G)$ matrices representing the discretization of the operators $I$ and $Q$ and defined as follows. 
We denote by $I_h^{(i,j)}= \left(I_h^{(i,j),1},\ldots,I_h^{(i,j),N_I+N_G} \right)$ and $Q_h^{(i,j)}= \left(Q_h^{(i,j),1},\ldots,Q_h^{(i,j),N_I+N_G} \right)$ the rows of $I_h$ and $Q_h$, respectively, associated with the grid point $(x_i,y_j)$.

\begin{figure}[htp]
\centering
\hfill
\begin{minipage}[b]
	{.45\textwidth}
	\centering
	\includegraphics[width=0.75\textwidth]{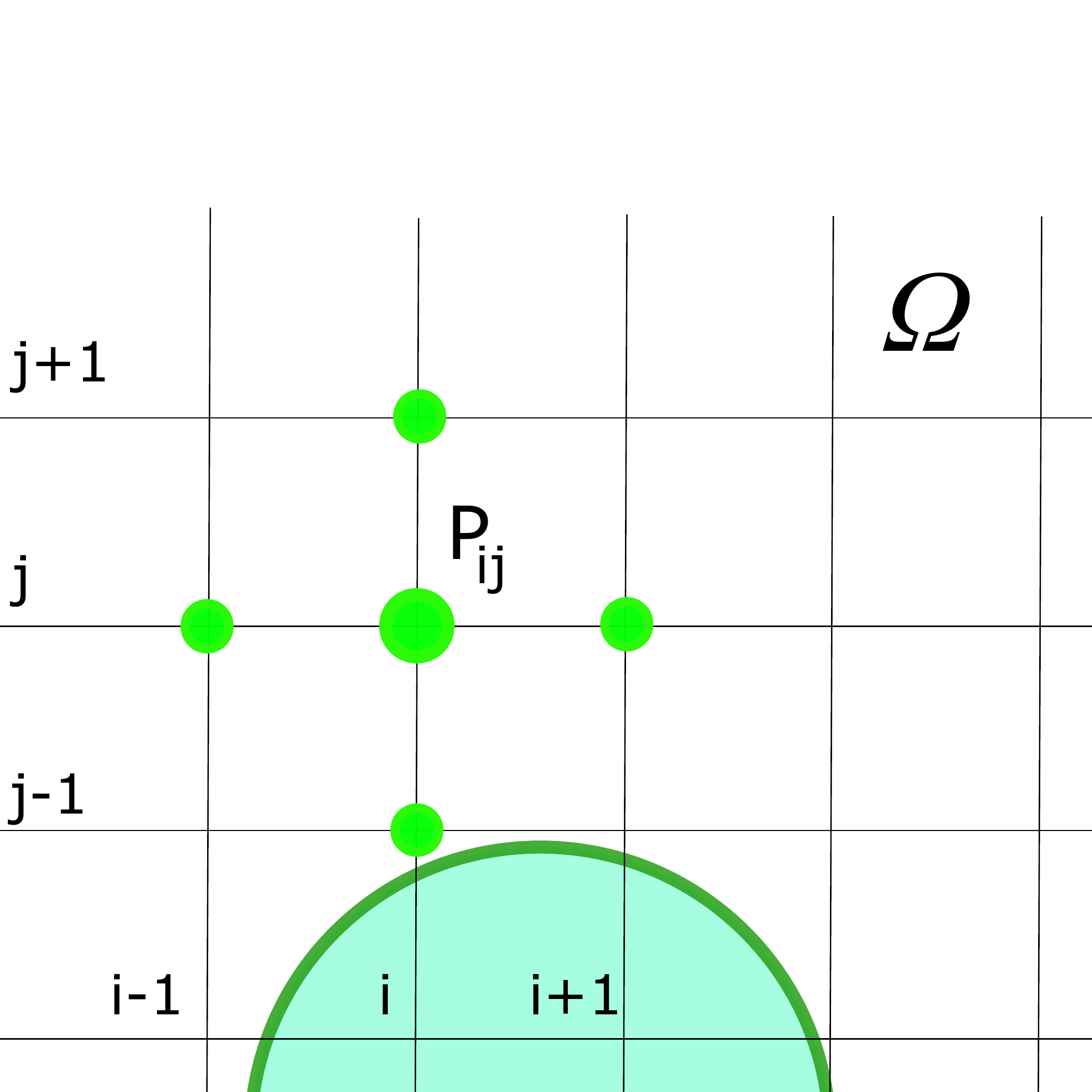}
\end{minipage}\hfill
\begin{minipage}[b]
	{.45\textwidth}
	\centering
	\includegraphics[width=0.75\textwidth]{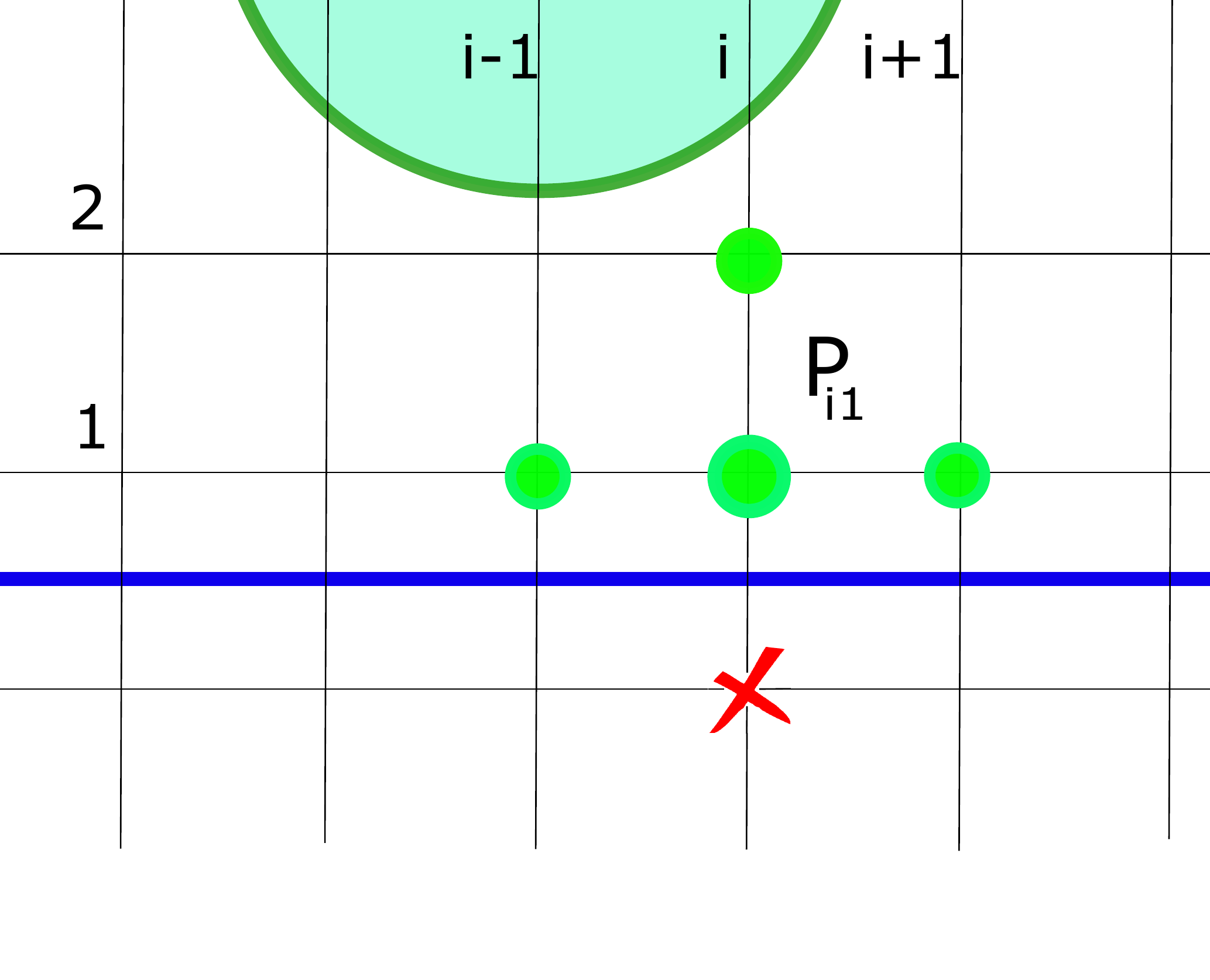}
\end{minipage}
\hspace*{\fill}
\caption{\textit{Representation of the five-point stencil for the discretization of internal points $P_{ij} = (x_i,y_j)$ (left panel) and the reduced stencil when $P_{ij}$ is close to the wall $\Gamma_{\mathcal{S}}$ (right panel). In the latter case the stencil is composed by four points.}}
\label{stencil0}
\end{figure}

If $P_{ij} = (x_i,y_j) \in \Omega_h$ in an internal grid point (as in Fig.~\ref{stencil0}, left panel), then the equation of the linear system is obtained from the discretization of the internal equation \eqref{pde2d} and the standard central difference on a five-point stencil is used to discretize the Laplace operator on $(x_i,y_j)$, and therefore $I_h^{(i,j)}$ and $Q_h^{(i,j)}$ are defined by
\begin{equation}\label{Ih_internal}
I_h^{(i,j)}  c _h =  c _{i,j}
\end{equation}
\[
Q_h^{(i,j)}  c _h = D \frac{ c _{i+1,j}+ c  _{i-1,j}+ c  _{i,j+1}+ c  _{i,j-1} - 4 c  _{i,j}}{h^2}.
\]
If $(x_i,y_j)$ is close the the wall and the five-point stencil contains grid points outside $\Omega$, we can use the boundary condition \eqref{bcwall} to reduce the five-point stencil and use only internal grid points. For example, looking at Fig.~\ref{stencil0} (right panel), we use the boundary condition
\[
\frac{ c _{i,0}- c _{i,1}}{h} = 0 \Longrightarrow  c _{i,0} =  c _{i,1}
\] 
and then 
\[
Q_h^{(i,1)}  c _h = D \frac{ c _{i+1,1}+ c  _{i-1,1}+ c  _{i,2}- 3 c  _{i,1}}{h^2}.
\]

\begin{figure}[htp]
\centering
	\centering
	\includegraphics[width=0.35\textwidth]{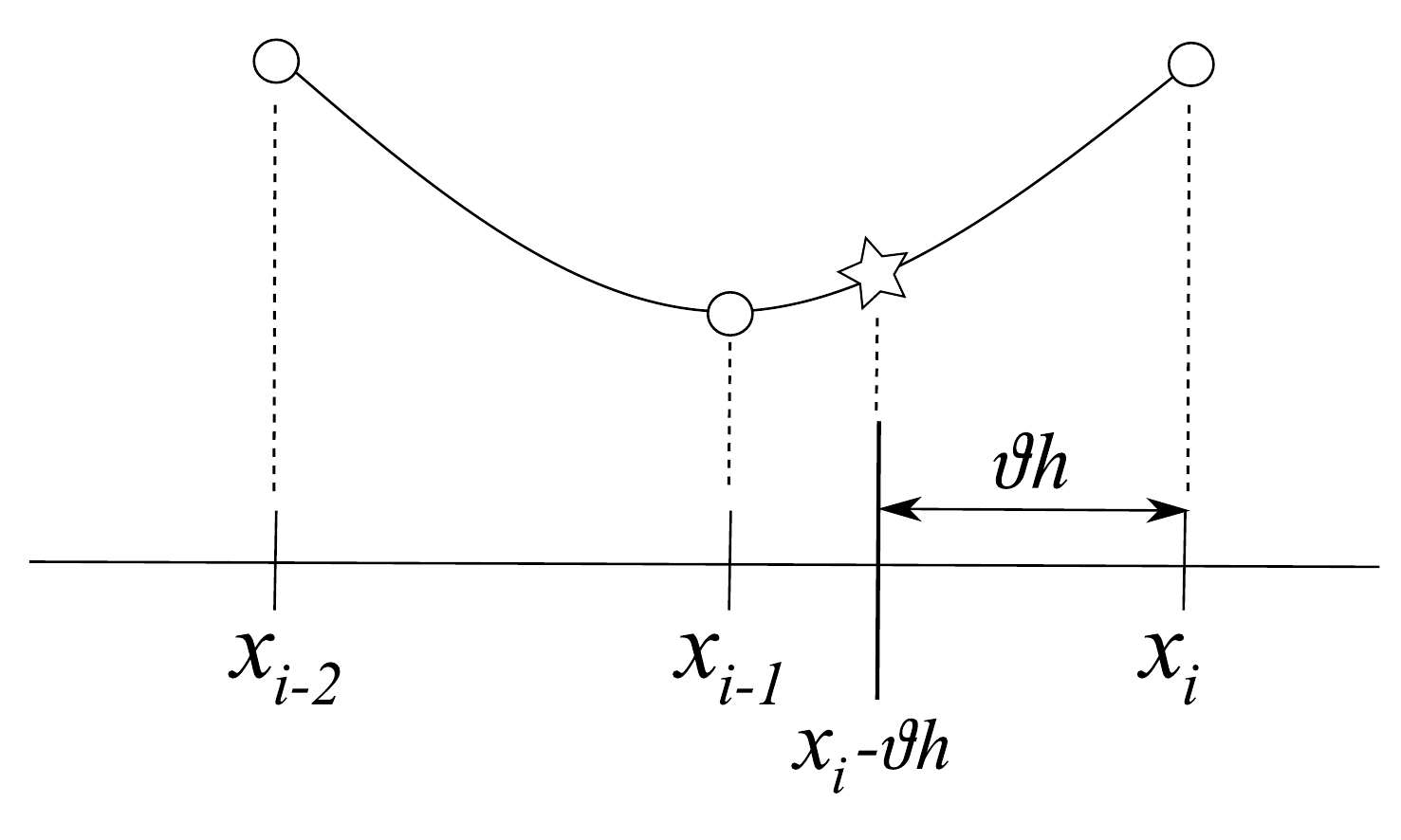}
\caption{\textit{1D interpolation on grid points $x_{i-2},x_{i-1},x_{i}$ (circle markers) to evaluate the function, its first derivative and the second derivative on $x_i - \vartheta h$ (star marker).
}}
\label{fig:1Dinterp}
\end{figure}

\begin{figure}[htp]
\centering
\hfill
\begin{minipage}[b]
	{.45\textwidth}
	\centering
	\includegraphics[width=0.8\textwidth]{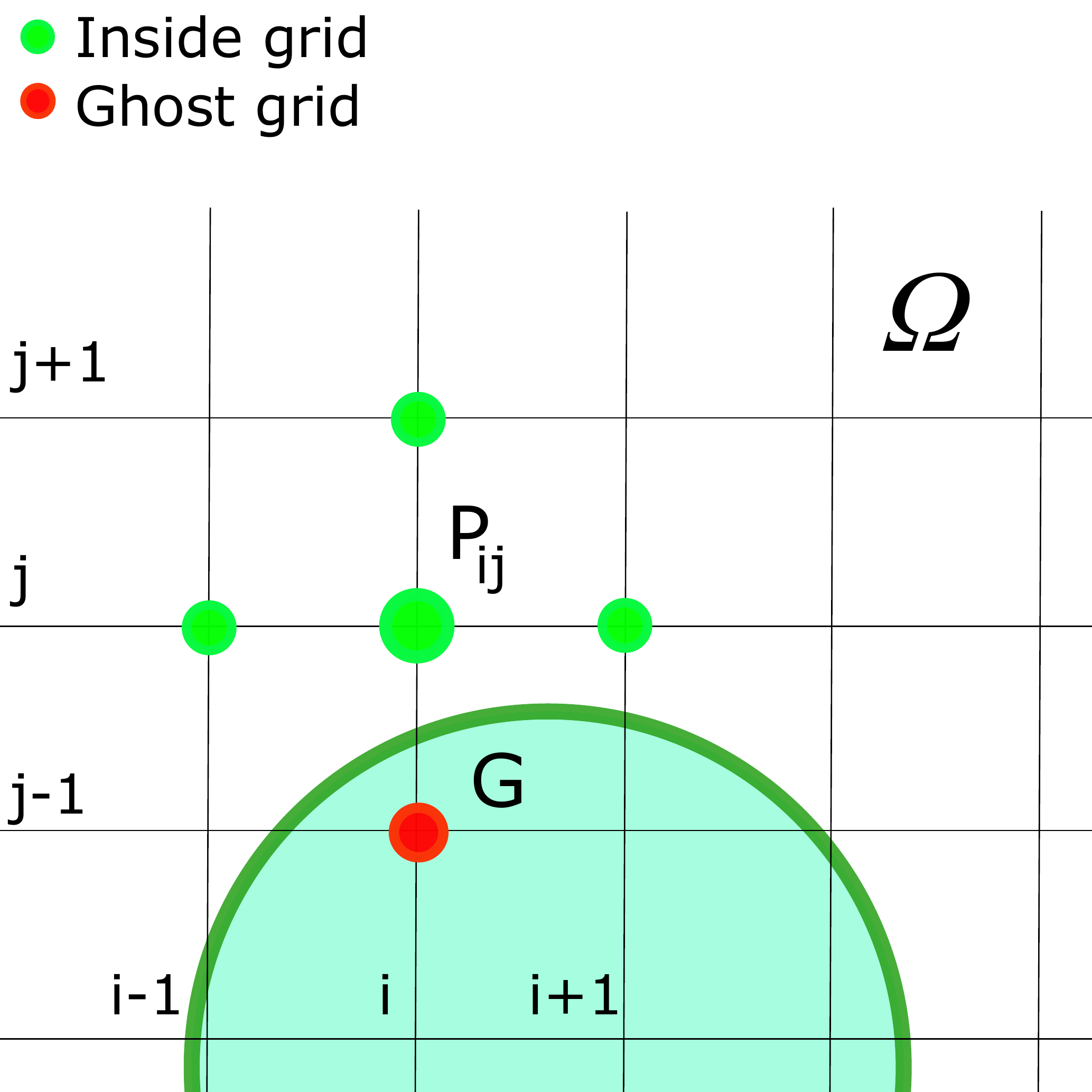}
\end{minipage}\hfill
\begin{minipage}[b]
	{.45\textwidth}
	\centering
	\includegraphics[width=0.6\textwidth]{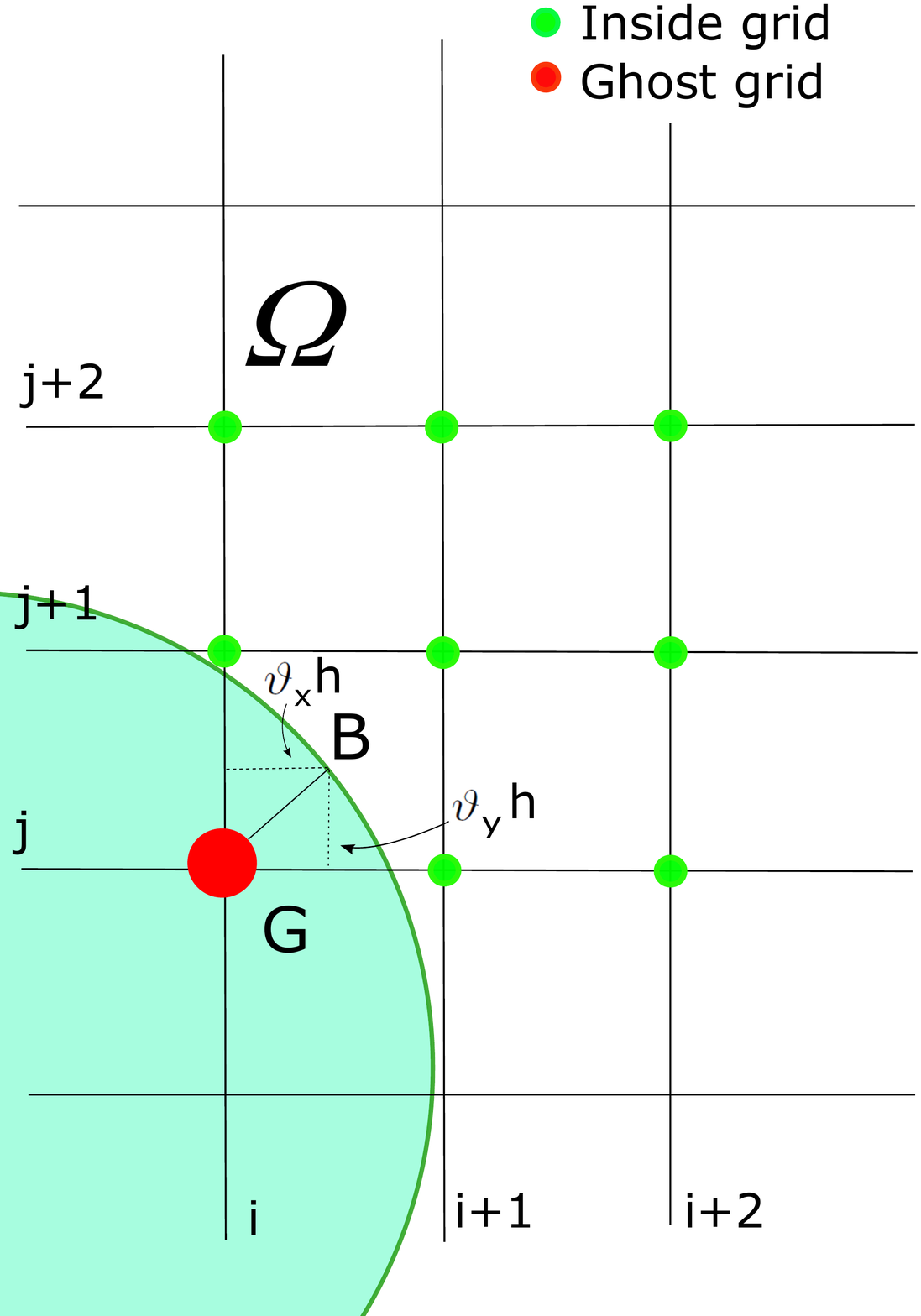}
\end{minipage}
\hspace*{\fill}
\caption{\textit{Representation of the five-point stencil when $P_{ij}$ is close to the boundary $\Gamma_{\mathcal{B}}$ (left panel). In this case the stencil contains a ghost point $G$. On the right panel we represent the upwind nine-point stencil associated with the ghost point $G$ and the boundary projection point $B$.}}
\label{stencil}
\end{figure}

If $G=(x_i,y_j) \in \mathcal{G}_h$ is a ghost point, then we discretize the boundary condition \eqref{bcbubble}, following a ghost-point approach similar to the one proposed in \cite{COCO2013464} and summarised as follows. We first compute the closest boundary point $B \in \Gamma_\mathcal{B}$ by
\[
B = O + R_\mathcal{B} \frac{G-O}{|G-O|},
\]
where $O$ is the centre of the bubble and $R_\mathcal{B}$ is the radius. Then, we identify the $3\times3$~--~point stencil having $G=(x_G,y_G)=(x_i,y_j)$ on one corner and whose convex hull contains $B=(x_B,y_B)$ (see Fig.~\ref{stencil}, right panel):
\[
\left\{ (x_{i+s_x m_x},x_{j+s_y m_y}) \colon m_x,m_y=0,1,2 \right\},
\]
where $s_x = {\rm SGN} (x_B-x_G)$ and $s_y = {\rm SGN} (y_B-y_G)$, with ${\rm SGN}(\alpha)=-1$ for $\alpha<0$ and ${\rm SGN}(\alpha)=1$ for $\alpha\geq0$. 
The solution $ c $ and its first and second derivatives are then interpolated at the boundary point $B$ using the discrete values $ c _{i,j}$ on the $3\times3$~--~point stencil.
The interpolations can be obtained as tensor products of 1D interpolations in the axis directions. 
In detail,
the 1D quadratic interpolations using the grid points $x_{i-2},x_{i-1},x_{i}$ to evaluate the function, its first derivative and the second derivative on $x_i - \vartheta h$
with $0\leq\vartheta<1$ (see Fig.~\ref{fig:1Dinterp}) are given by
\[
\tilde{ c }(x_i - \vartheta\,h) =  \sum_{m=0}^2 l_{m}(\vartheta) \,  c _{i-m},
\quad
\tilde{ c }'(x_i - \vartheta\,h) =  \sum_{m=0}^2 l'_{m}(\vartheta) \,  c _{i-m},
\quad
\tilde{ c }''(x_i - \vartheta\,h) =  \sum_{m=0}^2 l''_{m}(\vartheta)\, c _{i-m},
\]
where
\[
l(\vartheta) = \left( \frac{(1-\vartheta)(2-\vartheta)}{2}, \quad \vartheta (2-\vartheta), \quad \frac{\vartheta(\vartheta-1)}{2} \right)
\]
\[
l'(\vartheta) = \frac{1}{h} \left( \frac{(2\vartheta-3)}{2}, \quad 2(1-\vartheta), \quad \frac{(2\vartheta-1)}{2} \right)
\]
\[
l''(\vartheta) = \frac{1}{h^2} \left( 1, \quad -2, \quad 1 \right).
\]
In 2D, we define  (see Fig.~\ref{stencil}, right panel)
\[
\vartheta_x =  s_x (x_B-x_G)/h, \qquad
\vartheta_y =  s_y (y_B-y_G)/h.
\]
Observe that $0\leq \vartheta_x,\vartheta_y < 1$.
The 2D interpolation formulas are:
\begin{multline}
\label{coeffsLSstencil} 
\tilde{ c }(B) = \sum_{m_x,m_y=0}^2 l_{m_x}(\vartheta_x) l_{m_y}(\vartheta_y)  c _{i+s_x m_x,j+s_y m_y},
\\
\frac{\partial \tilde{ c }}{\partial x}(B) = s_x \sum_{m_x,m_y=0}^2 l'_{m_x}(\vartheta_x) l_{m_y}(\vartheta_y)  c _{i+s_x m_x,j+s_y m_y},
\\
\frac{\partial \tilde{ c }}{\partial y}(B) = s_y \sum_{m_x,m_y=0}^2 l_{m_x}(\vartheta_x) l'_{m_y}(\vartheta_y)  c _{i+s_x m_x,j+s_y m_y},
\\
\frac{\partial^2 \tilde{ c }}{\partial x^2}(B) = \sum_{m_x,m_y=0}^2 l''_{m_x}(\vartheta_x) l_{m_y}(\vartheta_y)  c _{i+s_x m_x,j+s_y m_y},
\\
\frac{\partial^2 \tilde{ c }}{\partial y^2}(B) = \sum_{m_x,m_y=0}^2 l_{m_x}(\vartheta_x) l''_{m_y}(\vartheta_y)  c _{i+s_x m_x,j+s_y m_y},
\\
\frac{\partial^2 \tilde{ c }}{\partial x \partial y}(B) = s_x\, s_y \sum_{m_x,m_y=0}^2 l'_{m_x}(\vartheta_x) l'_{m_y}(\vartheta_y)  c _{i+s_x m_x,j+s_y m_y}.
\end{multline}
Finally, the rows of $I_h$ and $Q_h$ associated with the ghost point $G=(x_G,y_G)$ are defined by evaluating the boundary condition on $B$, i.e.
\begin{equation}\label{IHghost}
I_h^{(i,j)}  c _h = \tilde{ c } (B)
\end{equation}
\begin{equation}\label{QHghost}
Q_h^{(i,j)}  c _h = D \left. \frac{\partial ^2 \tilde{ c }}{\partial \tau ^2} \right|_B - 
\frac{D}{M} \left. \frac{\partial \tilde{ c }}{\partial \bm{n}} \right|_B
\end{equation}
where
\begin{eqnarray}\label{normaleq}
\frac{\partial }{\partial n} = n_x\frac{\partial }{\partial x} +n_y\frac{\partial }{\partial y}, &\qquad
\displaystyle \frac{\partial ^2}{\partial \tau ^2} =
\displaystyle \tau_x^2\frac{\partial^2 }{\partial x^2} + 2\tau_x \tau_y\frac{\partial^2 }{\partial x \partial y} + \tau_y^2\frac{\partial^2 }{\partial y^2}, \\
(n_x,n_y) = \frac{O-G}{|O-G|}, & \qquad (\tau_x,\tau_y)= (-n_y,n_x).
\end{eqnarray}
Observe that, while the row $I_h^{(i,j)}$ associated with an internal point $(x_i,y_j)$ is the row of the identity matrix \eqref{Ih_internal}, for a ghost point $G=(x_i,y_j)$ this is not true (unless $G=B$, i.e.~$\vartheta_x=\vartheta_y=0$) since it contains $3^2=9$ values $l_{m_x}(\vartheta_x) l_{m_y}(\vartheta_y)$ for $m_x,m_y=0,1,2$.

\subsection{Complex-shaped bubbles: a level-set approach}
The discretization described in the previous sections for a spherical bubble can be extended to the case of more complex-shaped bubbles adopting a level-set approach. In detail, the bubble $\mathcal{B}$ can be implicit defined by a level set function $\phi(x,y)$ that is positive inside the bubble, negative outside and zero on the boundary $\Gamma_\mathcal{B}$ (\cite{Osher,book:72748}):
\begin{align}
\mathcal{B} &= \{(x,y): \phi(x,y) > 0\}&\\
\Gamma_\mathcal{B} &= \{(x,y): \phi(x,y) = 0\}.&
\end{align}
The unit normal vector $n$ in \eqref{normaleq} can be computed by:
\begin{equation}\label{LSnormal}
n = \frac{\nabla \phi }{|\nabla \phi|}
\end{equation}
provided that the level-set function is known explicitly. If it is known only at grid nodes, then the derivatives of $\phi$ in \eqref{LSnormal} are approximated by adopting a similar interpolation procedure as the one described in Eq.~\eqref{coeffsLSstencil}.

We observe that for a given bubble $\mathcal{B}$ there are infinite level-set functions. For example, $\phi=R_\mathcal{B}-\sqrt{x^2+y^2}$ and $\phi=R_\mathcal{B}^2-(x^2+y^2)$ describe the same circular  bubble. For a given bubble, the most convenient level-set function in terms of numerical stability is the signed distance function $\phi_d(x,y)$, i.e.~the distance between $(x,y)$ and $\Gamma_\mathcal{B}$ (positive inside the bubble, negative otherwise). The signed distance function $\phi_d$ can be computed from a generic level-set function $\phi$ by the reinitialization algorithm~\cite{sussman1994level, russo2000remark, du2008second}, 
consisting of finding the steady-state solution of:
\begin{equation}
\label{pde}
\frac{\partial \hat{\phi}}{\partial t} = {\rm sgn}(\phi)\left(1-|\nabla \hat{\phi} |\right), \qquad \hat{\phi}=\phi\quad \text{ at time } t=0
\end{equation}
where $t$ is a fictitious time.
A signed distance function is preferred to avoid numerical instabilities associated with sharp or shallow gradients close to the boundary (for a signed distance function we have $|\nabla \phi_d| =1$).
However, the cases investigated in the present paper involve steady bubbles or moving bubbles with a pre-determined evolution of the shape, then the instability issues of a generic level-set function are not observed. 



%


\section{Multigrid approach}\label{sect:MG}
The linear system \eqref{CNdiscspace} can be written as $A_h \, c _h^{n+1}=b_h$, where $A_h= \left(I_h - \frac{k}{2}Q_h\right)$ and $b_h=\left(I_h + \frac{k}{2}Q_h\right)  c _h^n$, and it is solved in this paper using an efficient multigrid approach that is an extension of the method proposed in \cite{COCO2013464} for elliptic equations on complex-shaped domains. In brief, a multigrid method is an iterative solver that starts by performing few steps of a suitable {\it relaxation scheme} to the linear system $A_h\, c _h^{n+1}=b_h$, obtaining an approximated solution $\bar{ c }^{n+1}$. The relaxation scheme is chosen in such a way that the high frequency Fourier modes of the residual $r_h = b_h-A_h\, c _h^{n+1}$ are dumped away much quicker than the low frequency Fourier modes. In other words, the relaxation operator smooths the residual $r_h$ after few relaxation steps (say $\nu_1$ steps). If so, it is said to have the {\it smoothing property}. Then, the residual $r_h$ is transferred to a coarser grid with spatial step $H=2h$ (without losing much information, as it is mainly composed of low frequency modes) by a suitable restriction operator $r_H = \mathcal{I}^h_H r_h$ and then the residual equation
$A_H e_H = r_H$ is solved on the coarse grid to obtain an approximation of the error $e_H$. Then, the error is transferred to the fine grid by an interpolation operator $e_h = \mathcal{I}^H_h e_H$ and the approximation $\bar{ c }^{n+1}$ is updated by $\bar{ c }^{n+1} \leftarrow \bar{ c }^{n+1} + e_h$. Few more steps (say $\nu_2$) of the relaxation operator are then performed on the fine grid to reduce the errors introduced by the interpolation procedure. The entire scheme is then performed iteratively until the residual falls below a certain tolerance. In addition, the residual equation $A_H e_H = r_H$ can be solved recursively by moving to a coarser grid with spatial step $2H$, and so on. Several types of multigrid schemes, such as $V-$cycle, $W-$cycle and Full Multigrid, can be adopted according to the different strategies that can be chosen.
In this paper we use a $W-$cycle approach and describe the main components of the multigrid method, i.e.~relaxation, restriction and interpolation operators, while we refer the reader to, for example,~\cite{Trottemberg:MG} for a comprehensive treatment of multigrid methods. 

\subsection{Relaxation scheme}
Standard relaxation schemes that show the smoothing property for elliptic equations in rectangular domains are the Gauss-Seidel scheme and the weighted Jacobi scheme (with weight $\omega = 2/3$ in 1D and $\omega = 4/5$ in 2D, see~\cite{Trottemberg:MG} ). While these relaxation schemes converge when the discretization is performed on a rectangular domain, they might not converge when using a ghost-point approach for curved boundaries (see \cite{COCO2013464}). To obtain a convergent scheme for the problem proposed in this paper, we modify the relaxation on ghost points as described below (while we keep a Gauss-Seidel scheme on the internal equations). The proposed relaxation scheme can be written in the Richardson form (with iterative index $k$)
\begin{equation}\label{rich}
c _h^{n+1,k+1} =  c _h^{n+1,k} + P_h^{-1} (b_h-A_h  c _h^{n+1,k})
\end{equation}
where $P_h$ is a $(N_I+N_G) \times (N_I+N_G)$ matrix called {\it preconditioner} and is chosen as a suitable approximation of $A_h$. A standard Gauss-Seidel scheme (on both internal and ghost points) corresponds to $P_h=(D_h+L_h)$, where $D_h$ and $L_h$ are the diagonal and lower part of $A_h$, respectively. 
To modify the scheme for ghost points, we change the $N_G$ diagonal values of $D_h$ that corresponds to ghost points, obtaining a new diagonal matrix $\tilde{D}_h$. Finally, the diagonal entries $\tilde{D}_h^{(i,j)}$ of $\tilde{D}_h$ are 
\[
\tilde{D}_h^{(i,j)} =
\left\{
\begin{matrix}
D_h^{(i,j)} = 1+\displaystyle\frac{2kD}{h^2} & \text{ if } & (x_i,y_j) \in \Omega_h \\
\beta & \text{ if } & (x_i,y_j) \in \text{ Ghost } \\
\end{matrix}
\right.
\]
where $\beta \in \mathbb{R}$ is a suitable value that we determine later.

We observe that in practice the relaxation \eqref{rich} is performed without storing the entire matrices $P_h$ and $A_h$ (matrix-free fashion). 
This has the advantage of avoiding the explicit construction of the sparse matrices with great simplification of implementation aspects and savings in computational time, especially for moving domains when the discrete operator depends on time. In fact, the vector $ c _h^{n+1}$ is computationally stored in a temporary array $ c _h$ and its components $ c _{i,j}$ are updated (overridden) in a Gauss-Seidel fashion by iterating all over the grid points.
We distinguish between internal and ghost points. We use the notation $a\leftarrow b$ to say that the variable $a$ is updated with the value $b$.

\paragraph{Internal points}
If $(i,j) \in \Omega_h$,
\[
c _{i,j}
\leftarrow
c _{i,j}+\frac{h^2}{h^2+2kD}
\left(b_{i,j}- c _{i,j}-\frac{kD}{2h^2} \left( c _{i+1,j}+ c  _{i-1,j}+ c  _{i,j+1}+ c  _{i,j-1}-4 c  _{i,j}\right) \right)
\]

\paragraph{Ghost points}
If $(i,j) \in \mathcal{G}_h$,
\begin{equation}\label{richghost}
c _{i,j}
\leftarrow
c _{i,j}+\beta^{-1}
\left(b_{i,j}-\left(I_h^{(i,j)}  c _h - \displaystyle \frac{k}{2} Q_h^{(i,j)}  c _h \right) \right)
\end{equation}
 where $I_h^{(i,j)}$ and $I_h^{(i,j)}$ are defined by \eqref{IHghost} and \eqref{QHghost}, respectively.
The iteration \eqref{richghost} can be written as
\begin{equation}\label{richghost2}
c _{i,j}
\leftarrow
\left(1-\beta^{-1} \left( I_h^{(i,j),(i,j)}- \displaystyle \frac{k}{2} Q_h^{(i,j),(i,j)} \right) \right)
c _{i,j}+ \ldots \text{ terms that do not depend on }  c _{i,j} \ldots
\end{equation}
and the value $\beta$ is chosen in such a way that the coefficient of $ c _{i,j}$ on the right-hand side of \eqref{richghost2} is not larger than one in absolute value, i.e.
\begin{equation}\label{condbeta}
\left| 1-\beta^{-1} A_h^{(i,j),(i,j)} \right|\leq 1,
\text{ with }
A_h^{(i,j),(i,j)} = \left( I_h^{(i,j),(i,j)}- \displaystyle \frac{k}{2} Q_h^{(i,j),(i,j)} \right).
\end{equation}
Using \eqref{IHghost} and \eqref{QHghost},
we have
\begin{equation}\label{diagI}
I_h^{(i,j),(i,j)} = l_{0}(\vartheta_x) \, l_{0}(\vartheta_y)
\end{equation}
%
\begin{align}\label{diagQ}
Q_h^{(i,j),(i,j)} = & D \left( \tau_x^2 l''_{0}(\vartheta_x) l_{0}(\vartheta_y) + 2 \tau_x \tau_y l'_{0}(\vartheta_x) l'_{0}(\vartheta_y) + \tau_y^2 l_{0}(\vartheta_x) l''_{0}(\vartheta_y) \right) \nonumber
\\& - \frac{D}{M} \left( n_x l'_{0}(\vartheta_x) l_{0}(\vartheta_y) + n_y l_{0}(\vartheta_x) l'_{0}(\vartheta_y) \right)
\end{align}
In order to satisfy condition \eqref{condbeta}, we require that $\beta$ is chosen in such a way that
\begin{equation}\label{condbeta1and2}
0 \leq \beta^{-1} A_h^{(i,j),(i,j)} \leq 2.
\end{equation}
The left inequality of \eqref{condbeta1and2} is satisfied by choosing $\text{sign } \beta = \text{sign } A_h^{(i,j),(i,j)}$ (we conventionally choose $\text{sign }(0) =1$). 
To satisfy the right inequality of \eqref{condbeta1and2}, we have
\[
\left| \beta \right| \geq \frac{\left| A_h^{(i,j),(i,j)} \right|}{2}.
\]
This condition is always satisfied (regardless of the ghost point) if we choose
\begin{equation}\label{condA}
\left| \beta \right| \geq \tilde{A}/2 \quad \text{ with } \quad \sup_{\vartheta_x,\vartheta_y \in [0,1]} \left| A_h^{(i,j),(i,j)} \right| \leq \tilde{A}.
\end{equation}
The estimate $\tilde{A}$ can be found as follows.
Using \eqref{diagI}, \eqref{diagQ} and \eqref{coeffsLSstencil},
we have
\begin{multline}\label{Abound}
\left| A_h^{(i,j),(i,j)} \right| \leq \left| I_h^{(i,j),(i,j)} \right| + \frac{k}{2} \left| Q_h^{(i,j),(i,j)} \right| \leq
\left| l_{0}(\vartheta_x) \right| \left|  l_{0}(\vartheta_y) \right| \\
+ \frac{Dk}{2} \left( \left| l''_{0}(\vartheta_x) \right| \left| l_{0}(\vartheta_y) \right| + 2 \left|l'_{0}(\vartheta_x) \right| \left| l'_{0}(\vartheta_y)\right| + \left| l_{0}(\vartheta_x) l''_{0}(\vartheta_y) \right|
+ \frac{1}{\left| M \right|} \left( \left|l'_{0}(\vartheta_x) \right| \left| l_{0}(\vartheta_y) \right| + \left| l_{0}(\vartheta_x) \right| \left| l'_{0}(\vartheta_y) \right| \right)
\right).
\end{multline}
Since
\begin{equation}
\begin{gathered}
	\sup_{\vartheta \in [0,1]} \left| l_0(\vartheta) \right| = 
	\sup_{\vartheta \in [0,1]} \left| \frac{(1-\vartheta)(2-\vartheta)}{2} \right| = 1, \quad
	\sup_{\vartheta \in [0,1]} \left| l'_0(\vartheta) \right| = 
	\sup_{\vartheta \in [0,1]} \left| \frac{(2\vartheta-3)}{2h} \right| = \frac{3}{2h}, \\
	\sup_{\vartheta \in [0,1]} \left| l''_0(\vartheta) \right| = 
	\sup_{\vartheta \in [0,1]} \left| \frac{1}{h^2} \right| = \frac{1}{h^2},
\end{gathered}
\end{equation}
then (from \eqref{Abound})
\begin{equation}
\label{supAtilde}
\sup_{\vartheta_x,\vartheta_y \in [0,1]}
\left| A_h^{(i,j),(i,j)} \right| \leq 1 + 
\frac{Dk}{2} \left( \frac{13}{2 h^2} + \frac{3}{\left| M \right| h} \right)  \tilde{A}
\end{equation}
and finally the condition on $\left| \beta \right|$ is (from \eqref{condA}, taking $\tilde{A}$ as the right side of \eqref{supAtilde})
\[
\left| \beta \right| \geq
\frac{1}{2}
\left( 1 + \frac{Dk}{2} \left( \frac{13}{2 h^2} + \frac{3}{\left| M \right| h} \right) \right).
\]

\subsection{Transfer operators}
In this section we define the transfer operators $\mathcal{I}^h_H$ (restriction) and $\mathcal{I}^H_h$ (interpolation). We adopt a
geometric multigrid, then the Galerkin conditions (required in the Algebraic multigrid) are not satisfied, meaning that the interpolation
and restriction operators are not one the transpose of the other (multiplied by a suitable constant) and the coarse
grid operator is constructed in the same way as in the fine grid, without taking into account the transfer operators. In this
way we can afford treating in a simpler manner more complex-shaped geometries, maintaining the same sparsity pattern of the system
at all levels.
We use a cell-centered discretization and denote by $\Omega_h$ and $\Omega_{H}$ the fine and coarse grids respectively, with $\Omega_H$ built in the same way  as $\Omega_h$ (Sect.~\ref{sect:discspace}) but with spatial step $H=2h$.

\subsubsection{Restriction operator}
The defect $r_h$ contains both the defect of the inner relaxations and the defect of the relaxation of the boundary conditions.
Since the discrete operators of inner equations and boundary conditions scale with different powers of $h$, the defect may
show a sharp gradient crossing the boundary. For this reason, the restriction operator of the inner equation should involve only inside grid points and not ghost or inactive nodes to prevent the degradation of the multigrid performance \cite{COCO2013464,COCO2018299}.

\begin{figure}[htp]
\centering
\begin{minipage}[b]
	{.50\textwidth}
	\centering
	\includegraphics[width=0.45\textwidth]{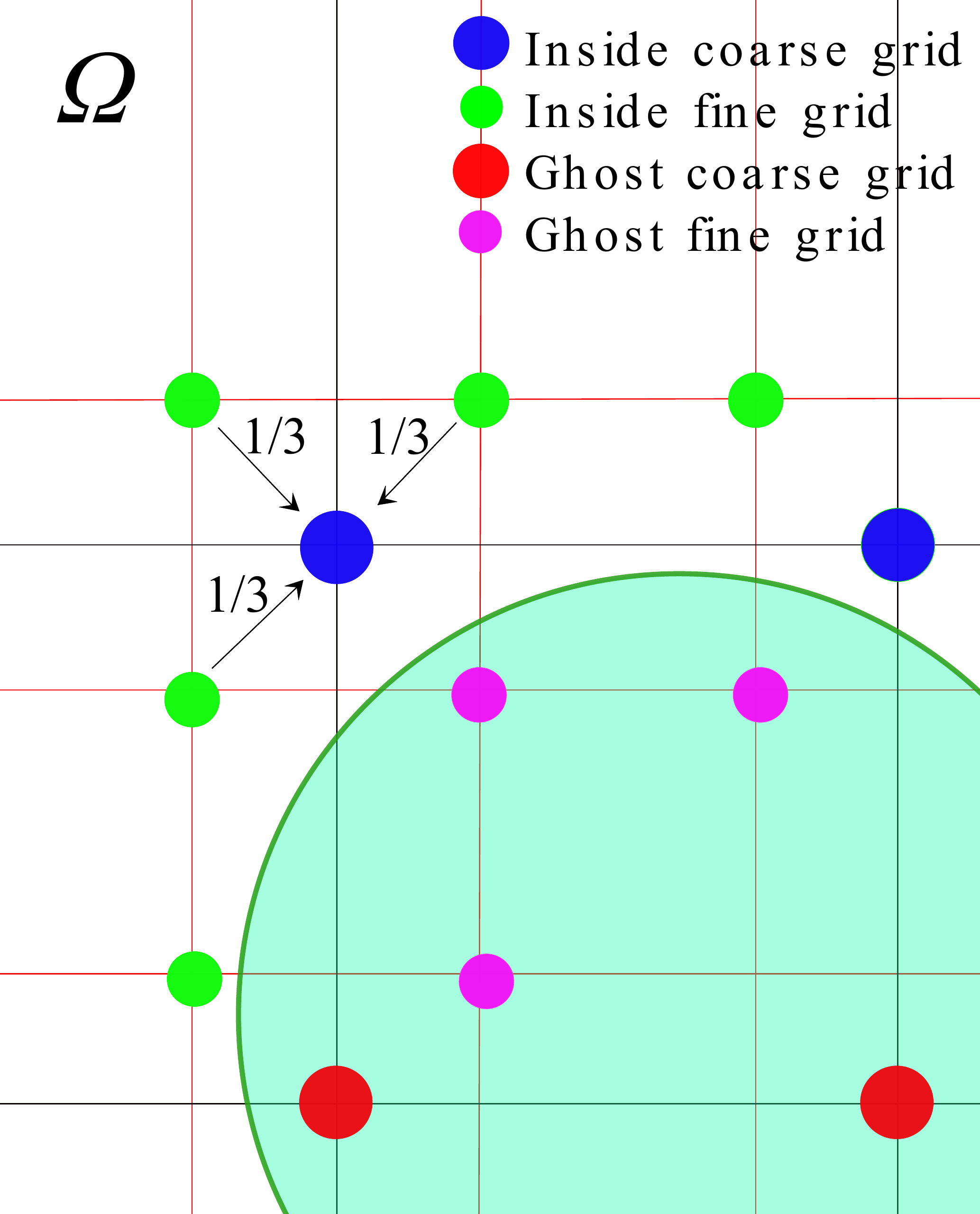}
\end{minipage}\hfill
\begin{minipage}[b]
	{.50\textwidth}
	\centering
	\includegraphics[width=0.45\textwidth]{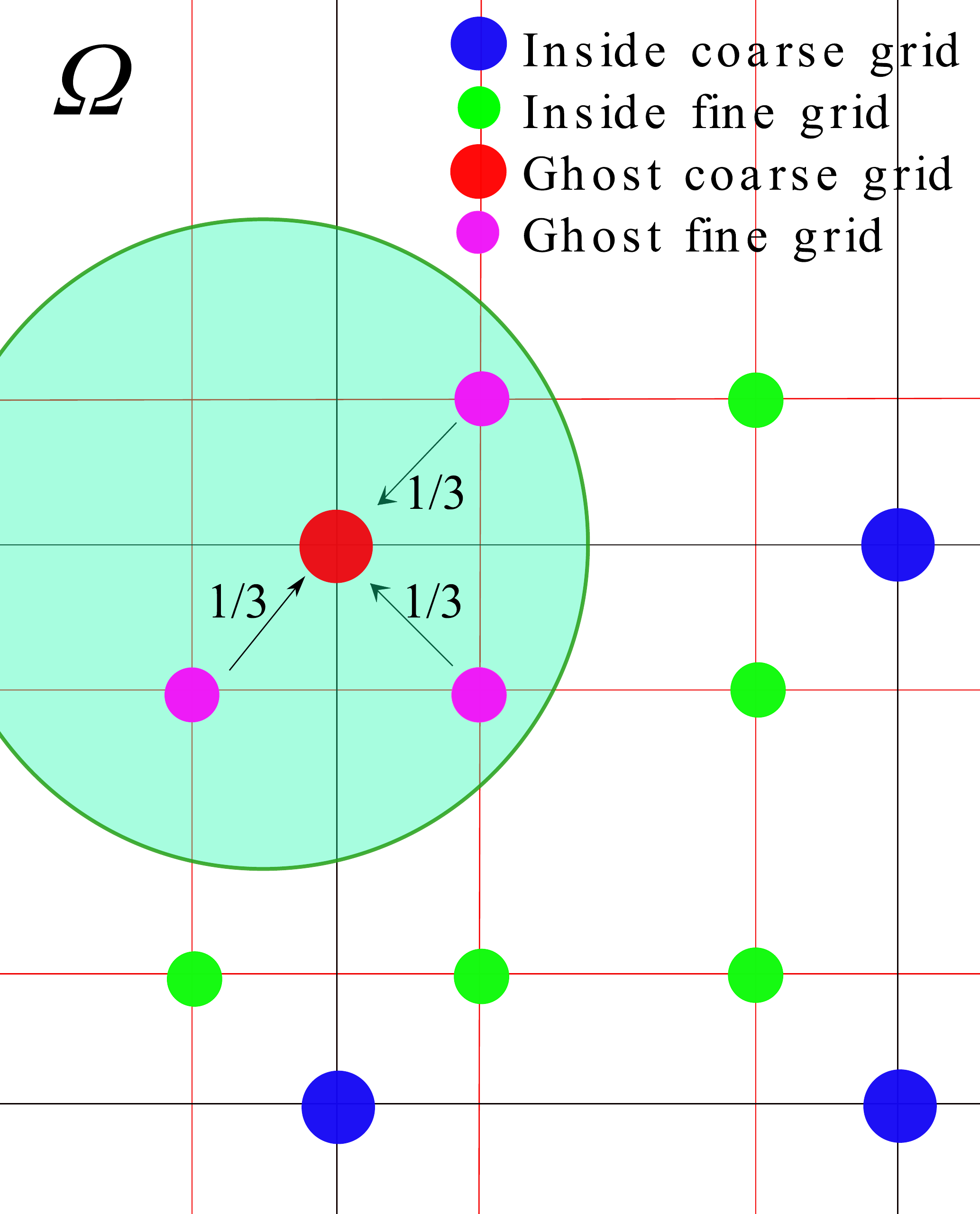}
\end{minipage}\hfill
\caption{\textit{Left panel: representation of the restriction operator and respective weights when the coarse grid point is close to the boundary $\Gamma_{\mathcal{B}}$.
		Right panel: restriction operator for ghost points and respective weights. 
}}
\label{restriction_operator_ghost}
\end{figure}

In practice, for each inner grid point $(x,y)$ of the coarse grid $\Omega_H$ we identify the four surrounding grid nodes of the fine grid $\Omega_h$, namely
$$
\mathcal{N}_{(x,y)} = \left\{ \left( x \pm \frac{h}{2},y\pm \frac{h}{2} \right) \right\}
$$
and we perform the restriction by averaging on the grid nodes of $\mathcal{N}_{(x,y)}$ that are inside $\Omega$ (see Fig.~\ref{restriction_operator_ghost}):
\begin{equation}
r_H(x,y) =
\mathcal{I}^{h}_{H} r_h(x,y)
= \frac{1}{\left| \mathcal{N}_{(x,y)} \cap \Omega_h \right|} \sum_{(x^*,y^*) \in \mathcal{N}_{(x,y)} \cap \Omega_h} r_h(x^*,y^*). 
\end{equation}
We observe that the restriction reverts to the classical restriction operator for cell-centered discretization and rectangular domains when $(x,y)$ is away from the boundary (see~\cite{Trottemberg:MG}):
\begin{equation*}
\mathcal{I}^{h}_{H} r_h(x,y)
= \frac{1}{4}\left[ r_h\left(x-\frac{h}{2},y - \frac{h}{2}\right) + r_h\left(x-\frac{h}{2},y + \frac{h}{2}\right)+r_h\left(x+\frac{h}{2},y - \frac{h}{2}\right) + r_h\left(x+\frac{h}{2},y + \frac{h}{2}\right) \right].
\end{equation*}
A similar approach is adopted to compute the restriction on a ghost point $(x,y) \in \mathcal{G}_H$:
\begin{equation}
\mathcal{I}^{h}_{H} r_h(x,y)
= \frac{1}{\left| \mathcal{N}_{(x,y)} \cap \mathcal{G}_h \right|} \sum_{(x^*,y^*) \in \mathcal{N}_{(x,y)} \cap \mathcal{G}_h} r_h(x^*,y^*). 
\end{equation}

\subsubsection{Interpolation operator}
\begin{figure}
\centering
\includegraphics[width=0.4\textwidth]{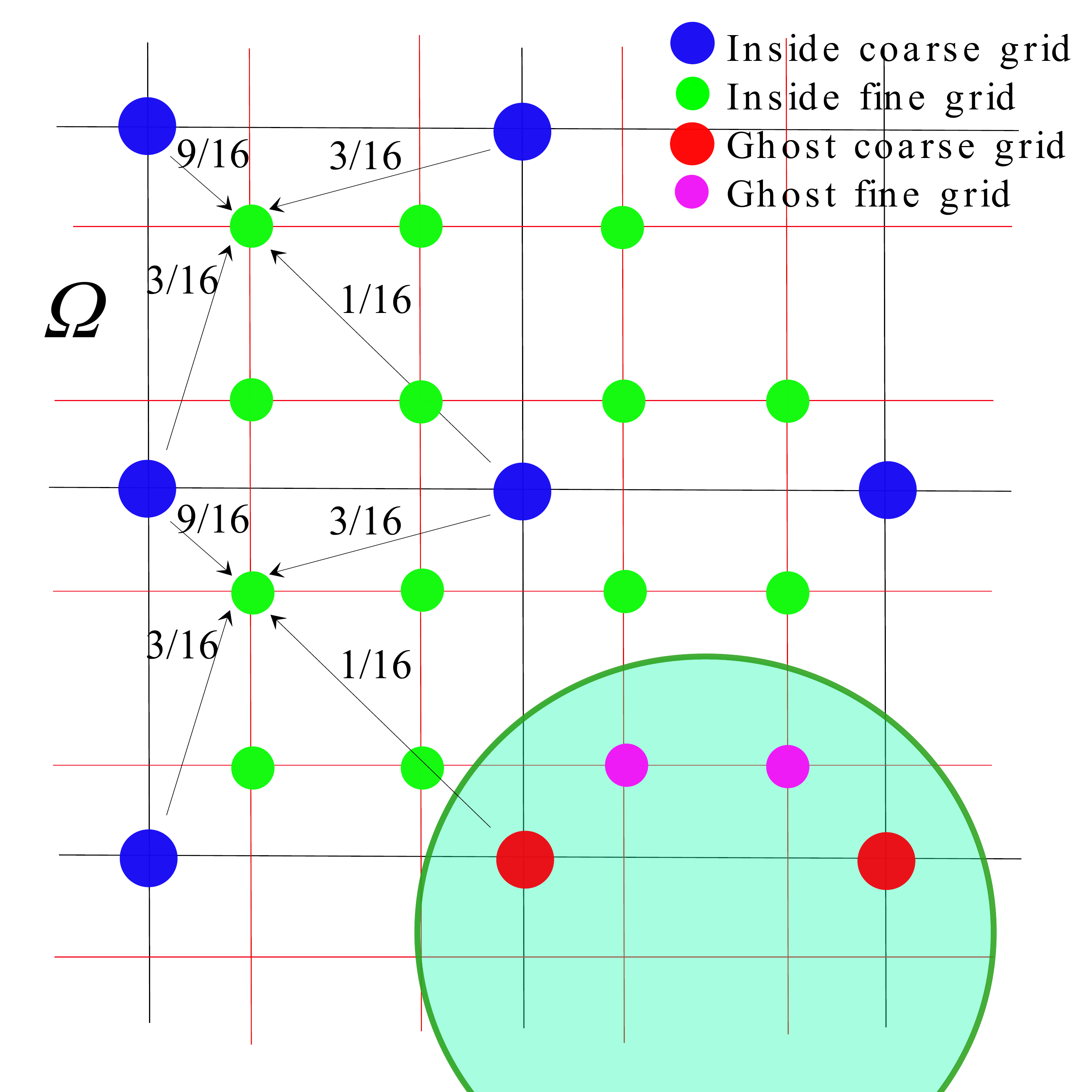}
\caption{\textit{Representation of the interpolation operator for cell-centered discretization and respective weights.
}}
\label{interpol}
\end{figure}


Once the residual equation $A_H e_H = r_H$ is solved on the coarse grid, the error $e_H$ is interpolated back to the fine grid $\mathcal{I}_h^H e_H = e_h$. Since the error $e_H$ is continuous across the boundary (it is the solution of the residual equation on the entire domain), then there is no need to separate the cases for inner and ghost points and the standard bilinear interpolation operator can be adopted across the entire domain. For example, for the inner grid node $(x,y)$ depicted in Fig.~\ref{interpol}, the interpolation reads:
\begin{multline}
e_h (x,y) = \mathcal{I}^H_h e_H (x,y) \\
= \frac{1}{16} \left( 9 e_H\left(x-\frac{h}{2},y+\frac{h}{2}\right) + 3 e_H\left(x+\frac{h}{2},y+\frac{h}{2}\right) + 3 e_H\left(x-\frac{h}{2},y-\frac{h}{2}\right) + e_H\left(x+\frac{h}{2},y-\frac{h}{2}\right) \right). 
\end{multline}

\section{Numerical results}\label{numtest}
\subsection{Accuracy test in 2D}\label{sect:2dtest}
In this section we test the accuracy of the method. We choose
an exact solution $ c _{exa}$ and augment the system~\eqref{Qexpr} as:
\begin{eqnarray}
\left\{
\begin{array}{l}
	\displaystyle \frac{\partial  c  }{\partial t} = D \Delta  c   + f \quad \text{ in } \Omega\\
	\displaystyle\frac{\partial  c  }{\partial n} + f_N = 0 \quad \text{ on } \Gamma_\mathcal{S}\\
	\displaystyle M\frac{\partial  c  }{\partial t} = MD\frac{\partial ^2 c  }{\partial \tau ^2}-D\frac{\partial  c  }{\partial n_\mathcal{B}} + f_B \quad \text{ on } \Gamma_\mathcal{B}\\
\end{array}
\right.
\label{accuracy2D}
\end{eqnarray}
choosing $\displaystyle f,f_N$ and $ f_B $ in such a way that $ c = c _{exa}$ is the exact solution. The computational domain is $\mathcal{S} =[-1,1]\times[-1,1]$, the radius of the bubble is $R_\mathcal{B} = 0.4$, while $M = 2 \times 10^{-4}$ 
and $ D = 0.1$. We choose the following exact solution:
\begin{align} 	\label{exact_rho}
c _{exa}(x,y,t) &= \cos(t)^2 c _0(x,y) + \sin(t)^2 c _1(x,y)&\\ \nonumber
c _0(x,y) &= \exp\left(-\frac{\left(x-x_0\right)^2+\left(y-y_0\right)^2}{\sigma_0}\right),\quad x_0 = 0,\quad y_0 = -0.6,\quad \sigma_0 = 0.1&\\ \nonumber
c _1(x,y) &= \exp\left(-\frac{\left(x-x_1\right)^2+\left(y-y_1\right)^2}{\sigma_1}\right),\quad x_1 = 0,\quad y_1 = -0.7,\quad \sigma_1 = 0.1.&
\end{align}
We compute the $L^1, {L}^2$ and ${L}^\infty$ norms of the relative error at $t=\pi/8$
\begin{align} 	\label{relativeerr}
e_\gamma = \frac{|| c - c _{exa}||_\gamma}{|| c _{exa}||_\gamma}, \quad \gamma = 1,2,\infty
\end{align}
for different values of $N$ and show the results in Table~\ref{table2D} and in the left panel of Fig.~\ref{fig:acc2D3D}, confirming numerically that the method is second order accurate.
\begin{table}
\centering
\begin{tabular}{||c||c||c||c||c||c||c||}
	\hline \hline
	$N$ & $e_1$ & {$p_1$}& $e_{2}$ & {$p_{2}$} & $e_{\infty}$ &  ${p_{\infty}}$  \\
	\hline \hline
	40 & 5.165E-02 & - & 4.377E-02 & - & 5.418E-02 & - \\
	80 & 1.234E-02 &  {2.066} & 1.045E-02 &  {2.066} & 1.453E-02 & {1.898}  \\
	160 & 3.054E-03 &  {2.014} & 2.584E-03 &  {2.016} & 3.696E-03 &  {1.975} \\
	320 & 7.766E-04 &  {1.975} & 6.472E-04 & {1.997} &  9.323E-04 & {1.987} \\
	\hline \hline
\end{tabular}
\caption{\textit{Relative errors $e_1$, $e_2$ and $e_\infty$  at $t=\pi/8$ and accuracy orders $p_1$, $p_2$ and $p_\infty$ in $\mathcal{L}^1, \mathcal{L}^2$ and $\mathcal{L}^\infty$ norms, respectively, for $ c $ in the 2D test of Sect.~\ref{sect:2dtest}. The exact solution is~\eqref{exact_rho}.}}
\label{table2D}
\end{table}

\subsection{Accuracy tests in 3D axisymmetric formulation}\label{sect:2dAStest}
\begin{figure}
\centering
\includegraphics[width=0.38\textwidth]{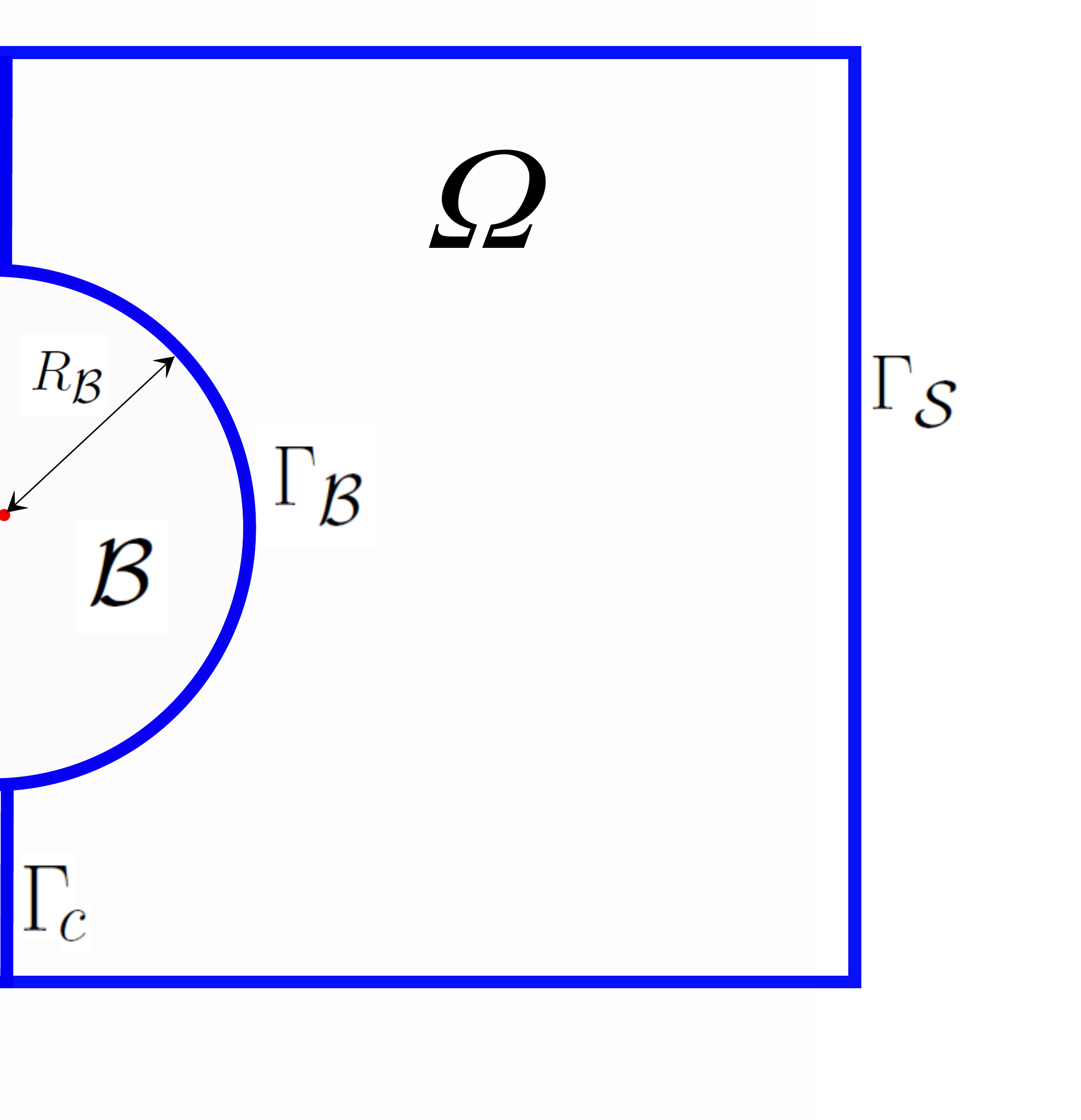}
\caption{\textit{Representation of the domain $\Omega$ in 3D axisymmetric: $\Gamma_{\mathcal{S}}$ is the external wall (top, right and bottom boundaries); $\mathcal{B} $ is the bubble with boundary $\Gamma_{\mathcal{B}}$ and radius $R_\mathcal{B}$; $\Gamma_{c}$ is the axis of symmetry (left boundary). }}
\label{3Daxisym}
\end{figure}
In this section we test the accuracy of the method for a 3D axisymmetric model (see Fig.~\ref{3Daxisym}).
The computational domain is $\mathcal{S} = [0,2]\times[-1,1]$, the radius of the bubble $R_\mathcal{B} = 0.4$ and $M = 2 \times 10^{-4}$.
The coordinates are the radial distance $\xi$ and the vertical coordinate $z$.
and $ D = 0.1$.
The problem reads (see Fig.~\ref{3Daxisym}):
\begin{eqnarray}
\left\{
\begin{array}{l}
	\displaystyle \frac{\partial  c }{\partial t} =D \left( \frac{\partial^2  c }{\partial \xi^2} + \frac{1}{\xi} \frac{\partial  c }{\partial \xi} + \frac{\partial^2  c }{\partial z^2} \right) \quad \text{ in } \Omega \\
	\displaystyle \nabla  c  \cdot n_\mathcal{S} = 0 \quad  \text{ on } \Gamma_\mathcal{S} \cup \Gamma_c\\
	\displaystyle M\frac{\partial  c }{\partial t} = MD\frac{\partial ^2 c }{\partial \tau ^2}-D\frac{\partial  c }{\partial n_\mathcal{B}}\quad \text{ on } \Gamma_\mathcal{B}\\
\end{array}
\right.
\label{system3D}
\end{eqnarray}

\begin{figure}[htp]
\centering
\hfill
\begin{minipage}[b]
	{.5\textwidth}
	\includegraphics[width=\textwidth]{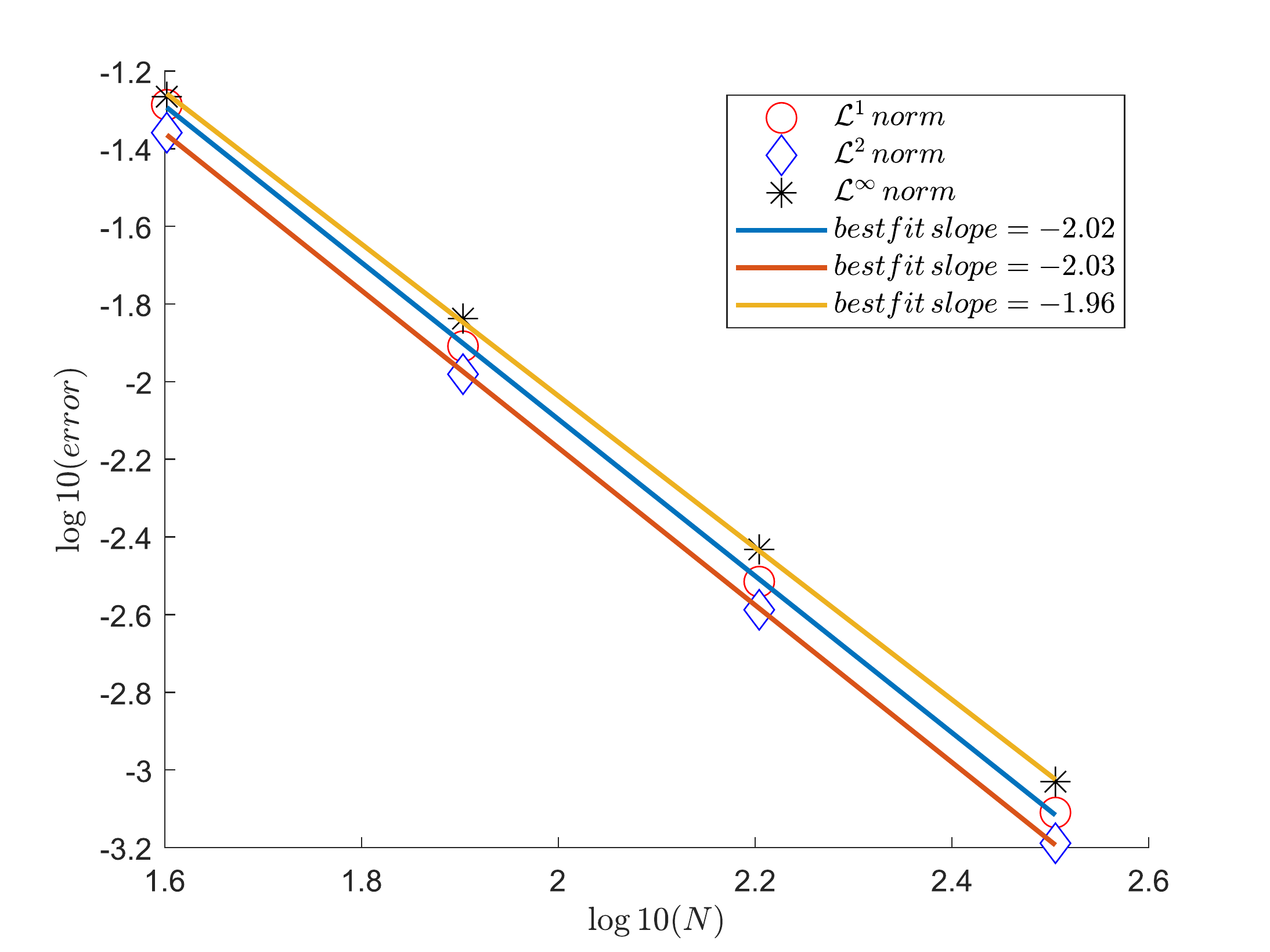}
\end{minipage}\hfill
\begin{minipage}[b]
	{.5\textwidth}
	\includegraphics[width=\textwidth]{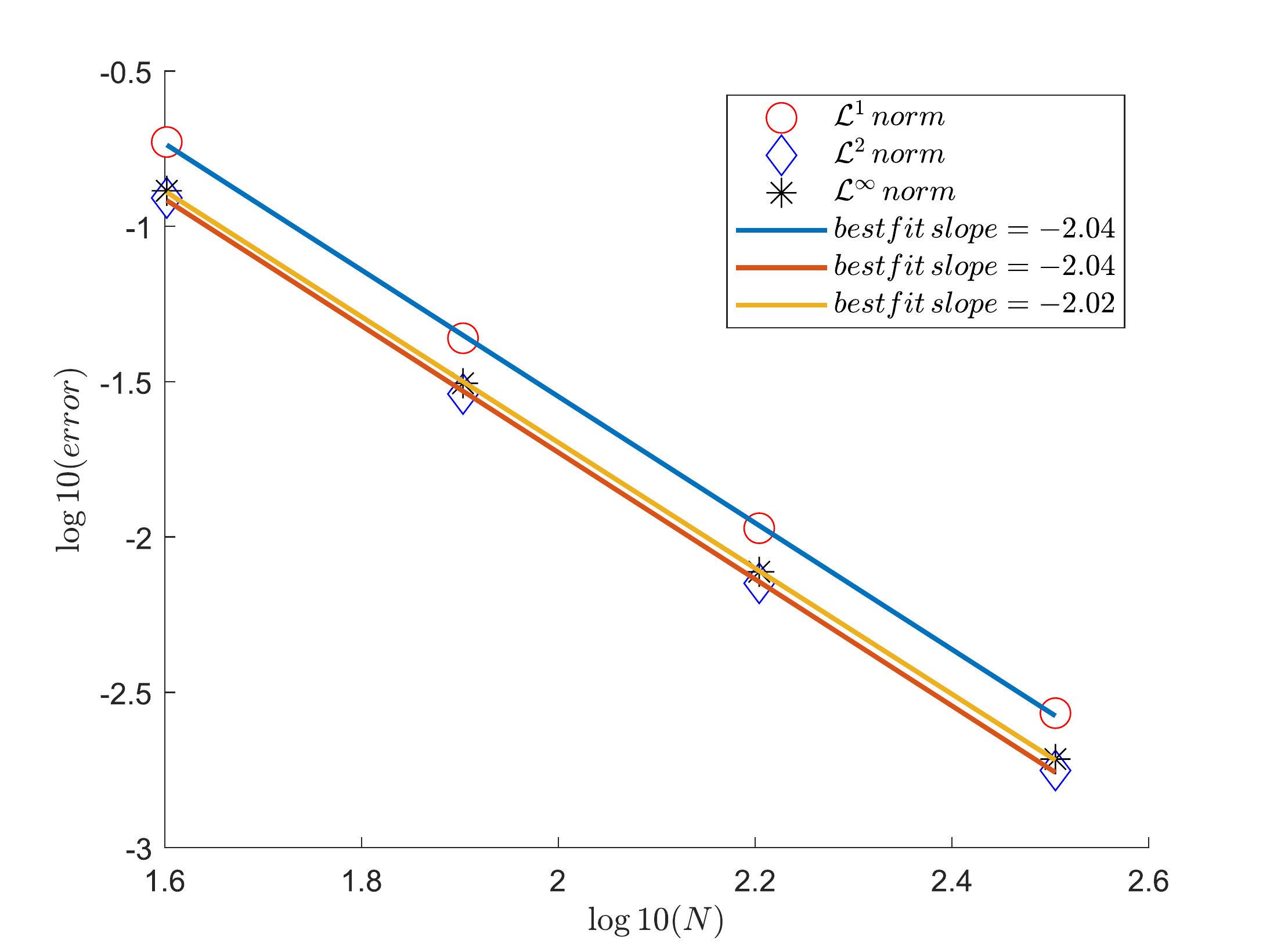}
\end{minipage}
\hspace*{\fill}
\caption{\textit{Representation of the relative errors at $t =\pi/8$ in logarithmic scale against the value of $N \in \{40,80,160,320\}$ in 2D (left panel, test of Sect.~\ref{sect:2dtest}, Table~\ref{table2D}) and 3D axisymmetric (right panel, test of Sect.~\ref{sect:2dAStest}, Table~\ref{table3D}).}}
\label{fig:acc2D3D}
\end{figure}

We choose the following exact solution:
\begin{align} 	\label{exact_rho3D}
c _{exa}(\xi,z,t) &= \cos(t)^2 c _0(\xi,z) + \sin(t)^2 c _1(\xi,z)&\\ \nonumber
c _0(\xi,z) &= \exp\left(-\frac{\left(\xi-\xi_0\right)^2+\left(z-z_0\right)^2}{\sigma_0}\right),\quad \xi_0 = 0,\quad z_0 = -0.6,\quad \sigma_0 = 0.1&\\ \nonumber
c _1(\xi,z) &= \exp\left(-\frac{\left(\xi-\xi_1\right)^2+\left(z-z_1\right)^2}{\sigma_1}\right),\quad \xi_1 = 0.1,\quad z_1 = -0.7,\quad \sigma_1 = 0.1&
\end{align}
and then we calculate the relative errors at $t=\pi/8$ as in Eq.~\eqref{relativeerr}.
Results are presented in table~\ref{table3D} and in the right panel of Fig.~\ref{fig:acc2D3D}, showing second order accuracy.

\begin{table}
\centering
\begin{tabular}{||c||c||c||c||c||c||c||}
	\hline
	\hline
	$N$ & $e_1$ & {$p_1$}& $e_{2}$ & {$p_{2}$} & $e_{\infty}$ &  ${p_{\infty}}$  \\
	\hline \hline
	40 & 1.868E-01 & - & 1.233E-01 & -& 1.299E-01 & -\\
	80 & 4.361E-02 &  {2.099} & 2.885E-02 &  {2.096} & 3.123E-02 &  {2.057}  \\
	160 & 1.075E-03 &   {2.021} & 7.103E-02 &  {2.022} & 7.728E-02 &  {2.015} \\
	320 & 2.708E-03 &  {1.989} & 1.771E-03 & {2.004} &  1.926E-03 &  {2.004} \\
	\hline \hline
\end{tabular}
\caption{\textit{Errors $e_1$, $e_2$ and $e_\infty$ at $t=\pi/8$ and accuracy orders $p_1$, $p_2$ and $p_\infty$ in $\mathcal{L}^1, \mathcal{L}^2$ and $\mathcal{L}^\infty$ norms, respectively,
		for $ c $ in the 3D axisymmetric test of Sect.~\ref{sect:2dAStest}. The exact solution is Eq.~\eqref{exact_rho3D}.}}
\label{table3D}
\end{table}

\section{Moving bubbles}\label{sect:moving_bubble}
In the following sections, unless otherwise specified, we consider the 3D axisymmetric model.

If the bubble moves over time it generates a fluid motion around it. In that case, the fluid domain depends on time $\Omega(t)$. Particles are subjected not only to the diffusion process but they are also transported by the moving fluids.
The motion of a fluid past an oscillating bubble is governed by the incompressible Navier-Stokes equations.
At low Reynolds numbers, the viscous forces are dominant and the convective term of the Navier-Stokes equations can be neglected so that the motion can then be described by the Stokes equations:

\begin{align}	\label{stokes_cylindrical}
\frac{\partial \textbf{u}}{\partial t} + \nabla p &= \frac{1}{Re} \nabla^2 \textbf{u} \quad \text{ in } \Omega(t) &\\ \nonumber
\nabla \cdot \textbf{u} &= 0 \quad \text{ in } \Omega(t) &
\end{align}
where \textbf{u} is the fluid velocity, $p$ is the pressure, $Re$ is the Reynolds number.

Driven by the application to Sorption Kinetics, we are interested in modelling the fluid dynamics generated by an oscillating bubble at extremely small amplitudes ($\sim 10^{-8} \text{ m}$), resulting in a low Reynolds number $Re<0.1$ (see~\cite{multiscale_mod}). Therefore, Stokes equations provide a reasonable approximation of the fluid dynamics for the problems investigated in this paper.

The 3D axisymmetric formulation (the coordinates are the radial distance $\xi$ and the vertical coordinate $z$) is completed by the following boundary conditions (see Fig.~\ref{3Daxisym}):

\begin{align}	\label{stokes_cylindricalBC}
\frac{\partial \textbf{u}}{\partial n} &= 0\quad  \text{ on } \Gamma_c &\\ \nonumber
\textbf{u} &= 0\quad  \text{ on } \Gamma_\mathcal{S} & \\ \nonumber
\textbf{u} \cdot \textbf{n} &=\textbf{u}_b \cdot \textbf{n} \quad  \text{ on } \Gamma_\mathcal{B}(t) & \\ \nonumber
\frac{\partial (\textbf{u} \cdot \tau)}{\partial n} & = 0 \quad  \text{ on } \Gamma_\mathcal{B}(t). &
\end{align}

The first one is the homogeneous boundary condition dictated by the axis symmetry, the second one is a no-slip boundary condition on the external wall, while the third and fourth ones are the free-slip boundary conditions at the boundary of the bubble ($\textbf{n}$ and $\tau$ are the normal and tangential vectors, respectively), where $\textbf{u}_b$ is the velocity of the bubble surface.


The Stokes equations \eqref{stokes_cylindrical} are discretized in time using Crank-Nicholson:
\begin{align}\label{stokesCN}
\displaystyle \frac{\textbf{u}^{(n+1)}-\textbf{u}^{(n)}}{\Delta t}+\nabla p^{(n+1/2)} &=\frac{1}{2Re}\left(\nabla^2\textbf{u}^{(n)} + \nabla^2 \textbf{u}^{(n+1)}\right) \quad \text{ in } \Omega^{(n+1)}&  \\ \nonumber
\nabla \cdot \textbf{u}^{(n+1)} &= 0 \quad \text{ in } \Omega^{(n+1)}& 
\end{align}

\begin{figure}[htp]
\centering
\includegraphics[width=0.5\textwidth]{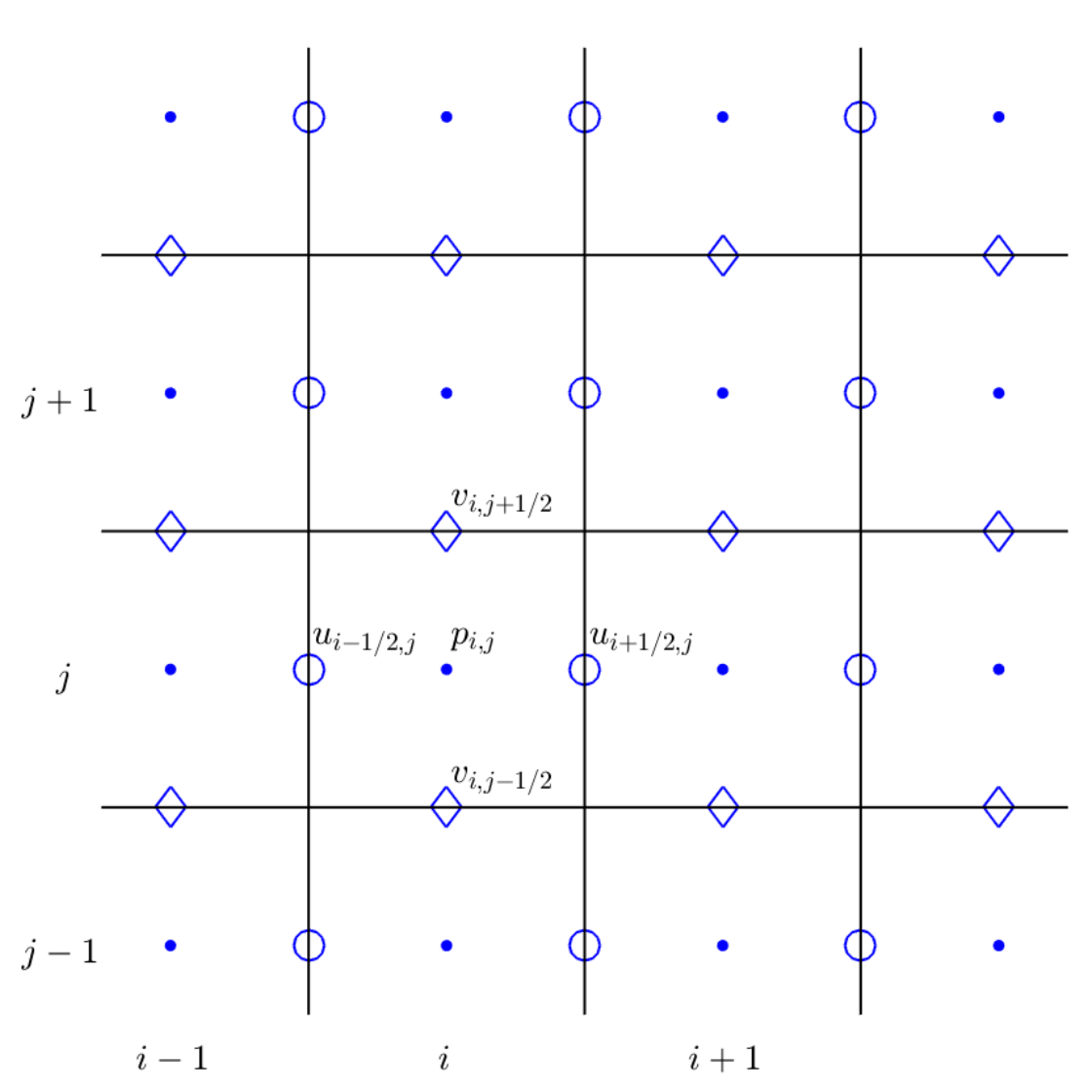}
\caption{\textit{Staggered grid for the Stokes problem. Horizontal velocity $u$ is defined on the middle points of the vertical edges of the cells (circle markers), vertical velocity $v$ is defined on the middle points of the horizontal edges of the cells (diamond markers), pressure $p$ is defined on the centers of the cells (dot markers).}}
\label{fig:staggered}
\end{figure}

The pressure $p$ and the velocity components $\textbf{u} = (u, v)$ are defined on a staggered grid (see Fig.~\ref{fig:staggered}): $p$ is defined at the centre of each cell, while $u$ and $v$ are defined at the mid points of vertical and horizontal side cells, respectively, i.e.~the so called Marker-and-cell discretization introduced by Harlow in the sixties \cite{harlow1965numerical}. 

The differential operators are discretized in space using central difference:
\begin{eqnarray*}
\displaystyle \frac{\partial p}{\partial \xi}\Big|_{i+1/2,j} = \frac{p_{i+1,j}-p_{i,j}}{h}, \quad \frac{\partial p}{\partial z}\Big|_{i,j+1/2} = \frac{p_{i,j+1}-p_{i,j}}{h} \\
\displaystyle \frac{\partial u}{\partial \xi} \Big|_{i+1/2,j} = \frac{u_{i+3/2,j}-u_{i-1/2,j}}{2h}, \quad \frac{\partial u}{\partial z}\Big|_{i+1/2,j} = \frac{u_{i+1/2,j+1}-u_{i+1/2,j-1}}{2h} \\
\displaystyle \frac{\partial v}{\partial \xi}\Big|_{i,j+1/2} = \frac{v_{i+1,j+1/2}-v_{i-1,j+1/2}}{2h}, \quad \frac{\partial v}{\partial z}\Big|_{i,j+1/2} = \frac{v_{i,j+3/2}-v_{i,j-1/2}}{2h} \\
\displaystyle \nabla^2u\Big|_{i+1/2,j} = \frac{u_{i+3/2,j}+u_{i-1/2,j}+u_{i+1/2,j+1}+u_{i+1/2,j-1}-4u_{i+1/2,j}}{h^2}\\
\displaystyle \nabla^2v\Big|_{i,j+1/2} = \frac{v_{i+1,j+1/2}+v_{i-1,j+1/2}+v_{i,j+3/2}+v_{i,j-1/2}-4v_{i,j+1/2}}{h^2}\\
\displaystyle \nabla \cdot \textbf{u}\Big|_{i,j} = \frac{u_{i+1/2,j}-u_{i-1/2,j}+v_{i,j+1/2}-v_{i,j-1/2}}{h}
\end{eqnarray*}

This discretization results in a linear system to be solved for $(\textbf{u}^{(n+1)},p^{(n+1/2)})$. This linear system is singular, due to the non uniqueness of $p^{(n+1/2)}$ (it is defined up to an additive constant) and it has to satisfy a compatibility condition to guarantee the existence of the solution.
Following the approach presented in~\cite{COCO2020109623}, the issue is circumvented by augmenting the problem \eqref{stokesCN} with an additional scalar unknown $\zeta \in \mathbb{R}$ and an additional equation for $p$ as follows:

\begin{align}\label{stokes_aug}
\displaystyle \frac{\textbf{u}^{(n+1)}-\textbf{u}^{(n)}}{\Delta t}+\nabla p^{(n+1/2)} &=\frac{1}{2Re}\left(\nabla^2\textbf{u}^{(n)} + \nabla^2 \textbf{u}^{(n+1)}\right) \quad \text{ in } \Omega^{(n+1)}& \\ \nonumber
\nabla \cdot \textbf{u}^{(n+1)} &= \zeta \quad \text{ in } \Omega^{(n+1)}&\\ \nonumber
\int_\Omega p\, d\Omega &= 0 \quad \text{ in } \Omega^{(n+1)}&
\end{align}
The problem consists then of finding $ (\textbf{u}^{(n+1)}, p^{(n+1/2)} , \zeta)$ that satisfies Eq.~\eqref{stokes_aug}.
We observe that the free divergence condition is not guaranteed, namely $\zeta$ is usually different from zero. However, the divergence decades with the same order of the method, i.e.~$\zeta = O(h^2)$, where $h$ is the spatial step, and then the overall accuracy order is not degraded (see~\cite{COCO2020109623} for more details).
We observed numerically that in the absence of a bubble, namely $\Omega = \mathcal{S}$, then $\xi = 0$ within machine precision. This is usually the case of domains without curved boundaries.

The third equation of (\ref{stokes_aug}) is discretized by standard mid-point rule, leading to the linear equation
\begin{equation}
\sum_{(\xi_i,z_j) \text{ internal points }} p_{i,j} = 0.
\end{equation}
The curved boundary is treated using the ghost-point technique described in~\cite{COCO2020109623}, similarly to the approach presented in Sect.~\ref{sect:discspace}.

\subsection{Test 1: pulsating bubble}\label{sect:pulsating}
\begin{figure}[htp]
\centering
\includegraphics[width=0.6\textwidth]{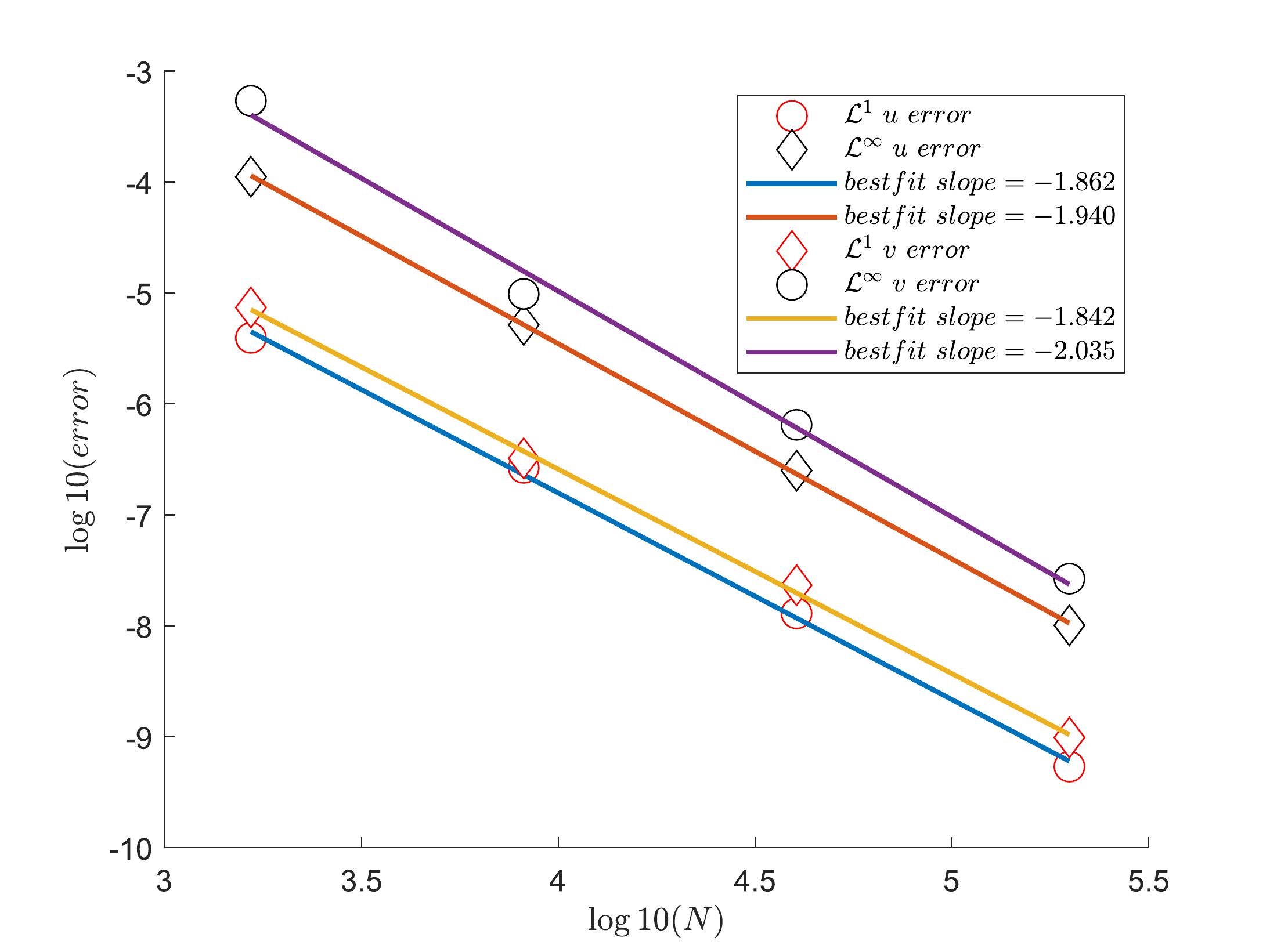}
\caption{\textit{Representation of the relative error for the oscillating bubble, Test 1, Sect~\ref{sect:pulsating}.
		We plot the errors for the two components of the velocity $u$ and $v$ in $\mathcal{L}^1$ and $\mathcal{L}^\infty$ norms. In this test $\Omega = [0,2]\times[-1,1]$, $A$ = 0.04, $nu$ = 1, $R_\mathcal{B}$ = 0.253 and final time $t_{fin} $ = 0.1.}}
\label{test_1_2}
\end{figure}

In this test we want to model the expansion/compression of a (pulsating) bubble, represented by a sphere $\mathcal{B}(t)$ centred at the origin and with radius:
\begin{equation}
R(t) = R_\mathcal{B} (1 + A\sin(\omega t))
\label{radius_t}
\end{equation} 
where $ R_\mathcal{B}$ is the radius of the bubble at time $t=0$.
The velocity of the bubble surface is 
\begin{equation}\label{bcub}
\textbf{u}_b (\xi,z) = R'(t) \, \textbf{n} =  A\, R_\mathcal{B} \, \omega \cos(\omega \, t) \, \textbf{n}, \text{ where } \textbf{n} = (\xi,z)/ \sqrt{\xi^2+z^2} \text{ and } \sqrt{\xi^2+z^2} = R(t). 
\end{equation}
The exact solution for the 3D axisymmetric Stokes problem~\eqref{stokes_cylindrical} with free-slip boundary conditions on the bubble surface (third and fourth equations of~\eqref{stokes_cylindricalBC}) in a semi-infinite domain $\Omega(t) = \left\{ (\xi,z) \colon 0<\xi<+\infty, \right.$ $ \left. \xi^2+z^2>(R(t))^2\right\}$ is:
\begin{equation}\label{exactu}
\textbf{u}_{\rm exa} = R'(t) \frac{(R(t))^2}{(\xi^2+z^2)^{3/2}} \cdot 
\begin{pmatrix}
\xi \\
z
\end{pmatrix}
, \quad p = R(t) (R''(t) R(t)+2 (R'(t))^2)/\sqrt{\xi^2+z^2}.
\end{equation}
In a finite domain $\Omega(t) = \mathcal{S} \backslash \mathcal{B}(t)$ we cannot prescribe the wall boundary conditions $\textbf{u}=0$ on the external boundary $\Gamma_\mathcal{S}$ otherwise the mass conservation is not guaranteed (since the volume of the bubble is not constant over time).
For this specific test we then prescribe the exact velocity \eqref{exactu} at $\Gamma_\mathcal{S}$.
We choose $R_\mathcal{B} = 0.253$, $A=0.04$ and $\omega= 2 \pi \, \nu$ with $\nu=50$ and we compute the numerical error at time $t_{fin}=0.1$ as the difference between the numerical solution and the exact solution \eqref{exactu}. The results are presented in Fig.~\ref{test_1_2}, where the second order accuracy is confirmed for both velocity components $u$ (radial component) and $v$ (vertical component) in norms $L^1$, $L^2$ and $L^\infty$.

\subsection{A steady computational bubble approach}\label{sect:fixed}
\begin{figure}[htp]
\centering
\hfill
\begin{minipage}[b]
	{.5\textwidth}
	\includegraphics[width=\textwidth]{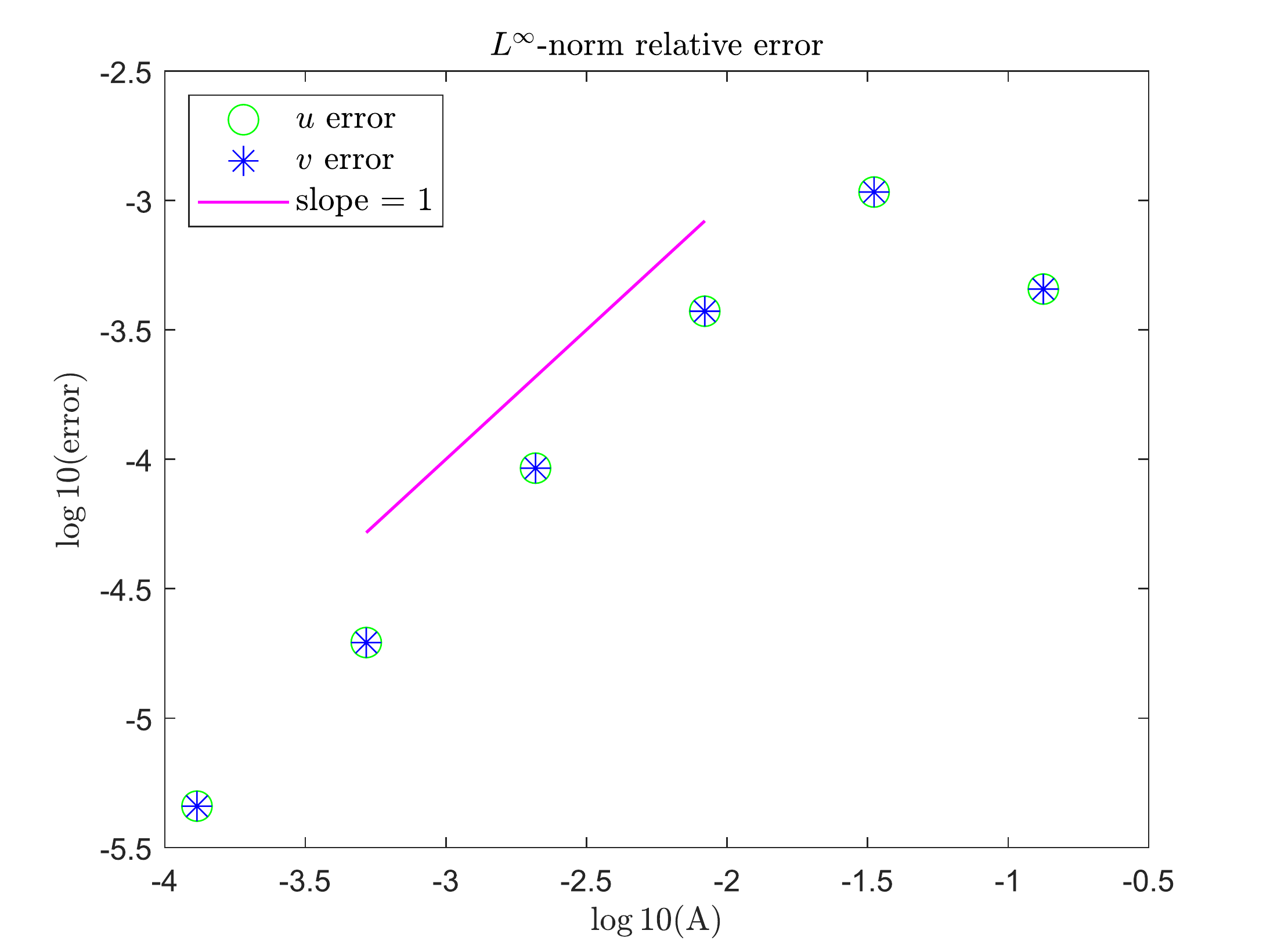}
\end{minipage}\hfill
\begin{minipage}[b]
	{.5\textwidth}
	\includegraphics[width=\textwidth]{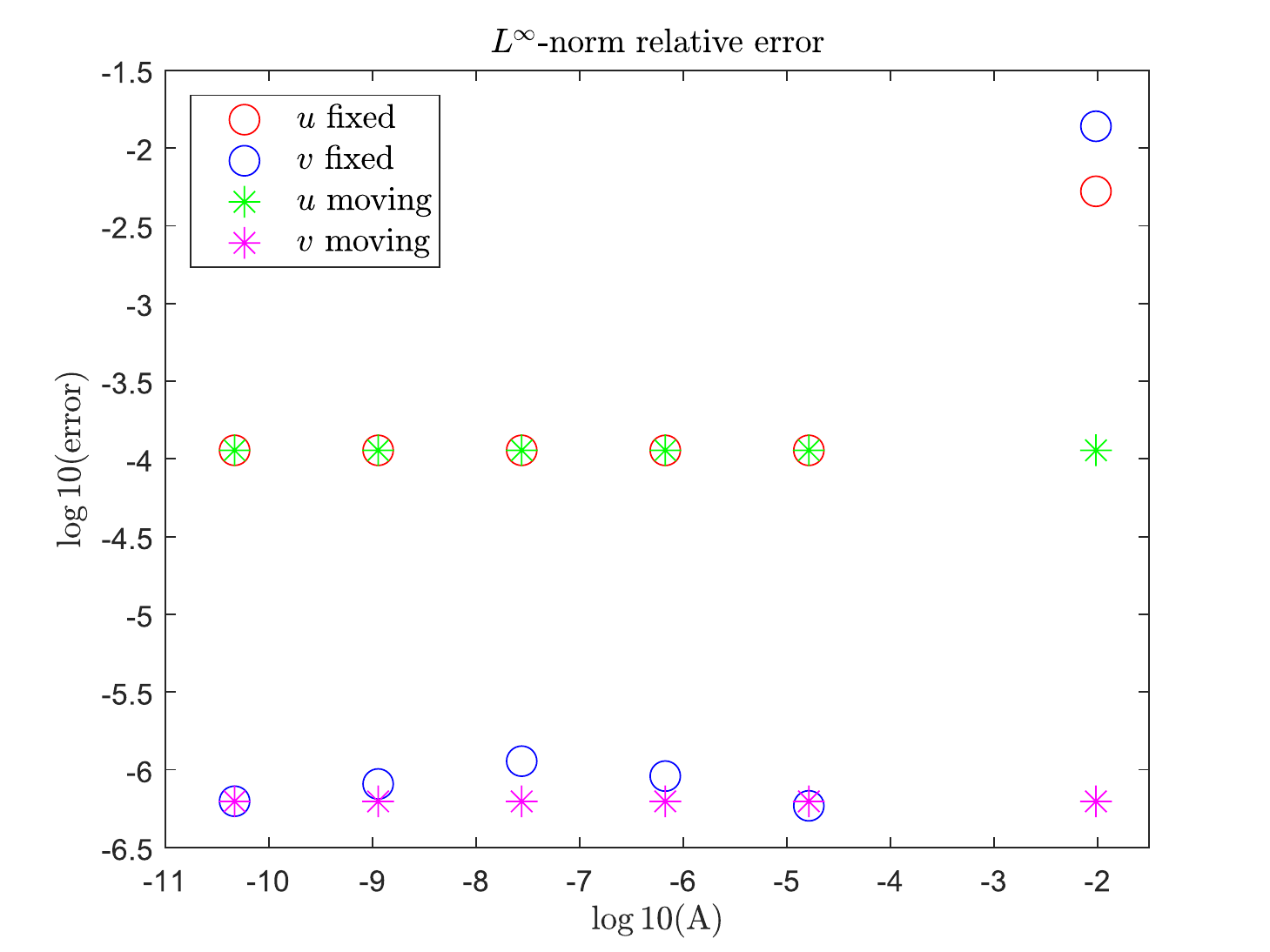}
\end{minipage}
\hspace*{\fill}
\caption{\textit{Left panel: relative difference between the numerical solutions of the two approaches described in Sections~\ref{sect:pulsating} and~\ref{sect:fixed} against the value of $A$.
		Right panel: relative errors of the two approaches computed as the difference between the numerical solutions and the respective exact solutions~\eqref{exactu} and~\eqref{exactufixed} against the value of $A$. In both plots we have $A \in \{h/64, h/16, h/4, h, 4h, 16h\}$, with $h=1/120$, and the time range is $[t_{\rm in},t_{\rm fin}] = [1,2]$, with $\nu = 1$.}
}
\label{2errors_1_2}
\end{figure}

When the amplitude of the bubble oscillation is sufficiently small compared to its dimensions, then $R(t) \approx R_\mathcal{B} $ and it is reasonable to simplify the model by assuming that the velocity of the surface bubble is assigned at a distance of $R_\mathcal{B}$ from the origin rather than $R(t)$: 
\begin{equation}\label{bcfixedtest1}
\textbf{u}_b(\xi,z) = A \, \omega \cos(\omega\,t) 
\cdot
\begin{pmatrix}
\xi \\
z
\end{pmatrix}
\text{ for } \sqrt{\xi^2+z^2} = R_\mathcal{B}.
\end{equation}
In this way, the computational domain does not move in time, since $\mathcal{B}(t)=\mathcal{B}(0)$ for any $t>0$ (steady computational bubble) and the fluid motion is generated purely from the boundary conditions.
The exact solution for the 3D axisymmetric Stokes problem~\eqref{stokes_cylindrical} with free-slip boundary conditions on the bubble surface (third and fourth equations of~\eqref{stokes_cylindricalBC}), with surface velocity defined by \eqref{bcfixedtest1}, in a semi-infinite domain $\Omega = \left\{ 0<\xi<+\infty, \xi^2+z^2>R_\mathcal{B}^2\right\}$ is:
\begin{equation}\label{exactufixed}
\textbf{u}^{\rm f}_{\rm exa} = R'(t) \frac{(R_\mathcal{B})^2}{(\xi^2+z^2)^{3/2}} \cdot 
\begin{pmatrix}
\xi \\
z
\end{pmatrix}, \quad p = R''(t) R_\mathcal{B}^2 /\sqrt{\xi^2+z^2}.
\end{equation}
The difference between the exact solutions \eqref{exactu} and \eqref{exactufixed} is $\textbf{u}_{\rm exa}-\textbf{u}^{\rm f}_{\rm exa} = \mathcal{O}(A)$. Therefore, for a fixed spatial step, the difference between the two approaches decades as $A \rightarrow 0$. This is confirmed numerically in Fig.~\ref{2errors_1_2} (left panel), where we compute the difference between the numerical solutions of the two approaches $\textbf{u}_h$ and $\textbf{u}_h^{f}$ at a fixed value of the spatial step $h=1/50$ and different values of $A$ over a period $T = 1/\nu$ as follows:
\begin{equation}
\displaystyle \frac{1}{T}\int_0^T\frac{\int_{\Omega(t)}\left|\textbf{u}_h-\textbf{u}^{f}_h\right|^p}{\int_{\Omega(t)}\left|\textbf{u}_h^{f}\right|^p}dt, \quad p = 1,2,\infty.
\label{period_error}
\end{equation}

In the right panel of Fig.~\ref{2errors_1_2} we compute the numerical error as the difference between the numerical and the exact solutions for both approaches at a fixed spatial step and different values of $A$. We observe that the second approach generally over performs the first one, although for sufficiently small values of $A$ the two errors reach a plateau, meaning that the overall error is dominated by the discretization error at the fixed spatial step. We conclude that when $A/R_\mathcal{B} \ll 1$ (as the cases investigated in this paper) it is more efficient to keep a steady computational bubble $\mathcal{B}$ and simulate the bubble motion by assigning a varying velocity at the computational surface
$\partial \mathcal{B}$. We will follow this approach in the following tests.

\subsection{Test 2: Oscillating bubble}\label{sect:oscbubble}
In the following tests we model the advection-diffusion process of particles in a moving fluid past an oscillating bubble:
\begin{equation}\label{pde2dcoupled}
\frac{\partial   c }{\displaystyle \partial t} = \nabla \cdot \left( D \nabla  c  -  c  \textbf{u} \right)
\end{equation}
where $\textbf{u} = (u,v)$ is the solution of the Stokes problem~\eqref{stokes_cylindrical}.
In general, we describe the motion of a bubble (that is not necessarily a sphere) by its parametric equations:
\begin{align*}
\xi(\theta,t) &= \xi_c(t)+\delta_\xi(t) \cos(\theta) &\\
z(\theta,t) &= z_c(t)+\delta_z(t) \sin(\theta) &
\end{align*}	
where $(\xi_c(t),z_c(t))$ is the centre of the bubble, while $\delta_\xi(t)$ and $\delta_z(t)$ regulate the deformation from a spherical shape.
We assume that at time $t=0$ the bubble is a sphere centred at the origin and with radius $R_\mathcal{B}$, then $\xi_c(0)=z_c(0)=0$ and $\delta_\xi(0)=\delta_z(0)=R_\mathcal{B}$.

\begin{figure}[htp]
\centering
\hfill
\begin{minipage}[b]
	{.249\textwidth}
	\centering
	\includegraphics[width=\textwidth]{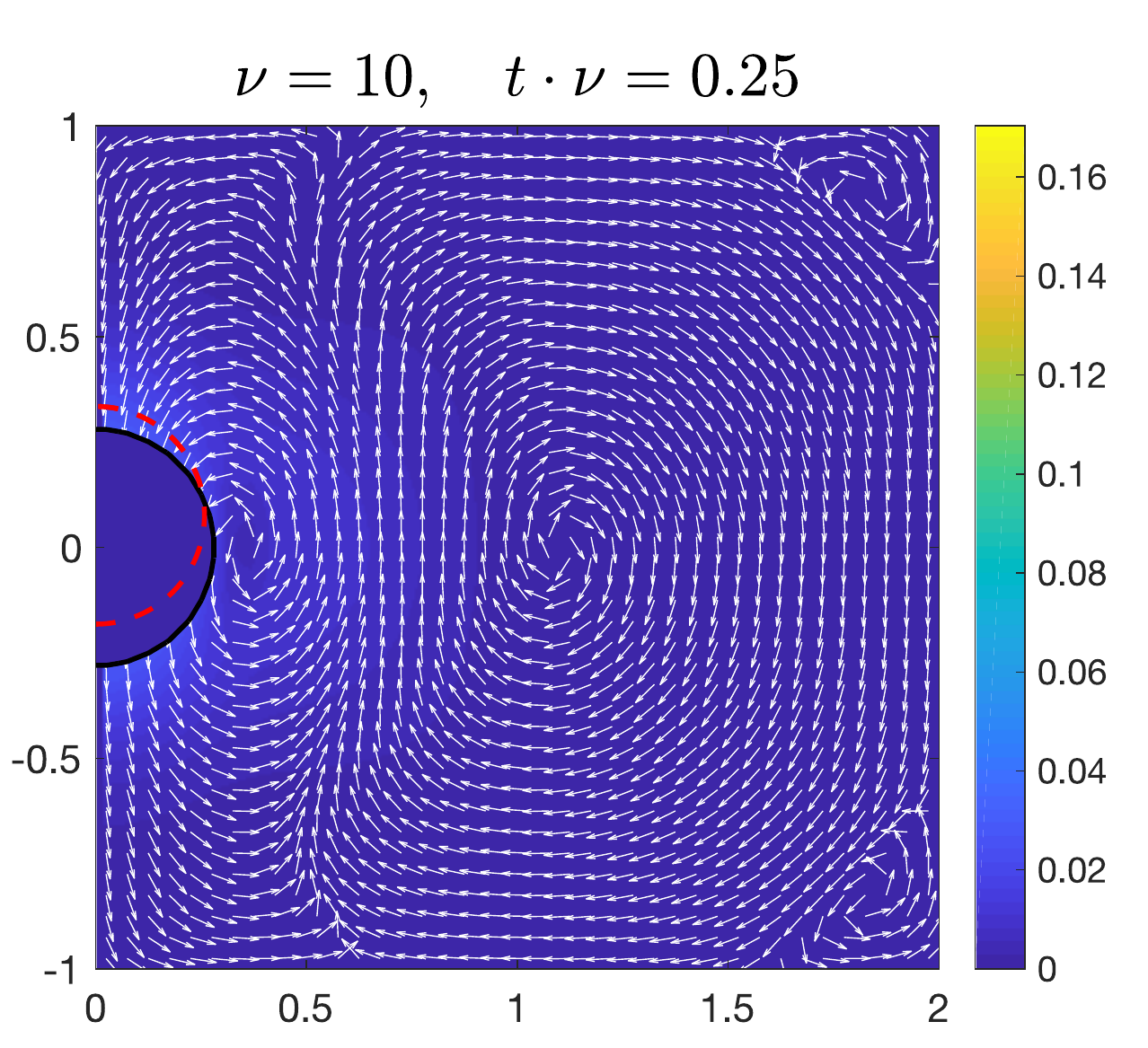}
\end{minipage}\hfill
\begin{minipage}[b]
	{.249\textwidth}
	\centering
	\includegraphics[width=\textwidth]{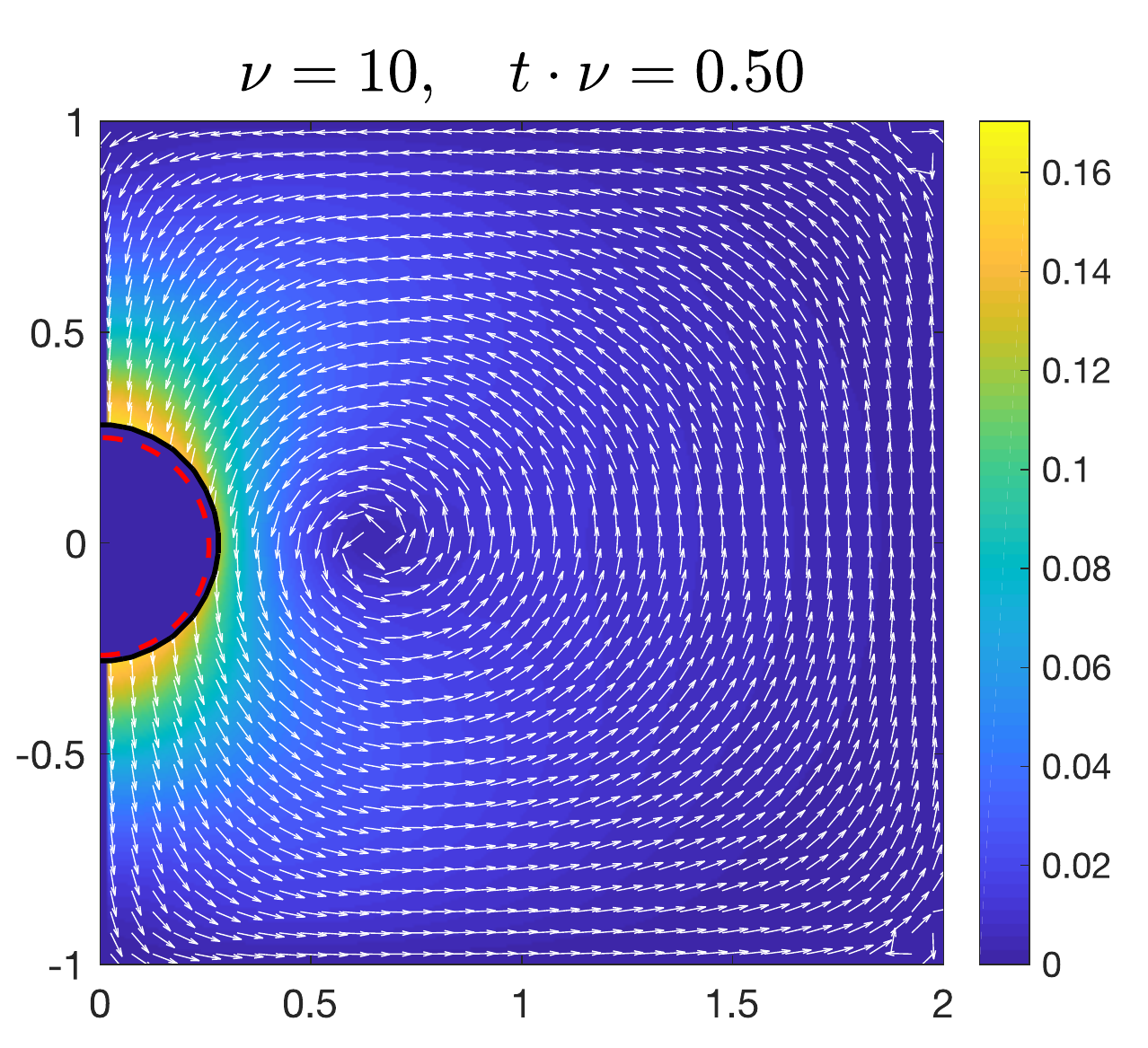}
\end{minipage}\hfill
\begin{minipage}[b]
	{.249\textwidth}
	\centering
	\includegraphics[width=\textwidth]{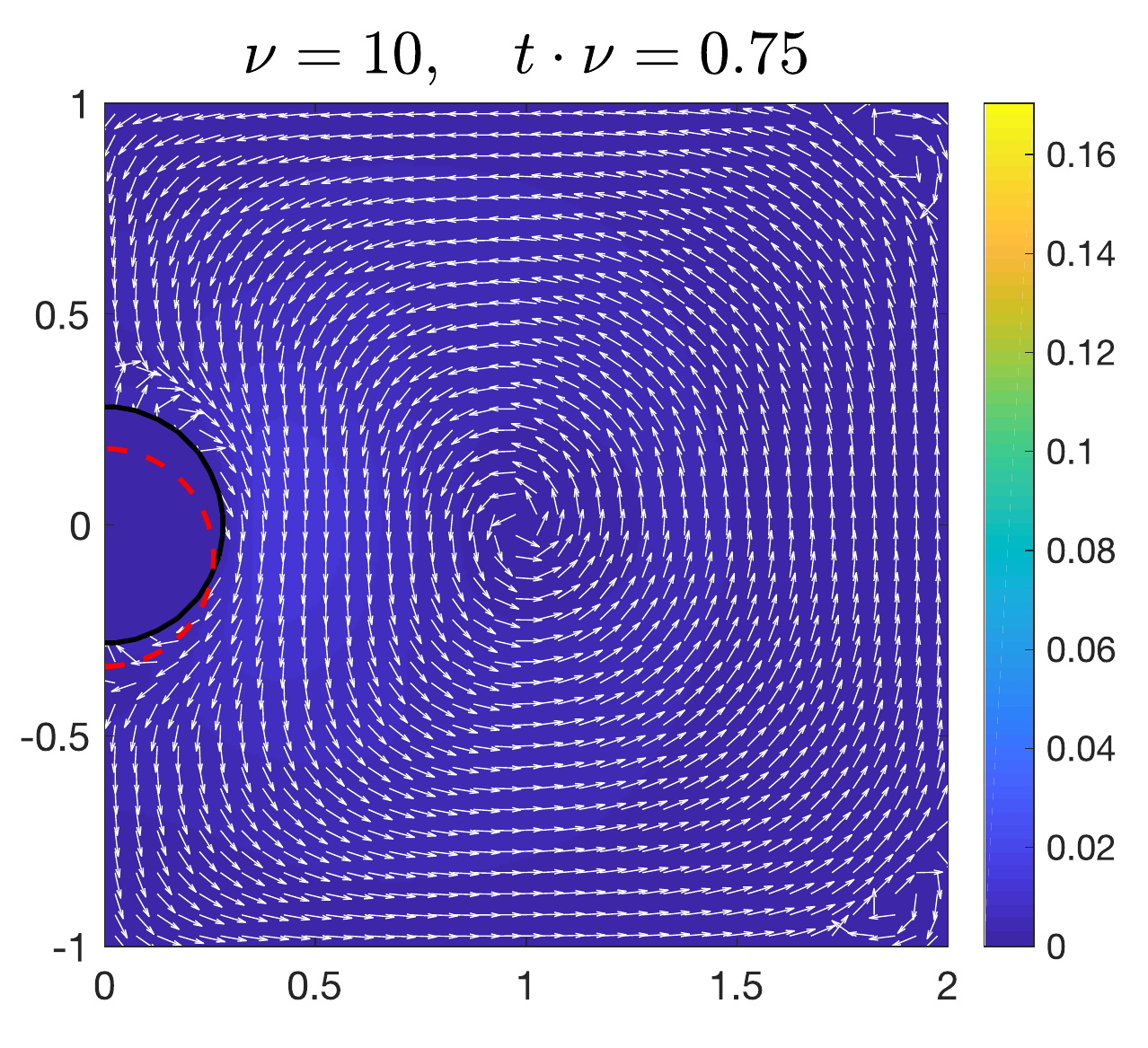}
\end{minipage}\hfill
\begin{minipage}[b]
	{.249\textwidth}
	\centering
	\includegraphics[width=\textwidth]{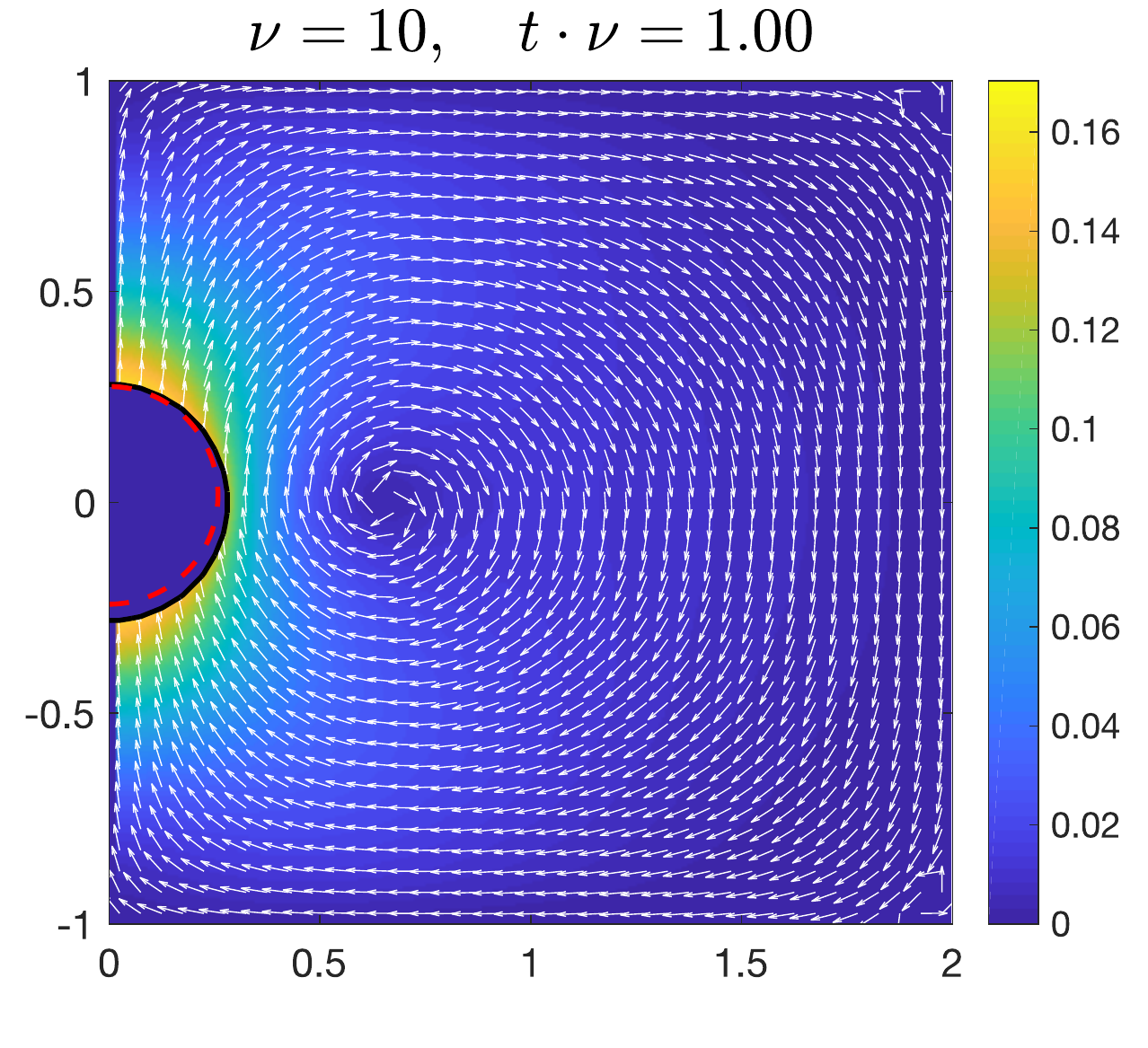}
\end{minipage}\hfill
\begin{minipage}[b]
	{.249\textwidth}
	\centering
	\includegraphics[width=\textwidth]{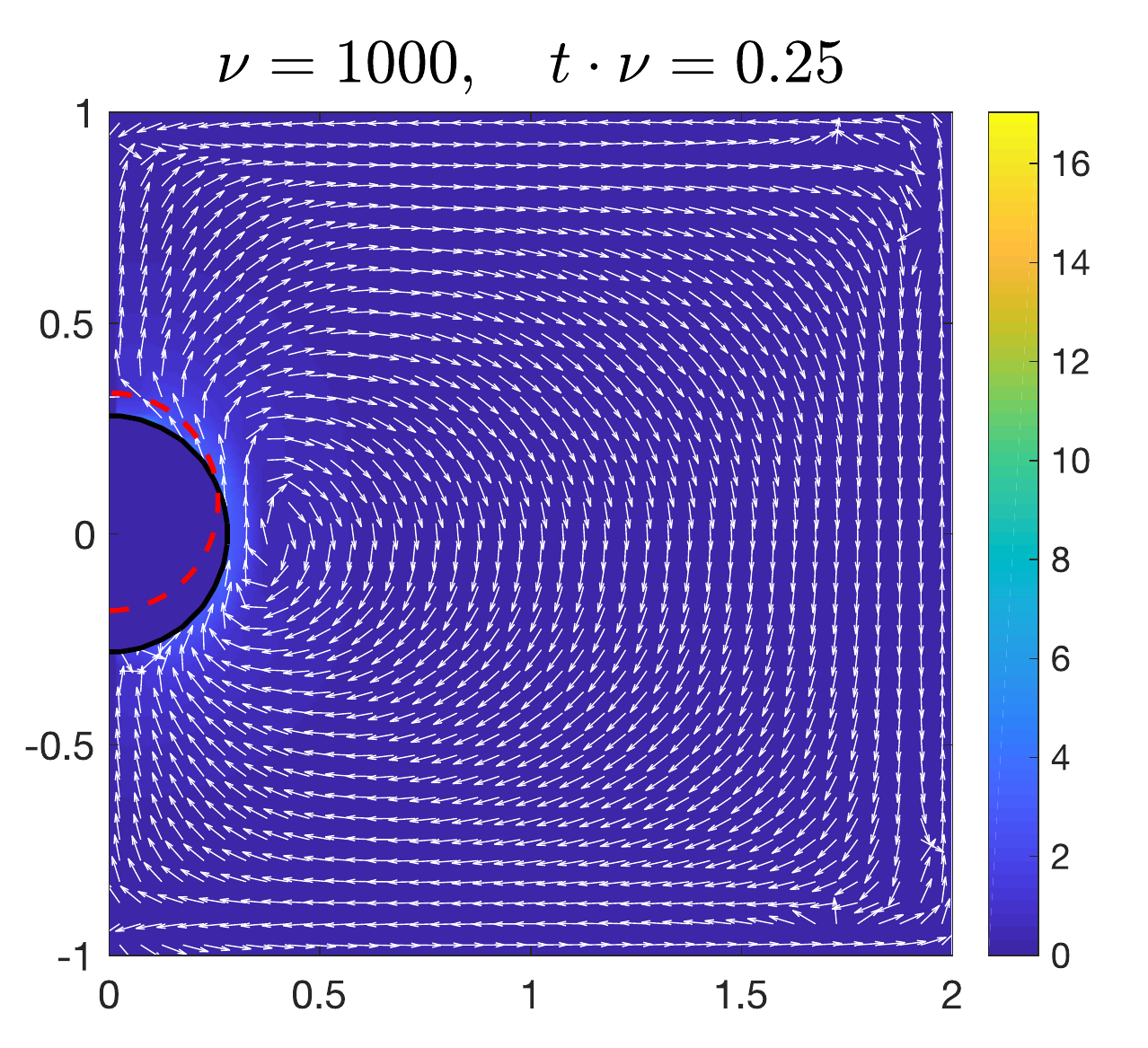}
\end{minipage}\hfill
\begin{minipage}[b]
	{.249\textwidth}
	\centering
	\includegraphics[width=\textwidth]{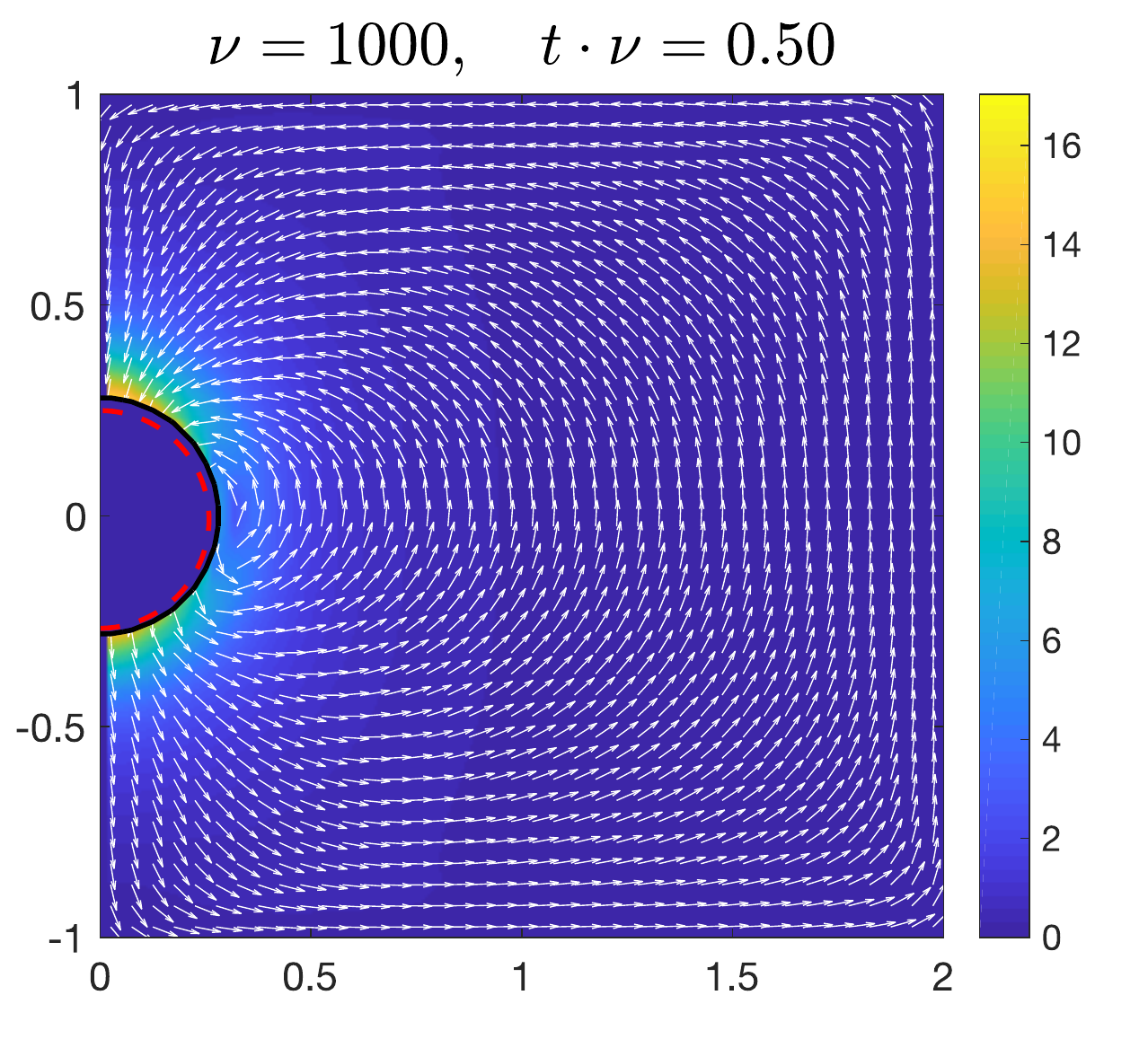}
\end{minipage}\hfill
\begin{minipage}[b]
	{.249\textwidth}
	\centering
	\includegraphics[width=\textwidth]{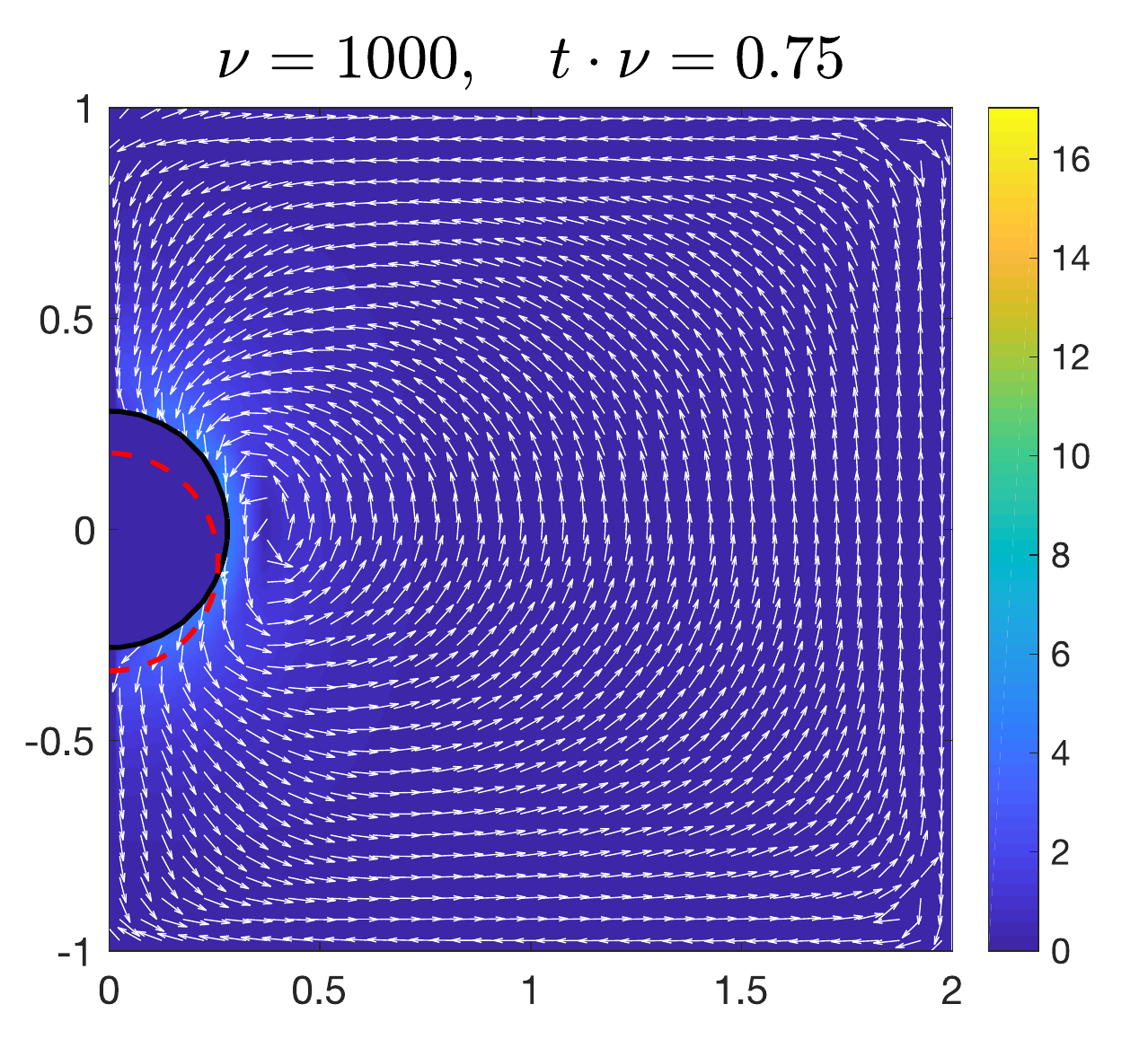}
\end{minipage}\hfill
\begin{minipage}[b]
	{.249\textwidth}
	\centering
	\includegraphics[width=\textwidth]{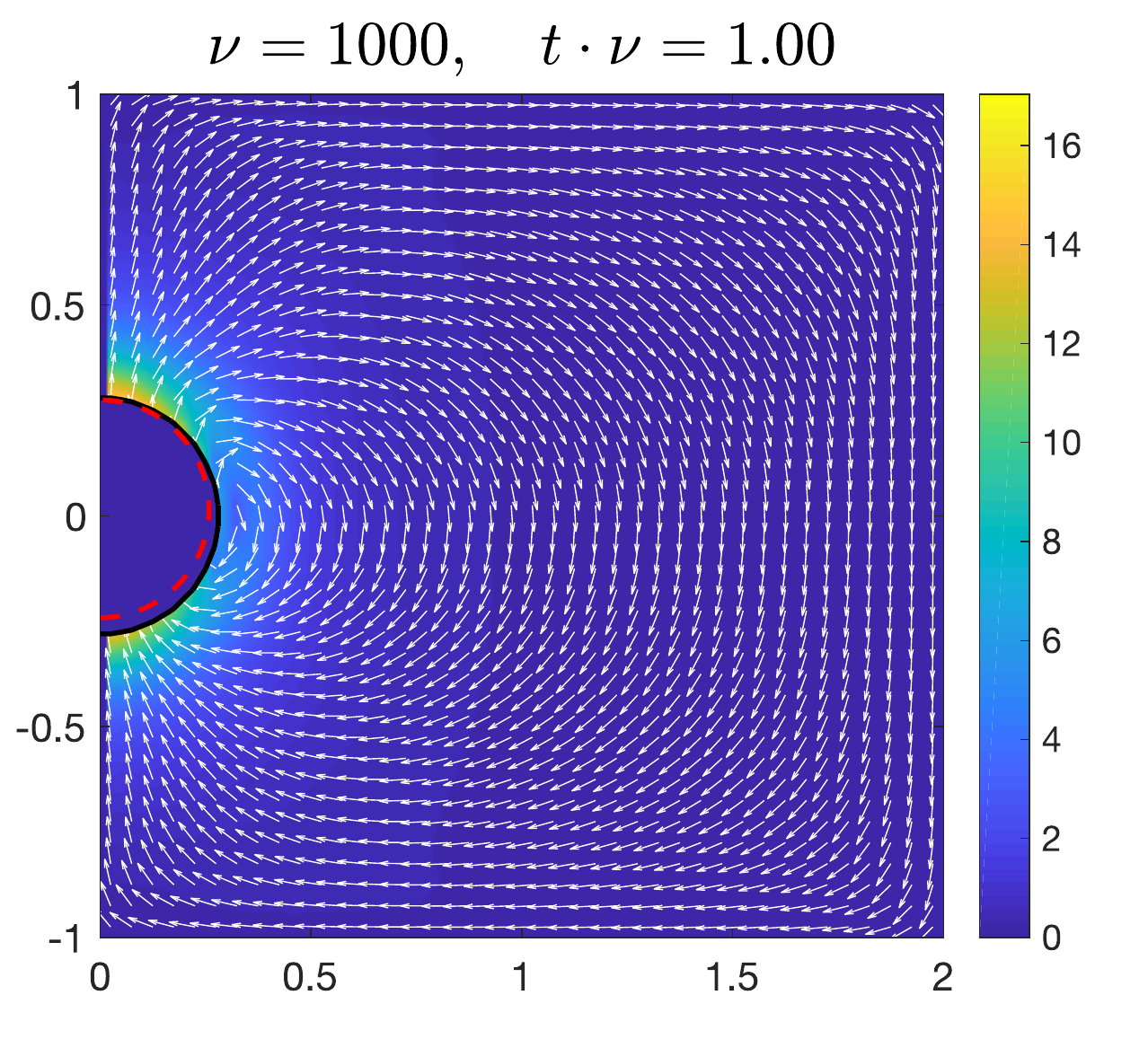}
\end{minipage}\hfill
\begin{minipage}[b]
	{.249\textwidth}
	\centering
	\includegraphics[width=\textwidth]{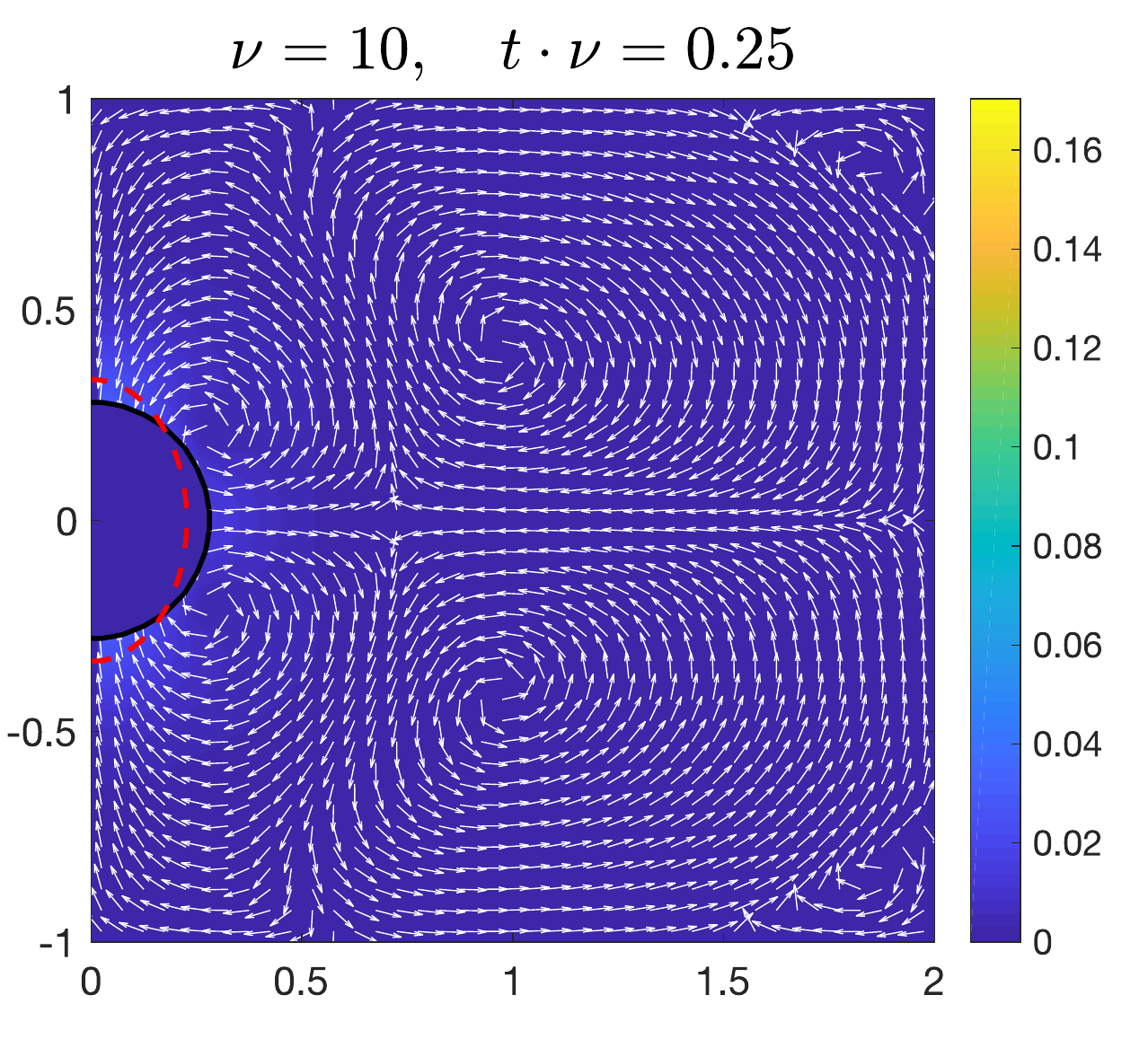}
\end{minipage}\hfill
\begin{minipage}[b]
	{.249\textwidth}
	\centering
	\includegraphics[width=\textwidth]{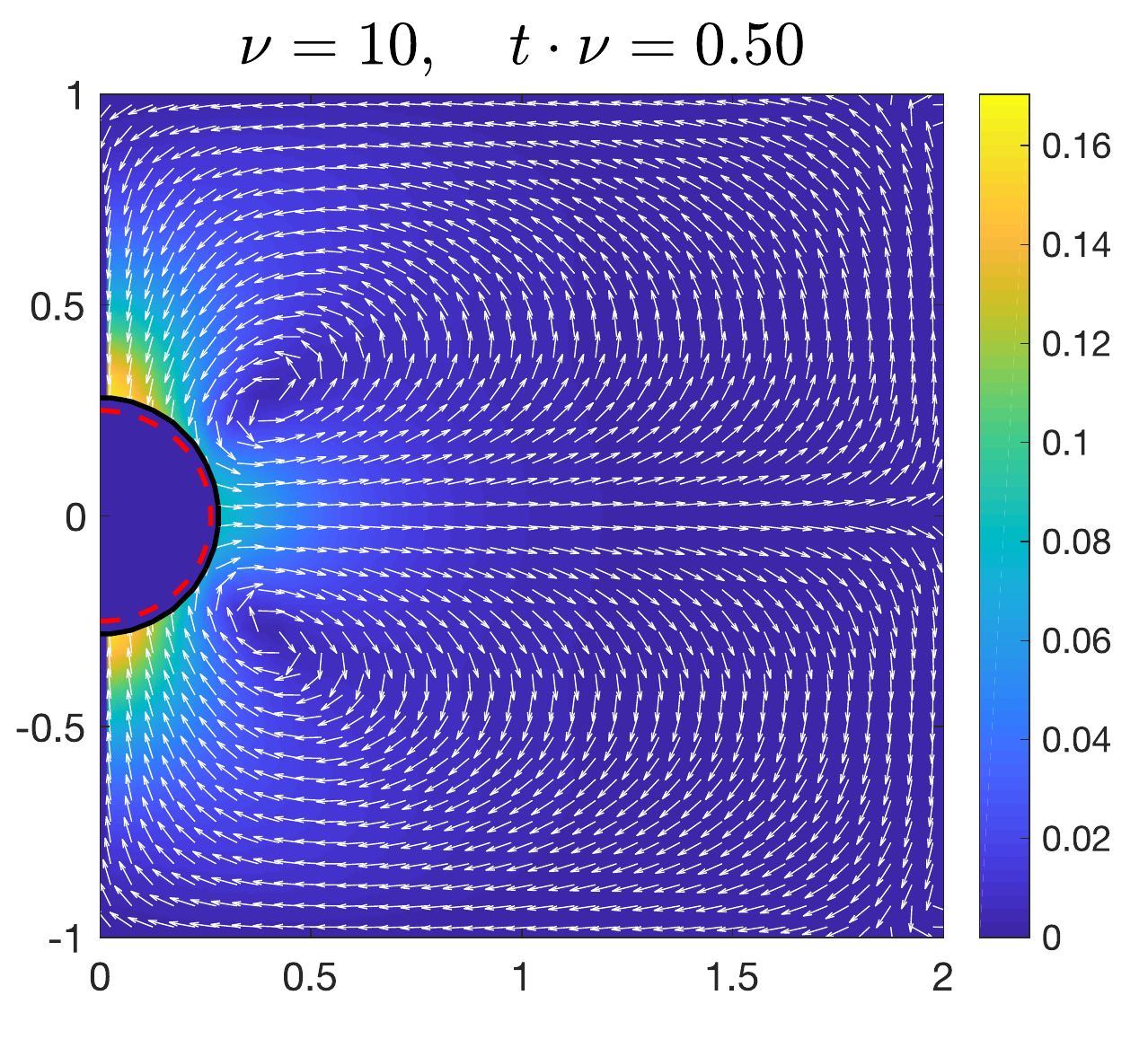}
\end{minipage}\hfill
\begin{minipage}[b]
	{.249\textwidth}
	\centering
	\includegraphics[width=\textwidth]{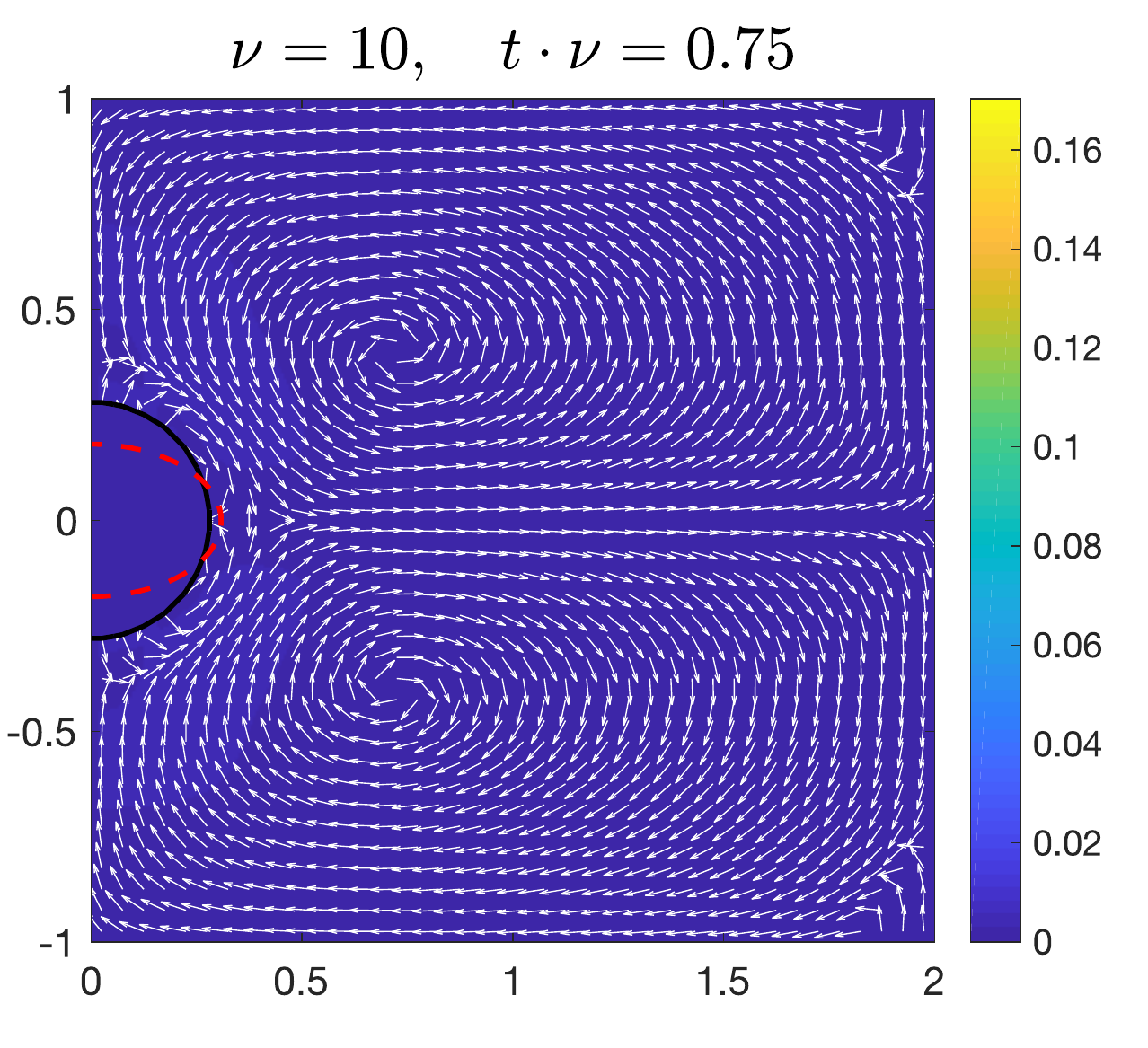}
\end{minipage}\hfill
\begin{minipage}[b]
	{.249\textwidth}
	\centering
	\includegraphics[width=\textwidth]{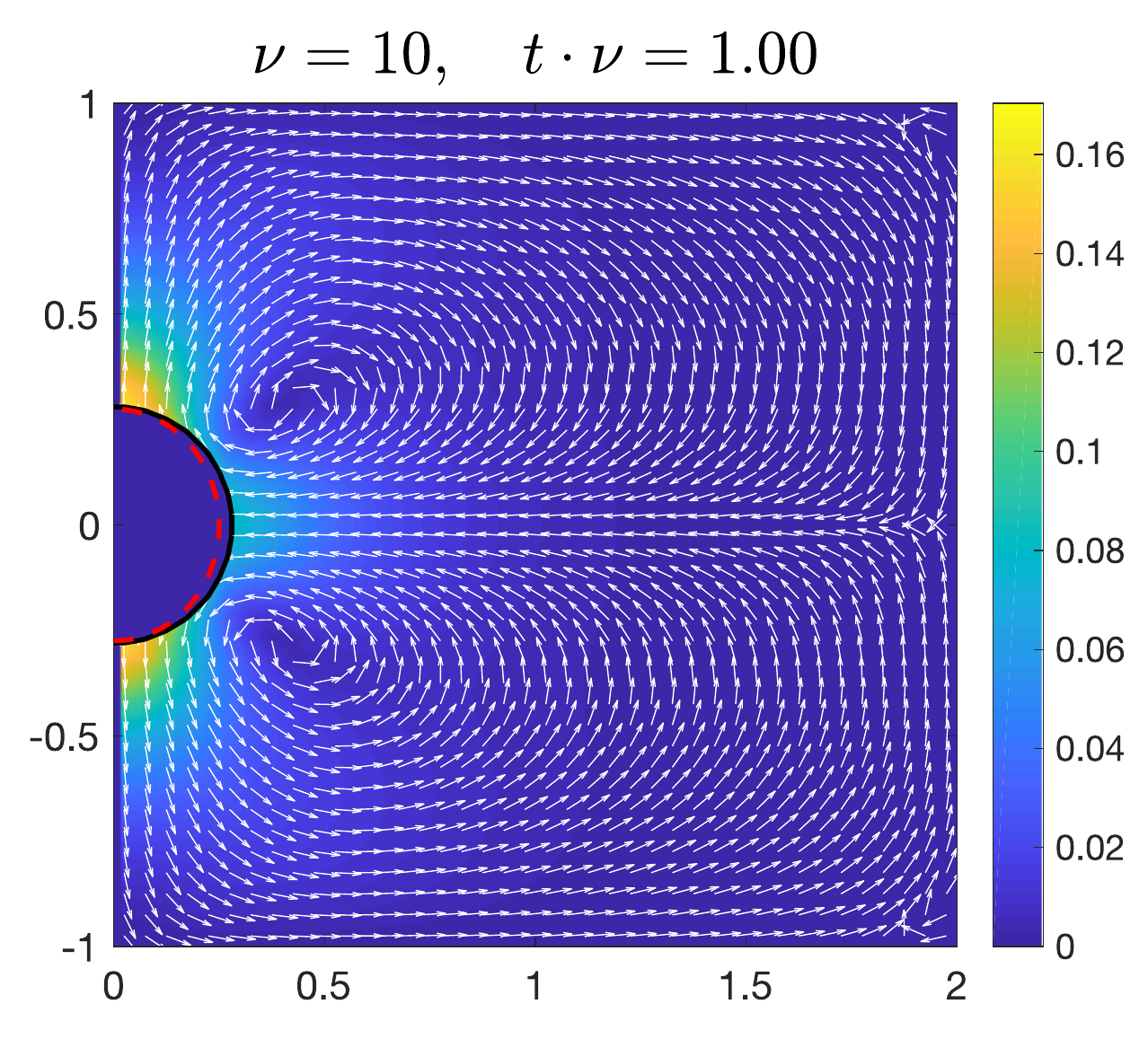}
\end{minipage}\hfill
\begin{minipage}[b]
	{.249\textwidth}
	\centering
	\includegraphics[width=\textwidth]{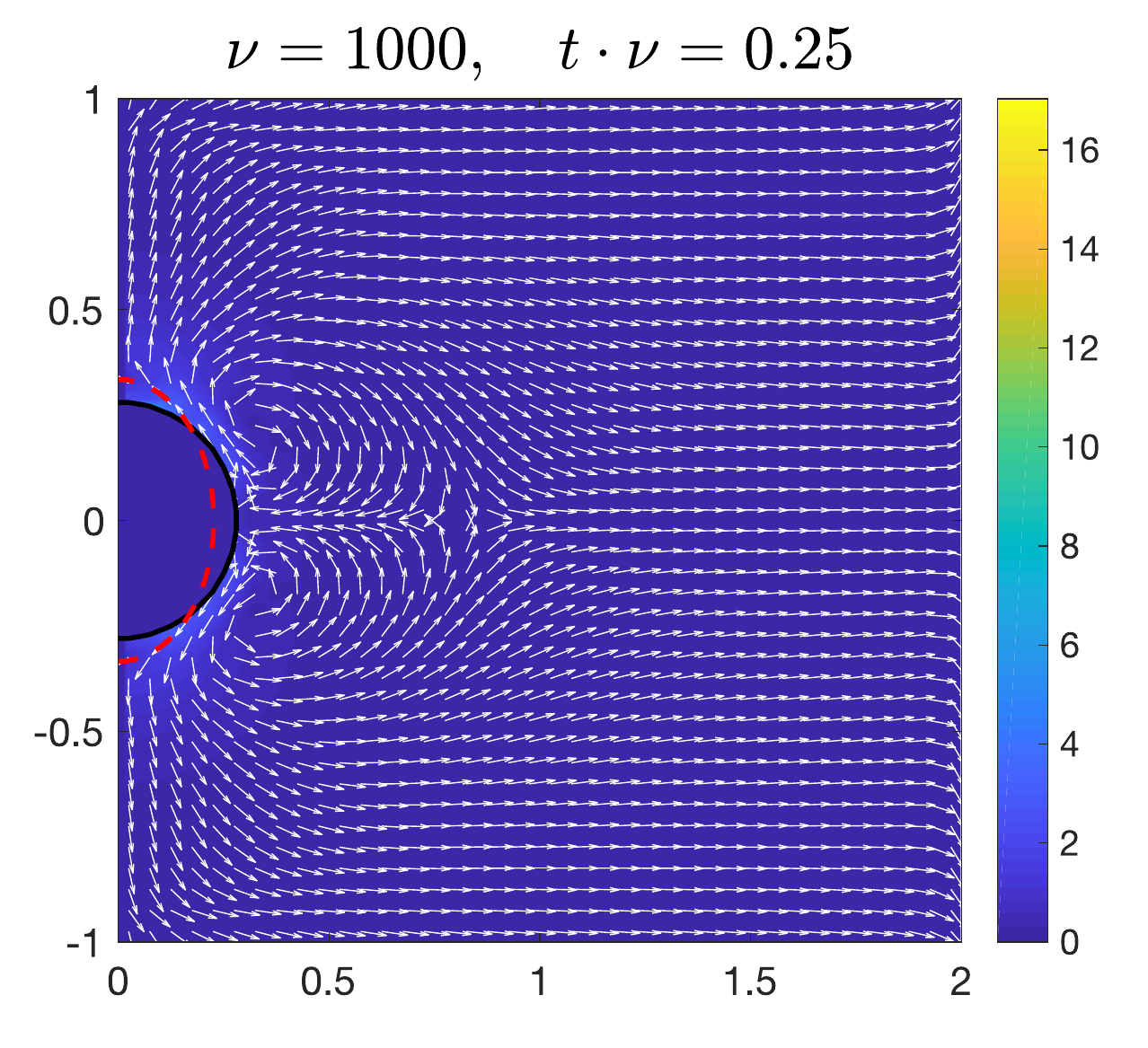}
\end{minipage}\hfill
\begin{minipage}[b]
	{.249\textwidth}
	\centering
	\includegraphics[width=\textwidth]{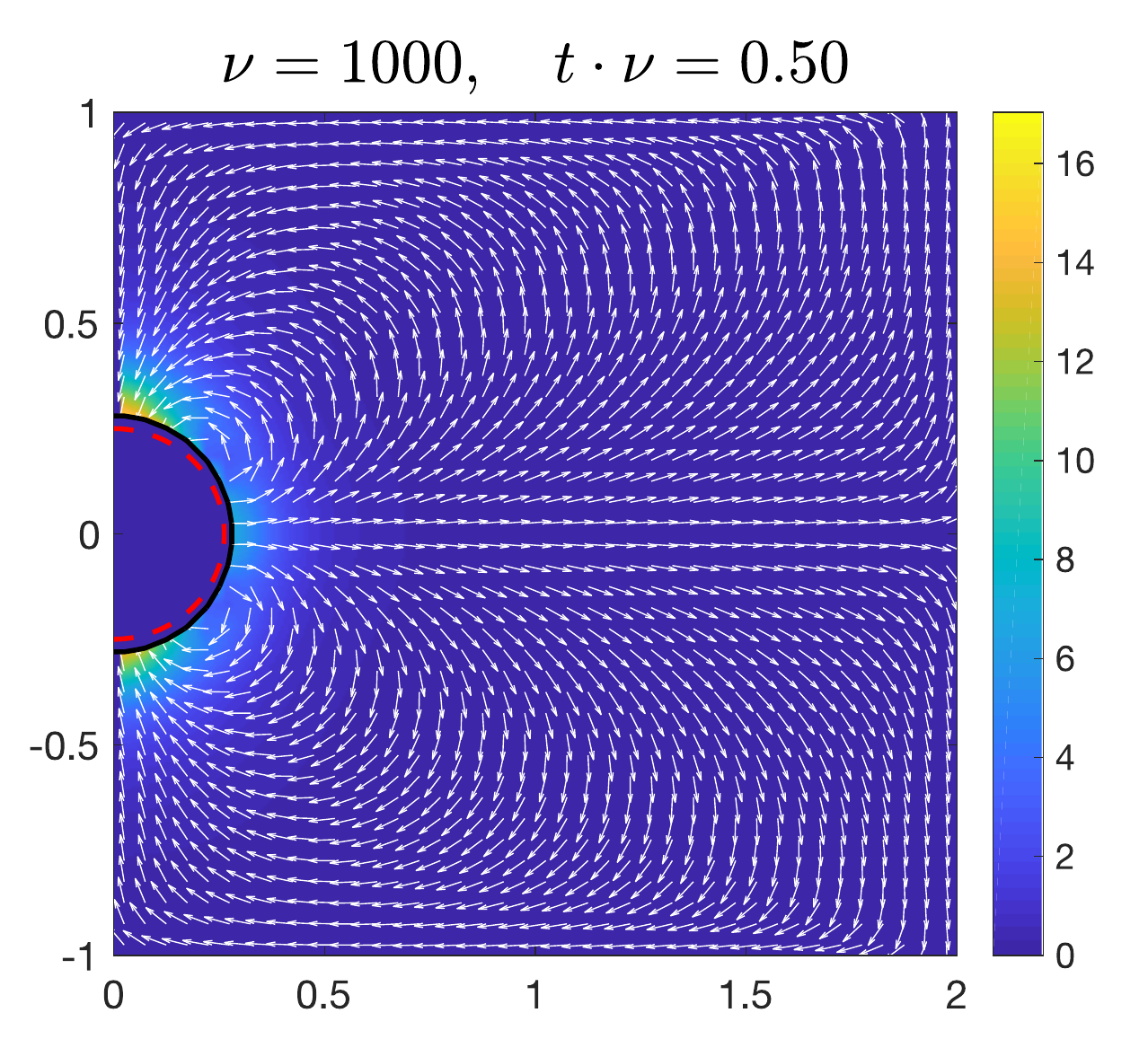}
\end{minipage}\hfill
\begin{minipage}[b]
	{.249\textwidth}
	\centering
	\includegraphics[width=\textwidth]{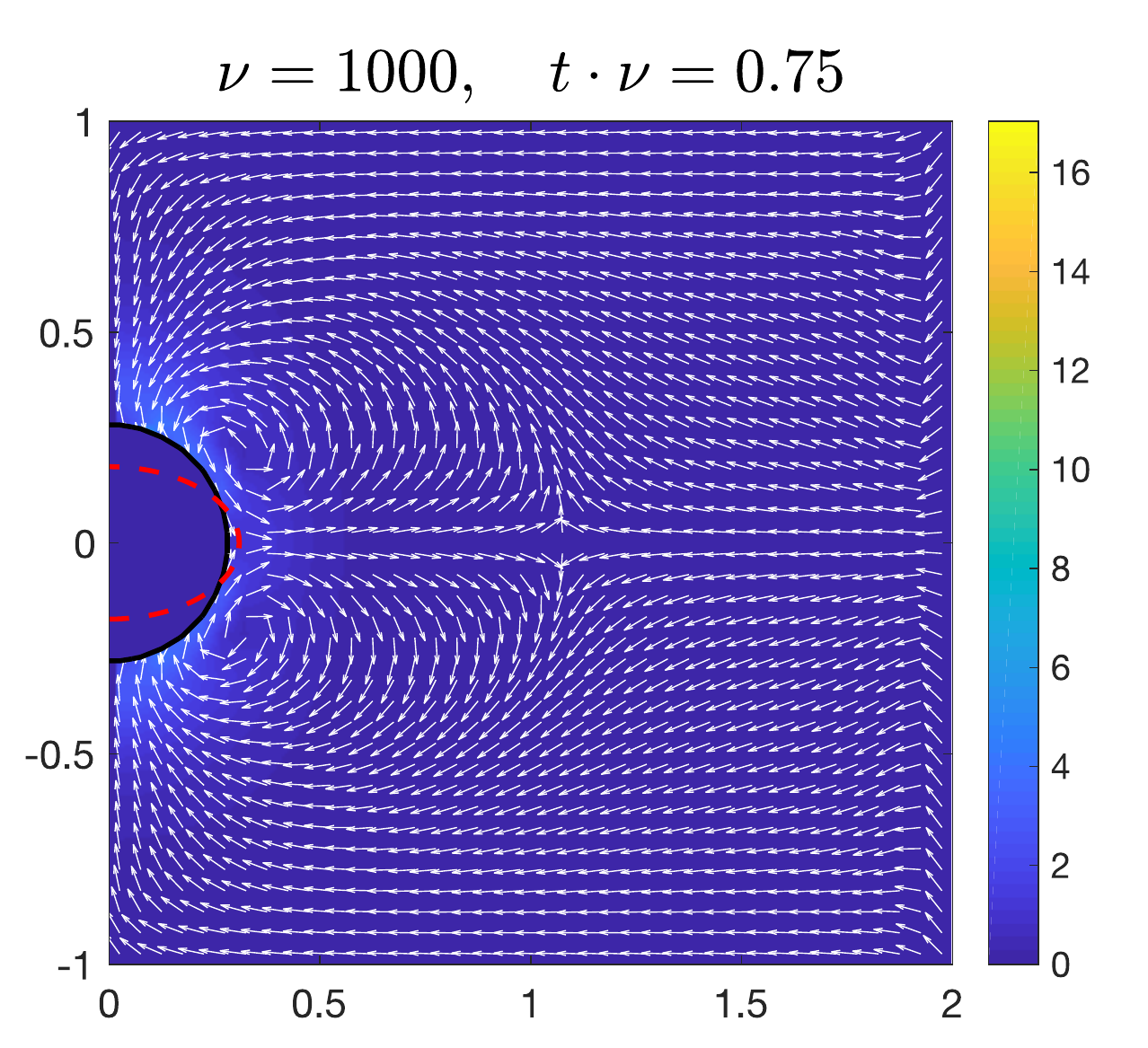}
\end{minipage}\hfill
\begin{minipage}[b]
	{.249\textwidth}
	\centering
	\includegraphics[width=\textwidth]{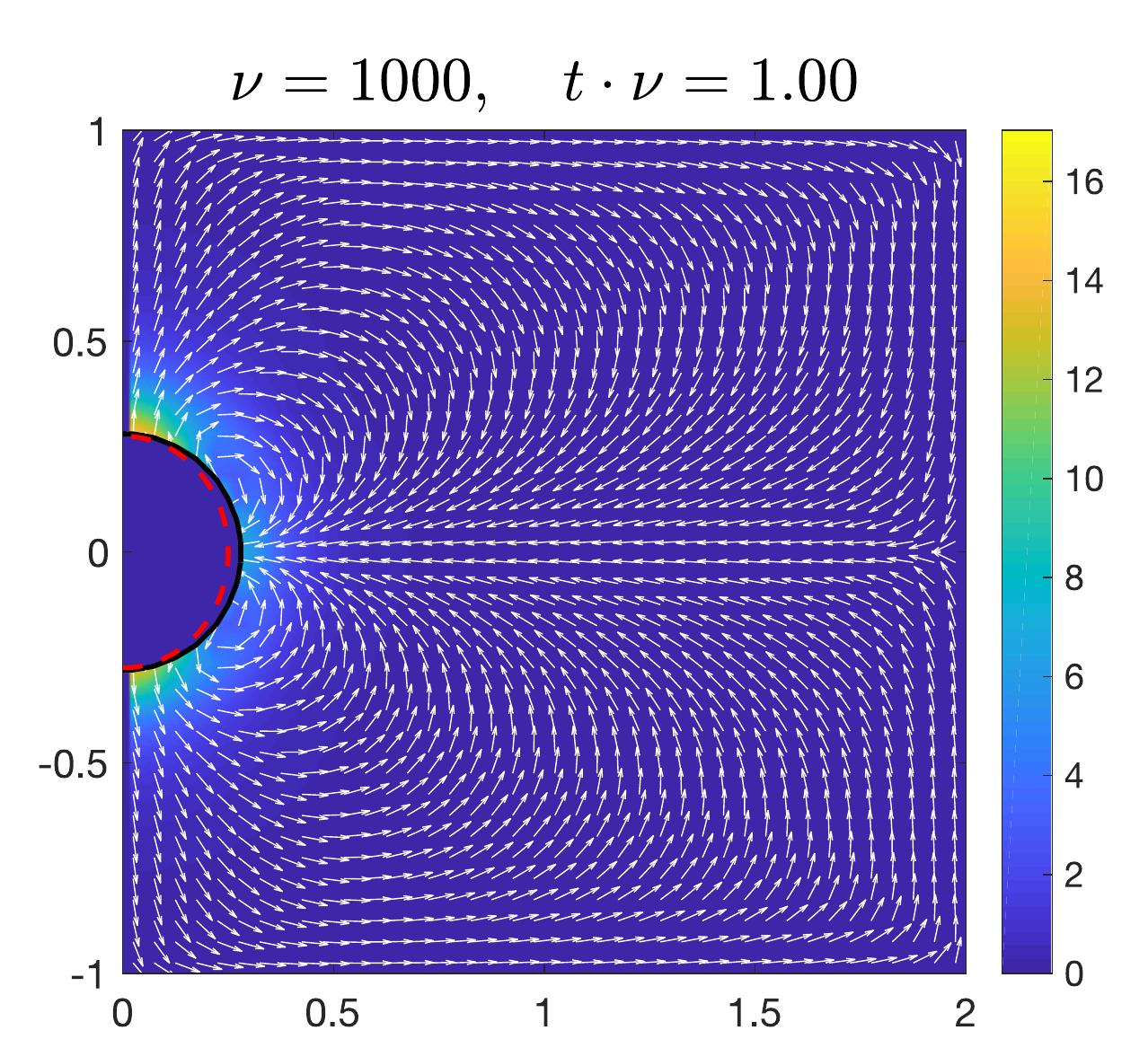}
\end{minipage}\hfill
\caption{Vector fields of the velocity at selected fractions of the first oscillation period $T$ ($\nu \cdot t = 0.25,0.50,0.75,1.00$) for the numerical tests of Sect.~\ref{sect:oscbubble}. Each row of plots corresponds to a numerical test. From top to bottom: \textsc{Test2a10}, \textsc{Test2a1000}, \textsc{Test2b10} and \textsc{Test2b1000}. The red dashed line is the fictitious representation of the bubble (where $A$ has been amplified by a factor 20 for graphical purposes).}
\label{fig:osc1and2}
\end{figure}

Following the same approach as in Sect.~\ref{sect:fixed}, we keep a steady computational bubble $\mathcal{B}(t)=\mathcal{B}(0)$ for $t>0$ and we model the velocity of the surface $\textbf{u}_b(\xi,z) = (u_b(\xi,z),v_b(\xi,z))$ at $(\xi,z)\colon \sqrt{\xi^2+z^2}=R_\mathcal{B}$ as:
\begin{align*}
u_b(\xi,z) &= \xi'_c(t)+\delta'_\xi(t) \cos(\theta) &\\
v_b(\xi,z) &= z'_c(t)+\delta'_z(t) \sin(\theta), \quad \text{ where } \quad  \theta = \arctan (z/\xi).&
\end{align*}	
In \textsc{Test2a}, we model an harmonic vertical oscillation of the spherical bubble:
\[
\xi(t)=0, \quad z(t)=A \sin(2 \pi \nu t), \quad \delta_\xi(t) = \delta_z(t) = R_\mathcal{B},
\]
while In \textsc{Test2b}, we model an ellipsoidal deformation of the bubble:
\begin{equation}\label{deftest2b}
\xi(t)=z(t)=0, \quad \delta_z(t) = R_\mathcal{B} (1+A \sin(2 \pi \nu t)), \quad \delta_\xi(t) = \sqrt{R_\mathcal{B}^3/\delta_z(t)}.
\end{equation}
In \eqref{deftest2b}, we observe that $\delta_\xi(t)$ has been defined to guarantee that the volume of the ellipsoid, $V(t)=4/3 (\pi \delta_\xi(t)^2 \delta_z(t))$, is constant over time.
We choose $R_\mathcal{B}=0.258$, $A = 0.01$ and $\nu=10$ (\textsc{Test2a10} and \textsc{Test2b10}) or $\nu=1000$ (\textsc{Test2a1000} and \textsc{Test2b1000}).
In Fig.~\ref{fig:osc1and2} we plot the vector fields of the fluid velocity at selected fractions of the first oscillation period $T$ ($\nu \cdot t = 0.25,0.50,0.75,1.00$). The colormap represents the magnitude of the velocity, while the red dashed line is the fictitious representation of the bubble (where $A$ has been amplified by a factor of 20 for graphical purposes).

In \textsc{Test2a10} we observe that a small vortex is generated next to the bubble, moving farther to the right after disappearing at around $\xi=1.5$ between $t \cdot \nu = 0.25$ and $t \cdot \nu = 0.50$. At that time, a new vortex is generated next to the bubble, disappearing between $t \cdot \nu = 0.75$ and $t \cdot \nu = 1.00$ and so on, approaching a periodic behaviour. A similar mechanism is observed in \textsc{Test2a1000}, except that the vortexes disappear when they are much closer to the bubble compared to \textsc{Test2a10}, say at around $\xi=0.5$. In \textsc{Test2b10} two vortexes are generated at the same time, moving towards the top right and bottom right corners of the domain, respectively. They disappear in favour of novel vortexes with the same timeline as in \textsc{Test2a10}. This phenomenon is observed in experimental results~\cite{tho2007cavitation}.
In \textsc{Test2b1000} a similar behaviour is observed, except that the two vortexes disappear when they are much closer to the bubble than in \textsc{Test2b10}.

In Table~\ref{table:divUaccuracy} we show that $\nabla \cdot \vec{u}_h$ decades with the same order of accuracy of the numerical error on $\vec{u}_h$, i.e.~$\nabla \cdot \vec{u}_h = \mathcal{O}(h^2)$, where $h$ is the spatial step.
For this purpose, we use \textsc{Test2a1000} and \textsc{Test2b1000}. Since we do not know the exact solution for these tests, we approximate the order of accuracy $q$ of $\vec{u}_h$ by using the Richardson extrapolation: 
\[
q \approx \log_2 \left( \frac{\left\| \vec{u}_h - \vec{u}_{h/2} \right\|_\infty}{\left\| \vec{u}_{h/2} - \vec{u}_{h/4} \right\|_\infty} \right).
\]
We observe numerically that $q \approx 2$.
The relative error on $\vec{u}_h$ is then approximated by
\[
{e}_h = \frac{\left\| \vec{u}_h - \vec{u}_{\text{exa}} \right\|_\infty}{\left\| \vec{u}_{\text{exa}} \right\|_\infty}
\approx 
\frac{4}{3} \frac{\left\| \vec{u}_h - \vec{u}_{h/2} \right\|_\infty}{\left\| \vec{u}_{h} \right\|_\infty}
\]
The relative error on $\nabla \cdot \vec{u}_h$ is computed by normalization with $\nabla \vec{u}_h$:
\[
{e}^{\text{div}}_h = \frac{\left\| \nabla \cdot \vec{u}_h \right\|_\infty}{\left\| \nabla \vec{u}_h \right\|_\infty}
\]
with $\left\| \nabla \vec{u}_h \right\|_\infty = \max \left\{ \left\| \, \left| \nabla u_h \right| \, \right\|_\infty, \left\| \, \left| \nabla v_h \right| \, \right\|_\infty \right\}$,
where $\left| \nabla u_h \right|$ and $\left| \nabla v_h \right|$ represent the central finite-difference approximations of $(u_x^2+u_y^2)^{1/2}$ and $(v_x^2+v_y^2)^{1/2}$, respectively.
The order of accuracy $q^\text{div}$ is approximated by
\[
q^\text{div} \approx
\log_2 \left( \frac{e_h^\text{div}}{e_{h/2}^\text{div}} \right) 
\]
We also observe from Table~\ref{table:divUaccuracy} that the relative error on the divergence ${e}^{\text{div}}_h$ is about one order of magnitude smaller than the relative error on velocity ${e}_h$.

  \begin{table}[H]
\centering      
\begin{tabular}{|| c | c | c | c | c | c ||} 
\hline \hline 
No.~of points & $\left\| \vec{u}_h - \vec{u}_{h/2} \right\|_\infty$ & order $q$ & $e_h$ & $e_h^\text{div}$ & order $q^\text{div}$\\ 
\hline 
32 $\times$ 32 & 2.12 $\cdot 10^{-1}$ & 1.79 & 1.62 $\cdot 10^{-3}$ & 2.02 $\cdot 10^{-4}$ & 1.74 \\ 
64 $\times$ 64 & 6.15 $\cdot 10^{-2}$ & 2.11 & 4.69 $\cdot 10^{-4}$ & 6.05 $\cdot 10^{-5}$ & 1.96 \\ 
128 $\times$ 128 & 1.42 $\cdot 10^{-2}$ & 2.02 & 1.08 $\cdot 10^{-4}$ & 1.55 $\cdot 10^{-5}$ & 1.99 \\ 
256 $\times$ 256 & 3.50 $\cdot 10^{-3}$ & - & 2.67 $\cdot 10^{-5}$ & 3.91 $\cdot 10^{-6}$ & 1.92 \\ 
512 $\times$ 512 & - & - & - & 1.03 $\cdot 10^{-6}$ & -\\ 
\hline \hline 
No.~of points & $\left\| \vec{u}_h - \vec{u}_{h/2} \right\|_\infty$ & order $q$ & $e_h$ & $e_h^\text{div}$ & order $q^\text{div}$\\ 
\hline 
32 $\times$ 32 & 4.51 $\cdot 10^{-1}$ & 1.88 & 3.44 $\cdot 10^{-3}$ & 6.12 $\cdot 10^{-4}$ & 1.79 \\ 
64 $\times$ 64 & 1.23 $\cdot 10^{-1}$ & 1.89 & 9.36 $\cdot 10^{-4}$ & 1.77 $\cdot 10^{-4}$ & 1.88 \\ 
128 $\times$ 128 & 3.32 $\cdot 10^{-2}$ & 1.94 & 2.53 $\cdot 10^{-4}$ & 4.84 $\cdot 10^{-5}$ & 1.84 \\ 
256 $\times$ 256 & 8.68 $\cdot 10^{-3}$ & - & 6.61 $\cdot 10^{-5}$ & 1.35 $\cdot 10^{-5}$ & 1.93 \\ 
512 $\times$ 512 & - & - & - & 3.55 $\cdot 10^{-6}$ & -\\ 
\hline \hline 
\end{tabular}
\caption{\textit{Relative numerical errors $e_h$ (on $\vec{u}_h$) and $e_h^\text{div}$ (on $\nabla \cdot \vec{u}_h$) and respective orders of accuracy $q$ and $q^\text{div}$ for \textsc{Test2a1000} (top) and \textsc{Test2b1000} (bottom). }}
    \label{table:divUaccuracy}  
 \end{table}

\begin{figure}[htp]
\centering
\hfill
\begin{minipage}[b]
	{.5\textwidth}
	\includegraphics[width=\textwidth]{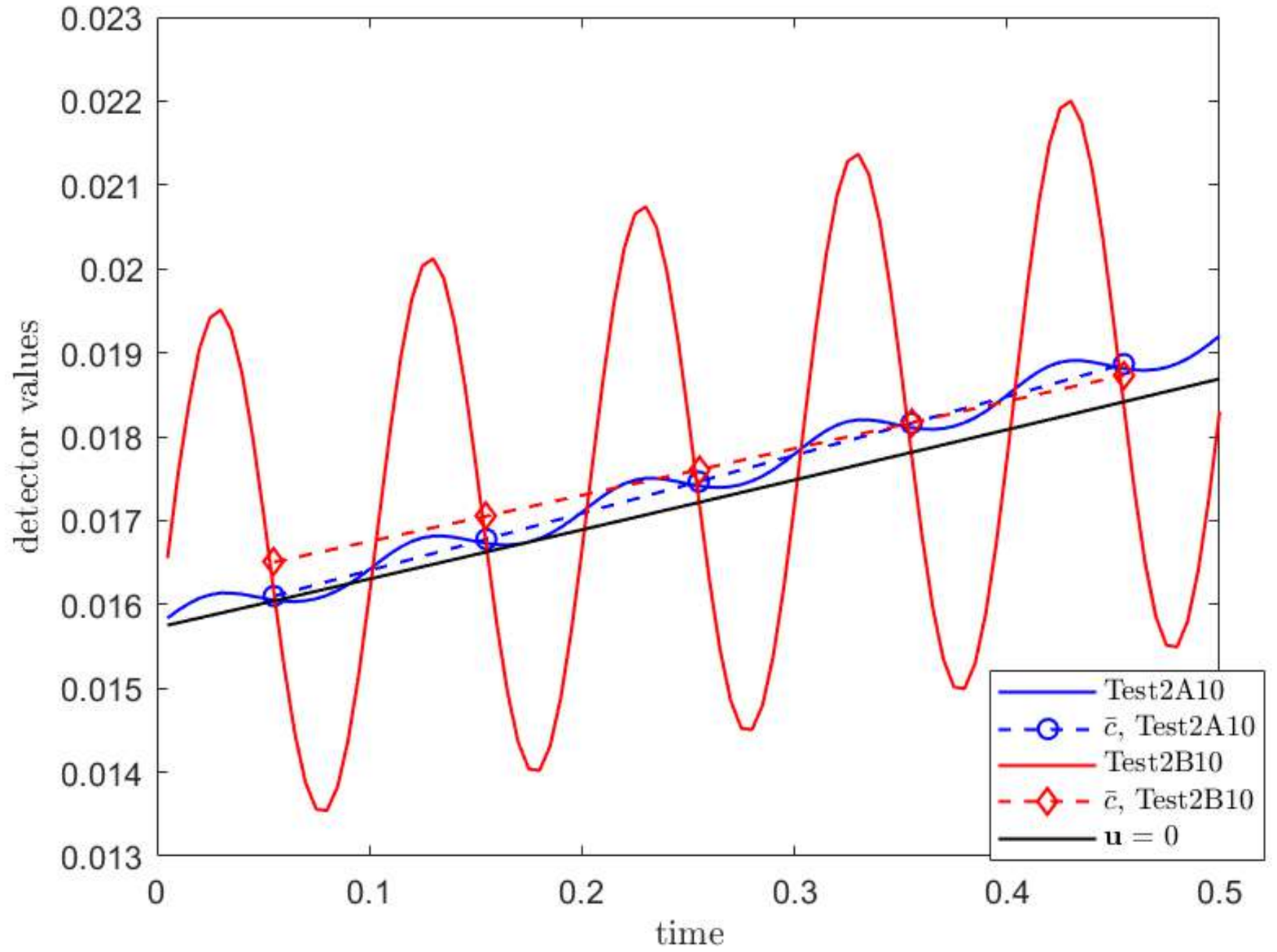}
\end{minipage}\hfill
\begin{minipage}[b]
	{.5\textwidth}
	\includegraphics[width=\textwidth]{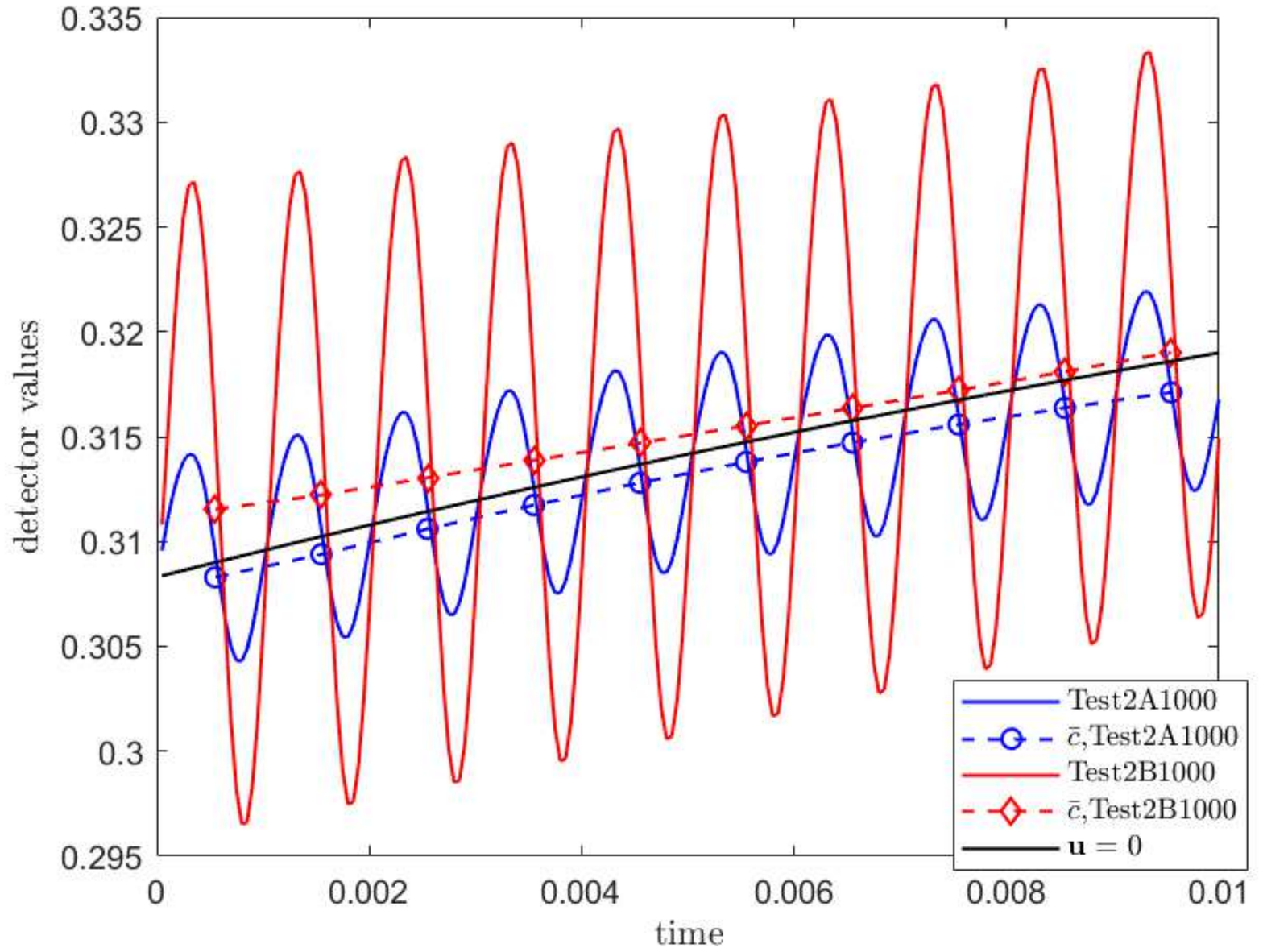}
\end{minipage}
\hspace*{\fill}
\caption{\textit{Detector values of the particle concentration $c$ of Eq.~\eqref{pde2dcoupled} at $(\xi_d=0.4,z_d=0)$. Initial condition is~\eqref{eq:IC}. The spatial step is $h = 1/120$. On the left we plot the comparison between \textsc{Test2a10} (blue line), \textsc{Test2b10} (red line) and steady-bubble case with $\textbf{u}=0$ (black line). Analogously, on the right, we show the comparison between \textsc{Test2a1000} and \textsc{Test2b1000}. The dashed lines represent the mean values of the respective tests.
}}
\label{fig:detector}
\end{figure}

In Fig.~\ref{fig:detector} we plot the particle concentration $c$ at a specific \textit{detector\/} point $(\xi_d=0.4,z_d=0)$ over time $t$ for different numerical tests (\textsc{Test2a10},\textsc{Test2a1000},\textsc{Test2b10},\textsc{Test2b1000}).
At the initial time $t=0$ the particles follow a Gaussian distribution centered in $(\xi_0=0.8,z_0=0)$, as follows (initial condition):
\begin{equation}\label{eq:IC} 
c(\xi,z,0) = a_1 e^{-a_2((\xi-\xi_0)^2+(z-z_0)^2)}, \quad a_1 = \left(2\pi\sigma^2\right)^{-3/2}, \quad a_2 = \left(2\sigma^2\right)^{-1}, \sigma = 0.4.
\end{equation}
The black line refers to the case of steady bubble (${\bf u}=0$). All tests present an oscillating behaviour of the particle concentration in the vicinity of the black line. The oscillation frequency of the particle concentration is strictly related with the oscillation frequency of the bubble, while the amplitude depends on the type of bubble oscillation (harmonic or ellipsoidal). However, the temporal average of the particle concentration of the proposed tests, namely
$\overline{c}(x,t) = \frac{1}{T}\int_{t-T/2}^{t+T/2} c(x,\tau) \, d \tau = \nu \int_{t-1/(2\nu)}^{t+1/(2\nu)} c(x,\tau) \, d \tau$, 
does not seem to match the green line, suggesting that the bubble oscillation actually changes the particle distribution in the vicinity of the bubble. The problem presents a temporal multiscale effect and a rigorous mathematical explanation of this phenomenon should provide the associated PDEs for $\overline{c}$, showing that these equations differ from the simple diffusion equations with ${\bf u}=0$.
This is part of our ongoing effort.

\section{Conclusions}
\label{sec:conclusions}

We have presented a second order accurate numerical method for the recently developed multiscale model of sorption kinetic \cite{multiscale_mod}, in the single carrier approximation. The problem consists of a concentration of particles that diffuse in a fluid agitated by an oscillating bubble. In addition, the bubble attracts the particles that are in the vicinity of its surface. This attraction is modelled by a time-dependent boundary condition on the bubble surface, where normal and (second order) tangential derivatives are involved. The region occupied by the bubble is implicitly described by a level-set function and the boundary conditions on the curved boundary (bubble surface) are discretized by a proper ghost-point technique. A multigrid method is designed to efficiently solve the linear system arising from the discretization of the equations.
The complexity of the boundary condition on the bubble surface leads to a specific stability condition that must be satisfied by the relaxation scheme of the multigrid method. 

The fluid motion is modelled by the Stokes equations, solved by a monolithic approach (continuity and momentum equations are solved simultaneously). A simplified model is provided for the treatment of the boundary conditions in the case of moving bubble, in which the bubble is computationally steady and its oscillations are modeled by a time-dependent fluid velocity imposed on its surface. This simplification is justified by the small amplitude of the bubble oscillations, and a comparison with the more realistic moving domain model confirms that the differences between the two approaches are negligible for the problems investigated in this paper. Furthermore, this approximation makes the 
computation more efficient, because it avoids evolving the domain and its discretization, so the coefficients of the linear system \eqref{CNdiscspace} are time independent.

Two types of bubble oscillations are implemented: harmonic and ellipsoidal oscillation, the latter providing a better representation of the experimental results existing in literature.

The particle concentration at a specific point of the domain is reconstructed over time for the two types of oscillation. The same test is repeated for a steady fluid case (actual steady bubble). We observed that the particle concentration for the oscillating bubble oscillates around an average function that is close to (but not the same as) the one obtained for the actual steady bubble. The discrepancy between the two functions can be mathematically described by numerical approaches developed for temporal multiscale problems.
A more rigorous analysis of this phenomenon is part of our ongoing effort. Future works will include a saturation effect (already modelled in~\cite{multiscale_mod}, when high concentrations is reached near the surface of the bubble) and a full two carrier model that describes the interaction between the two species of ions (with a potential obtained from a self-consistent Poisson equation).

\section*{Acknowledgments}
G.R.~and C.A.~acknowledge partial support from ITN-ETN Horizon 2020 Project ModCompShock, Modeling and Computation on Shocks and Interface (Project Reference 642768), and from PRIN Project 2017 entitled ''Innovative numerical methods for evolutionary partial differential equations and applications'' (No.2017KKJP4X), coming from the Italian Ministry of Education, University and Research (MIUR).
All authors acknowledge support from GNCS--INDAM (National Group for Scientific Computing, Italy).

\bibliographystyle{unsrt}
\bibliography{main}

\end{document}